\documentclass[a4paper,12pt]{article}
\usepackage{amsmath}
\usepackage{amscd}
\usepackage[xdvi]{graphicx}
\usepackage[dvips]{color}
\usepackage{times}
\usepackage{txfonts}

\newcommand{\bm}[1]{\mbox{\boldmath $#1$}}

\makeatletter

\@addtoreset{equation}{section}

\begin{document}
\setlength{\baselineskip}{24pt}
\begin{center}
\textbf{\LARGE{A proof of the Kuramoto conjecture for a bifurcation structure of the infinite dimensional Kuramoto model}}
\end{center}

\setlength{\baselineskip}{14pt}

\begin{center}
Faculty of Mathematics, Kyushu University, Fukuoka,
819-0395, Japan

\large{Hayato CHIBA} \footnote{E mail address : chiba@imi.kyushu-u.ac.jp}
\end{center}
\begin{center}
Revised Oct 3, 2012
\end{center}

\begin{center}
\textbf{Abstract}
\end{center}

The Kuramoto model is a system of ordinary differential equations for describing synchronization phenomena
defined as a coupled phase oscillators.
In this paper, a bifurcation structure of the infinite dimensional Kuramoto model is investigated.
A purpose here is to prove the bifurcation diagram of the model conjectured by Kuramoto in 1984;
if the coupling strength $K$ between oscillators, which is a parameter of the system,
is smaller than some threshold $K_c$, the de-synchronous state (trivial steady state) 
is asymptotically stable, while if $K$ exceeds $K_c$, a nontrivial stable solution,
which corresponds to the synchronization, bifurcates from the de-synchronous state.
One of the difficulties to prove the conjecture is that a certain non-selfadjoint linear operator, 
which defines a linear part of the Kuramoto model, has the continuous spectrum on the imaginary axis.
Hence, the standard spectral theory is not applicable to prove a bifurcation as well as the asymptotic stability
of the steady state.
In this paper, the spectral theory on a space of generalized functions is developed with the aid of a rigged Hilbert space
to avoid the continuous spectrum on the imaginary axis.
Although the linear operator has an unbounded continuous spectrum on a Hilbert space,
it is shown that it admits a spectral decomposition consisting of a countable number of eigenfunctions
on a space of generalized functions.
The semigroup generated by the linear operator will be estimated with the aid of the spectral theory
on a rigged Hilbert space to prove the linear stability of the steady state of the system.
The center manifold theory is also developed on a space of generalized functions.
It is proved that there exists a finite dimensional center manifold on a space of generalized functions,
while a center manifold on a Hilbert space is of infinite dimensional because of the continuous spectrum on the 
imaginary axis.
These results are applied to the stability and bifurcation theory of the Kuramoto model
to obtain a bifurcation diagram conjectured by Kuramoto.
\\[0.2cm]
\textbf{Keywords}: infinite dimensional dynamical systems; center manifold theory; continuous spectrum;
spectral theory; rigged Hilbert space; coupled oscillators; Kuramoto model

\newpage
\tableofcontents


\newpage

\section{Introduction}

Collective synchronization phenomena are observed in a variety of areas such as chemical reactions,
engineering circuits and biological populations~\cite{Pik}.
In order to investigate such phenomena, Kuramoto~\cite{Kura1} proposed the system of ordinary differential equations
\begin{equation}
\frac{d\theta _i}{dt} 
= \omega _i + \frac{K}{N} \sum^N_{j=1} \sin (\theta _j - \theta _i),\,\, i= 1, \cdots  ,N,
\label{KMN}
\end{equation}
where $\theta _i = \theta _i(t) \in [ 0, 2\pi )$ is a dependent variable which denotes the phase of an $i$-th oscillator on a circle,
$\omega _i\in \mathbf{R}$ denotes its natural frequency, $K>0$ is a coupling strength,
and where $N$ is the number of oscillators.
Eq.(\ref{KMN}) is derived by means of the averaging method from coupled dynamical systems having 
limit cycles, and now it is called the \textit{Kuramoto model}.

It is obvious that when $K=0$, $\theta _i(t)$ and $\theta _j(t)$ rotate on a circle at 
different velocities unless $\omega _i$ is equal to $\omega _j$, 
and this fact is true for sufficiently small $K>0$.
On the other hand, if $K$ is sufficiently large, it is numerically observed that
some of oscillators or all of them tend to rotate at the same velocity on average, which is called the 
\textit{synchronization}~\cite{Pik,Str1}.
If $N$ is small, such a transition from de-synchronization to synchronization may be well revealed
by means of the bifurcation theory~\cite{ChiPa,Mai1,Mai2}.
However, if $N$ is large, it is difficult to investigate the transition from the view point of
the bifurcation theory and it is still far from understood.

\begin{figure}
\begin{center}
\includegraphics[]{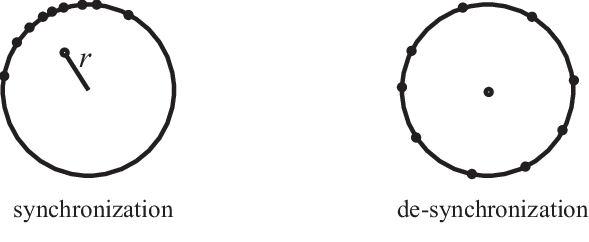}
\caption[]{The order parameter of the Kuramoto model.}
\label{fig1}
\end{center}
\end{figure}

In order to evaluate whether synchronization occurs or not, Kuramoto introduced
the \textit{order parameter} $r(t)e^{\sqrt{-1}\psi (t)}$ by
\begin{equation}
r(t)e^{\sqrt{-1}\psi (t)} := \frac{1}{N}\sum^N_{j=1} e^{\sqrt{-1} \theta _j(t)},
\label{order2}
\end{equation}
where $r, \psi \in \mathbf{R}$.
The order parameter gives the centroid of oscillators.
It seems that if synchronous state is formed, $r(t)$ takes a positive number, while
if de-synchronization is stable, $r(t)$ is zero on time average (see Fig.\ref{fig1}).
Further, this is true for every $t$ when $N$ is sufficiently large so that a statistical-mechanical 
description is applied.
Based on this observation and some formal calculations, Kuramoto conjectured a bifurcation diagram
of $r(t)$ as follows:
\\[0.2cm]
\textbf{Kuramoto conjecture}

\begin{figure}
\begin{center}
\includegraphics[]{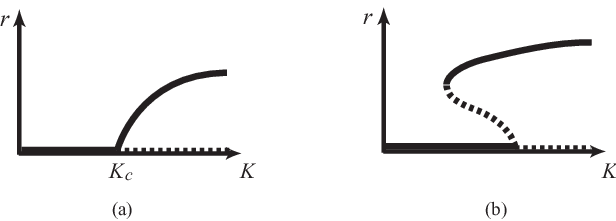}
\caption[]{Typical bifurcation diagrams of the order parameter for the cases that (a) $g(\omega )$
is even and unimodal (b) $g(\omega )$ is even and bimodal.
Solid lines denote stable solutions and dotted lines denote unstable solutions.}
\label{fig2}
\end{center}
\end{figure}

Suppose that $N\to \infty$ and natural frequencies $\omega _i$'s are distributed according to a probability
density function $g(\omega )$.
If $g(\omega )$ is an even and unimodal function such that $g''(0)\neq 0$, then the bifurcation
diagram of $r(t)$ is given as Fig.\ref{fig2} (a); that is, if the coupling strength $K$ is 
smaller than $K_c := 2/(\pi g(0))$, then $r(t) \equiv 0$ is asymptotically stable.
On the other hand, if $K$ is larger than $K_c$, the synchronous state emerges; there exists a positive constant $r_c$ such that 
$r(t) = r_c$ is asymptotically stable.
Near the transition point $K_c$, $r_c$ is of order $O((K- K_c)^{1/2})$.

A function $g(\omega )$ is called unimodal (at $\omega =0$) if $g(\omega _1) > g(\omega _2)$ for $0\leq \omega _1 < \omega _2$
and $g(\omega _1) < g(\omega _2)$ for $\omega _1 < \omega _2 \leq 0$.
Now the value $K_c = 2/(\pi g(0))$ is called the \textit{Kuramoto transition point}.
See \cite{Kura2} and \cite{Str1} for Kuramoto's discussion.
\\

In the present paper, the Kuramoto conjecture will be proved in the following sense:
At first, we will define the continuous limit of the model in Sec.2 to express the dynamics of the infinite 
number of oscillators ($N\to \infty$).
The trivial steady state of the continuous model corresponds to the de-synchronous state $r\equiv 0$.
For the continuous model, the following theorems will be proved.
\\[0.2cm]
\textbf{Theorem 1.1 (instability of the trivial state).}
Suppose that $g(\omega )$ is even, unimodal and continuous.
When $K > K_c:=2/(\pi g(0))$, then the trivial steady state of the continuous model is linearly unstable.
\\

This linear instability result was essentially obtained by Strogatz and Mirollo \cite{Str3}.
Although we do not give a proof of a local nonlinear instability, it is proved in the same way as the local nonlinear
stability result below.
\\[0.2cm]
\textbf{Theorem 1.2 (local stability of the trivial state).}
Suppose that $g (\omega )$ is the Gaussian distribution or a rational function which is even, unimodal
and bounded on $\mathbf{R}$.
When $0<K<K_c$,
there exists a positive constant $\delta$ such that if the initial condition $h(\theta )$ for the continuous model
(\ref{conti}) satisfies
\begin{equation}
\left| \int^{2\pi}_{0} \! e^{j\sqrt{-1}\theta } h(\theta ) d\theta  \right| \leq  \delta,
 \quad j=1,2,\cdots ,
\end{equation}
then the continuous limit $\eta (t)$ of the order parameter defined in (\ref{conti}) decays to zero 
exponentially as $t\to \infty$.
\\

This stability result will be stated as Thm.6.1 in more detail:
under the above assumptions, the trivial state of the continuous model proves to be locally stable 
with respect to a topology of a certain topological vector space constructing a rigged Hilbert space.
Thm.1.2 is obtained as a corollary of Thm.6.1.
\\[0.2cm]
\textbf{Theorem 1.3 (bifurcation).}
Suppose that $g (\omega )$ is the Gaussian distribution or a rational function
which is even, unimodal and bounded on $\mathbf{R}$.
For the continuous model,
there exist positive constants $\varepsilon _0$ and $\delta $ such that
if $K_c < K < K_c + \varepsilon _0$ and if the initial condition $h(\theta )$ satisfies
\begin{equation}
\left| \int^{2\pi}_{0} \! e^{\sqrt{-1}j \theta }h(\theta ) d\theta  \right| < \delta,
\quad j = 1,2,\cdots ,
\end{equation}
then the continuous limit $\eta (t)$ of the order parameter tends to the constant expressed as
\begin{equation}
r(t) = |\eta (t) |= \sqrt{\frac{-16}{\pi K_c^4 g''(0)}}\sqrt{K - K_c} + O(K-K_c), 
\end{equation}
as $t\to \infty$. In particular, the bifurcation diagram of the order parameter is given as Fig.\ref{fig2} (a).
\\

This result will be proved in Thm.7.10 with the aid of the center manifold theory on a rigged Hilbert space.
Again, a bifurcation of a stable nontrivial solution of the continuous model will be proved
with respect to a topology of a certain topological vector space.

A few remarks are in order.
\\
$\bullet$ Our bifurcation theory is applicable to a certain class of distribution functions $g(\omega )$.
It will turn out that one of the most essential assumptions is the holomorphy (meromorphy) of $g(\omega )$.
For example, let us slightly deform the Gaussian $g(\omega )$ so that it sags in the center
as it becomes bimodal function.
In this case, since $g''(0) >0$, $|\eta (t)|$ above is positive when $K<K_c$.
This means that a subcritical bifurcation occurs and the bifurcation diagram shown in Fig.\ref{fig2} (b) is obtained
at least near the bifurcation point $K=K_c$.
\\
$\bullet$ It is proved in \cite{Chi4} that the order parameter (\ref{order2}) for the $N$-dimensional Kuramoto 
model converges to that of the continuous model (\ref{conti}) as $N\to \infty$ in a certain probabilistic sense
for \textit{each} $t>0$.
\\
$\bullet$ In \cite{ChiNi}, bifurcation diagrams of the Kuramoto-Daido model (i.e. a coupling function
includes higher harmonic terms such as $\sin 2(\theta _j - \theta _i)$) are obtained in the same way 
as the present paper, although the existence of center manifolds has not been proved for the Kuramoto-Daido model.
\\
$\bullet$ In this paper, only local stability is proved and global one is still open.
\\

In the rest of this section, known results for the Kuramoto conjecture will be briefly reviewed
and our idea to prove the above theorems are explained.
See Strogatz \cite{Str1} for history of the Kuramoto conjecture.

In the last two decades, many studies to confirm the Kuramoto conjecture have been done.
Significant papers of Strogatz and coauthors~\cite{Str3, Str2} investigated the linear stability of the trivial solution,
which corresponds to the de-synchronous state $r\equiv 0$.
In \cite{Str3}, they introduced the continuous model for the Kuramoto model to describe the situation $N\to \infty$.
They derived the Kuramoto transition point $K_c = 2/(\pi g(0))$ and showed that if $K > K_c$,
the de-synchronous state is unstable because of eigenvalues on the right half plane.
On the other hand, when $0<K \leq K_c$, a linear operator $T_1$, which defines the linearized equation
of the continuous model around the de-synchronous state, has no eigenvalues;
the spectrum of $T_1$ consists only of the continuous spectrum on the imaginary axis.
This implies that the standard stability theory of dynamical systems is not applicable to this system.
However, in \cite{Str2}, they found that
an analytic continuation of the resolvent $(\lambda -T_1)^{-1}$ may have poles (\textit{resonance poles})
on the left half plane for a wide class of distribution functions $g(\omega )$.
They remarked a possibility that resonance poles induce a decay of the order parameter $r$ by a linear analysis.
This claim will be rigorously proved in this paper for a certain 
class of distribution functions by taking into account nonlinear terms (Thm.1.2).
In \cite{Mir2}, the spectra of linearized systems around other steady states, 
which correspond to solutions with positive $r=r_c$, are investigated.
They found that linear operators, which is obtained from the linearization of the system around synchronous states,
have continuous spectra on the imaginary axis.
Nevertheless, they again remarked that such solutions can be asymptotically stable because of the resonance poles.

Since results of Strogatz \textit{et al}. are based on a linearized analysis, effects of nonlinear terms are neglected.
To investigate nonlinear dynamics, the bifurcation theory is often used.
However, investigating the bifurcation structure near the transition point $K_c$ involves further difficult problems
because the operator $T_1$ has a continuous spectrum on the imaginary axis, that is, 
a center manifold in a usual sense is of infinite dimensional.
To avoid this difficulty, Bonilla \textit{et al.}~\cite{Ace2, Bon1, Bon2} and Crawford \textit{et al.}~\cite{Cra1, Cra2, Cra3}
added a perturbation (noise) with the strength $D>0$ to the Kuramoto model.
Then, the continuous spectrum moves to the left side by $D$, and thus the usual center manifold
reduction is applicable.
When $g(\omega )$ is an even and unimodal function, they obtained the Kuramoto bifurcation
diagram (Fig.2 (a)), however, obviously their methods are not valid when $D=0$.
For example, in Crawford's method, an eigenfunction of $T_1$ associated with a center subspace diverges as $D\to 0$
because an eigenvalue on the imaginary axis is embedded in the continuous spectrum as $D\to 0$.
Thus the original Kuramoto conjecture was still open.

Despite the active interest in the case that the distribution function $g(\omega )$ is even and unimodal,
bifurcation diagrams of $r$ for $g(\omega )$ other than the even and unimodal case are not understood well.
Martens \textit{et al}.~\cite{Mar} investigated the bifurcation diagram for a bimodal $g(\omega )$
which consists of two Lorentzian distributions.
In particular, they found that stable synchronous states can coexist with stable 
de-synchronous states if $K$ is slightly smaller than $K_c$ (see Fig.\ref{fig2} (b)).
Their analysis depends on extensive symmetries of the Kuramoto model found by Ott and Antonsen~\cite{Ott1, Ott2} (see also \cite{Marv})
and on the special form of $g(\omega )$, however,
such a diagram seems to be common for any bimodal distributions.

In this paper, the stability, spectral and bifurcation theory of the continuous model
of the Kuramoto model will be developed to prove the Kuramoto conjecture.
Let $T_1$ be a linear operator obtained by linearizing the continuous model (\ref{conti}) 
around the de-synchronous state.
The spectrum and the semigroup of $T_1$ will be investigated in detail.
The operator $T_1 = T_1(K)$ defined on the weighted Lebesgue space $L^2(\mathbf{R}, g(\omega )d\omega )$
has the continuous spectrum $\sigma _c(T_1)$ on the imaginary axis for any $K>0$.
For example, when $g$ is the Gaussian distribution, then $\sigma _c(T_1) = \sqrt{-1}\mathbf{R}$.
At first, we derive the transition point (bifurcation point) $K_c$ for any distribution function $g(\omega )$;
When $K>K_c$, $T_1$ has eigenvalues $\lambda = \lambda (K)$ on the right half plane.
As $K$ decreases, $\lambda (K)$ goes to the left side, and at $K=K_c$, the eigenvalues are absorbed into the 
continuous spectrum on the imaginary axis and disappear.
When $0<K<K_c$, there are no eigenvalues and the spectrum of $T_1$ consists of the continuous spectrum.
As a corollary, the Kuramoto transition point $K_c = 2/(\pi g(0))$ is obtained if $g(\omega )$ 
is an even and unimodal function.
When $K>K_c$, it is proved that the de-synchronous state is unstable because the operator 
$T_1$ has eigenvalues on the right half plane.

On the other hand, when $0<K\leq K_c$, the operator $T_1$ has no eigenvalues and the continuous spectrum
lies on the imaginary axis.
Thus the stability of the de-synchronous state is nontrivial.
Despite this fact, under appropriate assumptions for $g(\omega )$,
the order parameter proves to decay exponentially to zero as $t\to \infty$ 
because of the existence of resonance poles on the left half plane, as was expected by Strogatz \textit{et al}.~\cite{Str2}.
To prove it, the notion of spectrum is extended.
Roughly speaking, the spectrum is the set of singularities of the resolvent $(\lambda -T_1)^{-1}$.
However, if $g(\omega )$ has an analytic continuation, the resolvent has 
an analytic continuation if the domain is restricted to a suitable function space.
The analytic continuation has singularities, which are called resonance poles, on the second Riemann sheet.
By using the Laplace inversion formula for a semigroup, we will prove that the resonance poles induce
an exponential decay of the order parameter.
This suggests that in general, linear stability of a trivial solution of a linear equation on an
infinite dimensional space is determined by not only the spectrum of the linear operator but 
also its resonance poles.

Next purpose is to investigate a bifurcation at $K=K_c$.
To handle the continuous spectrum on the imaginary axis,
a spectral theory of the resonance poles is developed with the aid of a rigged Hilbert space
(Gelfand triplet).
A rigged Hilbert space consists of three topological vector spaces 
\begin{eqnarray*}
X \subset H \subset X',
\end{eqnarray*}
a space $X$ of test functions,
a Hilbert space $H$ (in our problem, this is the weighted Lebesgue space $L^2(\mathbf{R}, g(\omega )d\omega )$)
and the dual space $X'$ of $X$ (a space of continuous linear functionals on $X$ called generalized functions).
A suitable choice of $X$ depends on $g(\omega )$.
In this paper, two cases are considered: (i) $g(\omega )$ is the Gaussian distribution, (ii) $g(\omega )$
is a rational function (e.g. Lorentzian distribution $g(\omega ) = 1/(\pi (1 + \omega ^2))$).
For the case (i), $X := \mathrm{Exp}_+$ is a space of holomorphic functions $\phi (z)$ defined near the real axis and 
the upper half plane such that $\sup_{\mathrm{Im}(z) \geq -\varepsilon } |\phi (z)|e^{-\beta |z|}$ 
is finite for some $\varepsilon >0$ and $\beta \geq 0$.
For the case (ii), $X := H_+$ is a space of bounded holomorphic functions on the real axis and the upper half plane.
For both cases, we will show that if the domain of the resolvent $(\lambda -T_1)^{-1}$ is restricted to $X$,
then it has an $X'$-valued meromorphic continuation from the right half plane to the left half plane
beyond the continuous spectrum on the imaginary axis.
Although $(\lambda -T_1)^{-1}$ diverges on the imaginary axis as an operator on $H$ because of the continuous spectrum,
it has an analytic continuation from the right to the left as an operator from $X$ into $X'$.
Singularities of the continuation of the resolvent is called resonance poles $\lambda _n \, (n=0,1,\cdots )$.
We will show that there exists a generalized function $\mu_n \in X'$ satisfying
\begin{eqnarray*}
T_1^\times \mu_n = \lambda _n \mu_n,
\end{eqnarray*}
where $T_1^\times : X' \to X'$ is a dual operator of $T_1$ and $\mu_n$ is called 
the \textit{generalized eigenfunction} associated with the resonance pole.
Despite the fact that $T_1$ is \textit{not} a selfadjoint operator and it has the continuous spectrum,
it is proved that the operator $T_1$ admits the spectral decomposition
on $X'$ consisting of a countable number of generalized eigenfunctions:
roughly speaking, any element $\phi$ in $X$  is decomposed as 
\begin{eqnarray*}
\phi = \sum^\infty_{n=0}\mu_n(\phi)\cdot \mu_n.
\end{eqnarray*}
Further, it is shown that for the case (ii), the decomposition is reduced to a finite sum because of a certain 
degeneracy of the space $X = H_+$.
We further investigate the semigroup generated by $T_1$ and the projection to the eigenspace of $\mu_n$.
It is proved that the semigroup $e^{T_1t}$ behaves as 
\begin{eqnarray*}
e^{T_1t}\phi = \sum^\infty_{n=0}e^{\lambda _nt}\mu_n(\phi) \cdot \mu_n
\end{eqnarray*}
for any $\phi \in X$.
This equality completely determines the dynamics of the linearized Kuramoto model.
In particular, when $0 < K< K_c$, all resonance poles lie on the left half plane: $\mathrm{Re}(\lambda _n)<0$, which proves the
linear stability of the de-synchronous state.
When $K=K_c$, there are resonance poles on the imaginary axis.
We define a generalized center subspace $\mathbf{E}_c$ on $X'$ to be a space spanned by generalized eigenfunctions
associated with resonance poles on the imaginary axis.
It is remarkable that though the center subspace in a usual sense is of infinite dimensional
because of the continuous spectrum on the imaginary axis, the dimension of the generalized center subspace on $X'$
is finite in general.
The projection operator to the generalized center subspace will be investigated in detail.

Note that the spectral decomposition based on a rigged Hilbert space was originally proposed by Gelfand \textit{et al.}
~\cite{Gel1, Mau}.
They proposed a spectral decomposition of a selfadjoint operator by using a system of generalized eigenfunctions, 
however, it involves an integral; that is, eigenfunctions are uncountable.
Our results are quite different from Gelfand's one in that our operator $T_1$ is not selfadjoint and its spectral 
decomposition consists of a countable number of eigenfunctions.

Finally, we apply the center manifold reduction to the continuous Kuramoto model by 
regarding it as an evolution equation on $X'$.
Since the generalized center subspace is of finite dimensional, a corresponding center manifold on $X'$
seems to be a finite dimensional manifold.
However, there are no existence theorems of center manifolds on $X'$ because $X'$ is \textit{not}
a Banach space.
To prove the existence of a center manifold, we introduce a topology on $X$ in a technical way so that
the dual space $X'$ becomes a complete metric space.
With this topology, $X'$ becomes a topological vector space called Montel space, which is obtained as a projective limit
of Banach spaces. 
This topology has a very convenient property that every weakly convergent series
in $X'$ is also convergent with respect the metric.
By using this topology and the spectral decomposition, the existence of a finite dimensional center manifold for
the Kuramoto model will be proved.
The dynamics on the center manifold will be derived when $g(\omega )$ is Gaussian.
In this case, the center manifold on $X'$ is of one dimensional,
and we can show that the synchronous solution (a solution such that $r>0$) emerges through the pitchfork bifurcation,
which proves Thm.1.3.

This paper is organized as follows:
In Sec.2, the continuous model for the Kuramoto model is defined and its basic properties are reviewed.
In Sec.3, Kuramoto's transition point $K_c$ is derived and it is proved that
if $K>K_c$, the de-synchronous state is unstable because of eigenvalues on the right half plane.
In Sec.4, the linear stability of the de-synchronous state is investigated.
We will show that when $0< K < K_c$, the order parameter decays exponentially to zero as $t\to \infty$
because of the existence of resonance poles.
In Sec.5, the spectral theory of resonance poles on a rigged Hilbert space is developed.
We investigate properties of the operator $T_1$, the semigroup, eigenfunctions, projections by means of the rigged Hilbert space.
In Sec.6, the nonlinear stability of the de-synchronous state is proved as an application 
of the spectral decomposition on the rigged Hilbert space.
It is shown that when $0< K < K_c$, the order parameter tends to zero as $t\to \infty$ without neglecting the nonlinear term.
The center manifold theory will be developed in Sec.7.
Sec.7.1 to Sec.7.4 are devoted to the proof of the existence of a center manifold on the dual space $X'$.
In Sec.7.5, the dynamics on the center manifold is derived, and the Kuramoto conjecture is solved.


\section{Continuous model}

In this section, we define a continuous model of the Kuramoto model
and show a few properties of it.

For the $N$-dimensional Kuramoto model (\ref{KMN}), taking the continuous limit $N\to \infty$,
we obtain the continuous model of the Kuramoto model,
which is an evolution equation of a probability measure $\rho_t = \rho_t (\theta , \omega )$ 
on $S^1 = [0, 2\pi )$ parameterized by $t \in \mathbf{R}$ and $\omega \in \mathbf{R}$, defined as
\begin{eqnarray}
\left\{ \begin{array}{ll}
\displaystyle \frac{\partial \rho_t}{\partial t} + 
\frac{\partial }{\partial \theta }
\left( \Bigl(\omega  + \frac{K}{2\sqrt{-1}}(\eta (t) e^{-\sqrt{-1}\theta }-\overline{\eta (t)}
       e^{\sqrt{-1}\theta })\Bigr) \rho_t \right) = 0,  \\[0.3cm]
\displaystyle \eta (t) := \int_{\mathbf{R}}
 \! \int^{2\pi}_{0} \! e^{\sqrt{-1}\theta } \rho_t (\theta , \omega ) g(\omega ) d\theta d\omega ,  \\[0.3cm]
\rho_0 (\theta , \omega ) = h(\theta),
\end{array} \right.
\label{conti}
\end{eqnarray}
where $h(\theta)$ is an initial condition and $g(\omega )$ is a given
probability density function for natural frequencies.
We are assuming that the initial condition $h(\theta )$ is independent of $\omega $.
This assumption corresponds to the assumption for the discrete model (\ref{KMN}) that
initial values $\{\theta _j(0)\}_{j=1}^{N}$ and natural frequencies $\{\omega _j\}_{j=1}^{N}$
are independently distributed, and is a physically natural assumption often used in literature.
However, we will also consider $\omega $-dependent initial conditions $h(\theta ,\omega )$,
a probability measure on $S^1$ parameterized by $\omega $, for mathematical reasons, in Sec.7.
Roughly speaking, $\rho_t (\theta , \omega )$ denotes a probability that
an oscillator having a natural frequency $\omega $ is placed at a position $\theta $
(for example, see \cite{Ace, Cra3} for how to derive Eq.(\ref{conti})).
Since $h$ and $\rho_t$ are measures on $S^1$, they should be denoted as $dh(\theta )$ and $d\rho_t(\theta , \,\cdot \,)$,
however, we use the present notation for simplicity.
The $\eta (t)$ is a continuous version of (\ref{order2}), and we also call it the \textit{order parameter}.
$\overline{\eta (t)}$ denotes the complex conjugate of $\eta (t)$.
We can prove that Eq.(\ref{conti}) is a proper continuous model in the sense that the order parameter (\ref{order2})
of the $N$-dimensional Kuramoto model converges to $\eta (t)$ as $N \to \infty$ under some assumptions, 
see Chiba~\cite{Chi4}.
The purpose in this paper is to investigate the dynamics of Eq.(\ref{conti}).

A few properties of Eq.(\ref{conti}) are in order.
It is easy to prove the low of conservation of mass:
\begin{equation}
\int_{\mathbf{R}} \! \int^{2\pi}_{0} \! \rho_t (\theta , \omega ) g(\omega ) d\theta d\omega 
 =  \int_{\mathbf{R}} \! \int^{2\pi}_{0} \! h(\theta)g(\omega ) d\theta d\omega  =1.
\label{g}
\end{equation}
By using the characteristic curve method, Eq.(\ref{conti}) is formally integrated as follows:
Consider the equation
\begin{eqnarray}
\frac{dx}{dt} &=& \omega + \frac{K}{2\sqrt{-1}}(\eta (t) e^{-\sqrt{-1}x}- \overline{\eta (t)} e^{\sqrt{-1}x}),\,\, x\in  [0, 2\pi ),
\label{cha0} 
\end{eqnarray}
which defines a characteristic curve.
Let $x = x(t, s; \theta,\omega  )$ be a solution of Eq.(\ref{cha0}) 
satisfying the initial condition $x(s,s; \theta,\omega  ) = \theta $ at an initial time $s$.
Then, along the characteristic curve, Eq.(\ref{conti}) is integrated to yield
\begin{equation}
\rho_t (\theta , \omega ) = h(x(0,t; \theta,\omega )) \exp \Bigl[ 
\frac{K}{2} \int^t_{0} \! (\eta(s) e^{-\sqrt{-1} x(s,t; \theta,\omega ) }
                              + \overline{\eta(s)} e^{\sqrt{-1} x(s,t; \theta,\omega ) }) ds \Bigr],
\label{cha2}
\end{equation}
see \cite{Chi4} for the proof.
By using Eq.(\ref{cha2}), it is easy to show the equality
\begin{equation}
\int^{2\pi}_{0} \! a(\theta , \omega ) \rho_t (\theta , \omega ) d\theta
 =  \int^{2\pi}_{0} \! a(x(t,0; \theta,\omega ) , \omega ) h(\theta )d\theta,
\label{cha3}
\end{equation}
for any measurable function $a(\theta , \omega )$.
In particular, the order parameter $\eta (t)$ are rewritten as
\begin{equation}
\eta (t)
 =  \int_{\mathbf{R}} \! \int^{2\pi}_{0} \! e^{\sqrt{-1}x(t,0; \theta, \omega)}g(\omega ) h(\theta ) d\theta d\omega .
\label{cha4}
\end{equation}
These expressions will be often used for a nonlinear stability analysis.
Substituting it into Eqs.(\ref{cha0}) and (\ref{cha2}), we obtain
\begin{equation}
\frac{d}{dt}x(t,s; \theta , \omega )
 = \omega  + K \int_{\mathbf{R}} \! \int^{2\pi}_{0} \! 
       \sin \Bigl( x(t,0;\theta ', \omega ') - x(t,s; \theta , \omega )\Bigr) g(\omega ')h(\theta ')d\theta ' d\omega ',
\label{sol1}
\end{equation}
and
\begin{eqnarray}
\rho_t(\theta , \omega ) &=& h(x(0,t; \theta,\omega )) \times \nonumber \\
 & &  \!\!\!\!\! \!\!\!\!\! \!\!\!\!\! \exp \Bigl[ 
K \int^t_{0} \! ds \cdot \int_{\mathbf{R}} \! \int^{2\pi}_{0} \!\cos 
\Bigl( x(s,0;\theta ', \omega ') - x(s,t; \theta , \omega ) \Bigr) h(\theta')g(\omega' ) d\theta' d\omega' \Bigr],
\label{sol2}
\end{eqnarray}
respectively.
They define a system of integro-ordinary differential equations which is equivalent to Eq.(\ref{conti}).
Even if $h(\theta )$ is not differentiable, we consider Eq.(\ref{sol2})
to be a weak solution of Eq.(\ref{conti}).
Indeed, even if $h$ and $\rho_t$ are not differentiable, the quantity (\ref{cha3}) is differentiable with respect to $t$
when $a(\theta , \omega )$ is differentiable.
It is natural to consider the dynamics of weak solutions because $\rho_t$ is a probability measure
and we are interested in the dynamics of its moments, in particular the order parameter.
In \cite{Chi4}, the existence and uniqueness of weak solutions of Eq.(\ref{conti}) is proved.


\section{Transition point formula and the linear instability}

A trivial solution of the continuous model (\ref{conti}),
which is independent of $\theta $ and $t$, is given by the uniform distribution
$\rho_t (\theta , \omega ) = 1/(2\pi)$. In this case, $\eta (t) \equiv 0$. This solution
is called the \textit{incoherent state} or the \textit{de-synchronous state}.
In this section and the next section, we investigate the linear stability of the de-synchronous state.
The nonlinear stability will be discussed in Sec.6.
The analysis of the spectrum of a linear operator obtained from the Kuramoto-type model was first
reported by Strogatz and Mirollo \cite{Str3}.

Let
\begin{equation}
Z_j(t, \omega ) := \int^{2\pi}_{0} \! e^{\sqrt{-1}j\theta } \rho_t (\theta , \omega ) d\theta 
 = \int^{2\pi}_{0} \! e^{\sqrt{-1}jx(t,0;\theta ,\omega )} h(\theta) d\theta
\label{4+1}
\end{equation}
be the Fourier coefficients of $\rho_t (\theta , \omega )$.
Then, $Z_0(t, \omega ) = 1$ and $Z_j$ satisfy the 
differential equations
\begin{eqnarray}
\frac{dZ_1}{dt} = \sqrt{-1}\omega Z_1 + \frac{K}{2} \eta (t) - \frac{K}{2} \overline{\eta (t)} Z_{2},
\label{4-0}
\end{eqnarray}
and
\begin{eqnarray}
\frac{dZ_j}{dt} &=& j\sqrt{-1}\omega Z_j + \frac{jK}{2}( \eta (t) Z_{j-1} - \overline{\eta (t)} Z_{j+1}),
\label{4-0b}
\end{eqnarray}
for $j=2,3,\cdots $.
The order parameter $\eta (t)$ is the integral of $Z_1(t, \omega )$ with the weight $g(\omega )$.
The de-synchronous state corresponds to the trivial solution $Z_j \equiv 0$ for $j = 1, 2, \cdots $.
Eq.(\ref{4+1}) shows $|Z_j(t, \omega )| \leq 1$ and thus $Z_j(t, \omega )$ is
 in the weighted Lebesgue space $L^2 (\mathbf{R}, g(\omega )d\omega )$
for every $t\,$ :
\begin{eqnarray*}
|| Z_j(t, \,\cdot\, ) ||^2_{L^2 (\mathbf{R}, g(\omega )d\omega )} = \int_{\mathbf{R}} \! |Z_j(t, \omega )|^2 g(\omega )d\omega \leq 1.
\end{eqnarray*} 
In order to investigate the linear stability of the trivial solution,
the above equations are linearized around the origin as 
\begin{equation}
\frac{dZ_1}{dt} = \left( \sqrt{-1} \mathcal{M} + \frac{K}{2} \mathcal{P}\right) Z_1, 
\label{4-1}
\end{equation}
and
\begin{eqnarray}
\frac{dZ_j}{dt} = j\sqrt{-1} \mathcal{M}Z_j,
\label{4-1b}
\end{eqnarray}
for $j=2,3,\cdots $,
where $\mathcal{M}: q(\omega ) \mapsto \omega q(\omega )$ is the multiplication operator 
on $L^2 (\mathbf{R}, g(\omega )d\omega )$ and $\mathcal{P}$ is the projection on $L^2 (\mathbf{R}, g(\omega )d\omega )$
defined to be
\begin{equation}
\mathcal{P}q(\omega ) = \int_{\mathbf{R}} \! q(\omega )g(\omega ) d\omega . 
\end{equation} 
If we put $P_0(\omega ) \equiv 1$, $\mathcal{P}$ is also expressed as
$\mathcal{P}q(\omega ) = (q, P_0)$, where $(\,\, , \,\,)$ 
is the inner product on $L^2 (\mathbf{R}, g(\omega )d\omega )$ defined as
\begin{equation}
(q_1, q_2) := \int_{\mathbf{R}} \! q_1(\omega ) \overline{q_2(\omega )} g(\omega )d\omega .
\end{equation}
Note that the order parameter is given as $\eta (t) = \mathcal{P}Z_1 = (Z_1, P_0)$.
To determine the linear stability of the de-synchronous state and the order parameter, 
we have to investigate the spectrum and the semigroup of the operator
$\displaystyle T_1 := \sqrt{-1} \mathcal{M} + \frac{K}{2} \mathcal{P}$.
\\[0.2cm]
\textbf{Remark.}\, We need not assume that the Fourier series $\sum^\infty_{-\infty} Z_j(t, \omega )e^{\sqrt{-1}j \theta }$
converges to $\rho_t (\theta , \omega )$ in any sense. 
It is known that there is a one-to-one correspondence between a measure on $S^1$ and its Fourier coefficients 
(see Shohat and Tamarkin~\cite{Sho}).
Thus the dynamics of $\{Z_j(t, \omega )\}^\infty_{-\infty}$ uniquely determines the dynamics of 
$\rho_t (\theta , \omega )$, and vice versa.
In particular, since a weak solution of the initial value problem (\ref{conti}) is unique (Chiba~\cite{Chi4}),
so is Eqs.(\ref{4-0}),(\ref{4-0b}).
In what follows, we will consider the dynamics of $\{Z_j(t, \omega )\}^\infty_{-\infty}$ instead of $\rho_t$.


\subsection{Analysis of the operator $\sqrt{-1}\mathcal{M}$}

Before investigating the operator $T_1$, we give a few properties of the multiplication operator
$\mathcal{M} : q(\omega ) \mapsto \omega q(\omega )$ on $L^2 (\mathbf{R}, g(\omega )d\omega)$. 
The domain $\mathsf{D}(\mathcal{M})$ of $\mathcal{M}$ is dense in $L^2 (\mathbf{R}, g(\omega )d\omega)$. 
It is well known that its spectrum is given by
$\sigma (\mathcal{M}) = \mathrm{supp} (g) \subset \mathbf{R}$, where $\mathrm{supp} (g)$ is a support of the 
function $g$. Thus the spectrum of $\sqrt{-1} \mathcal{M}$ is 
\begin{equation}
\sigma (\sqrt{-1} \mathcal{M})= \sqrt{-1} \cdot \mathrm{supp} (g)
 = \{ \sqrt{-1} \lambda \, | \, \lambda \in \mathrm{supp} (g) \} \subset \sqrt{-1}\mathbf{R}.
\label{4-3}
\end{equation}
Since $\mathcal{M}$ is selfadjoint, $\sqrt{-1}\mathcal{M}$ generates a $C^0$ semigroup 
$e^{\sqrt{-1}\mathcal{M}t}$ given as $e^{\sqrt{-1}\mathcal{M}t} q(\omega ) = e^{\sqrt{-1}\omega t}q(\omega ) $.
In particular, we obtain
\begin{equation}
(e^{\sqrt{-1}\mathcal{M}t} q_1, q_2)
 = \int_{\mathbf{R}} \! e^{\sqrt{-1}\omega t}q_1(\omega )\overline{q_2(\omega )}g(\omega )d\omega ,
\label{4-4}
\end{equation}
for any $q_1, q_2 \in L^2 (\mathbf{R}, g(\omega )d\omega)$.
This is the Fourier transform of the function $q_1(\omega ) \overline{q_2 (\omega )} g(\omega )$.
Thus if $q_1(\omega ) \overline{q_2 (\omega )} g(\omega )$ is real analytic on $\mathbf{R}$
and has an analytic continuation to the upper half plane,
then $( e^{\sqrt{-1}\mathcal{M}t} q_1, q_2 )$ decays exponentially as $t\to \infty$,
while if $q_1(\omega ) \overline{q_2 (\omega )} g(\omega )$ is $C^r$, then it decays as $O(1/t^r)$
(see Vilenkin~\cite{Vil}).
This means that $e^{\sqrt{-1}\mathcal{M}t}$ does not decay in $L^2 (\mathbf{R}, g(\omega )d\omega)$,
however, it decays to zero in a suitable weak topology.
A weak topology will play an important role in this paper.
These facts are summarized as follows:
\\[0.2cm]
\textbf{Proposition 3.1.}\, A solution of the equation (\ref{4-1b}) with an initial value 
$q(\omega ) \in L^2 (\mathbf{R}, g(\omega )d\omega )$
is given by $Z_j(t, \omega ) = e^{j \sqrt{-1} \mathcal{M} t}q(\omega ) = e^{j \sqrt{-1} \omega  t}q(\omega )$.
The quantity $( e^{j\sqrt{-1}\mathcal{M}t}q_1, q_2)$
decays exponentially to zero as $t\to \infty$ if $g(\omega ), q_1(\omega )$ and $\overline{q_2(\omega )}$
have analytic continuations to the upper half plane.
\\[-0.2cm]

This proposition suggests that analyticity of $g(\omega )$ and initial conditions also plays an important role for 
an analysis of the operator $T_1$.
The resolvent $(\lambda - \sqrt{-1}\mathcal{M})^{-1}$ of the operator $\sqrt{-1}\mathcal{M}$
is calculated as
\begin{equation}
( (\lambda - \sqrt{-1}\mathcal{M})^{-1}q_1, q_2 )
 = \int_{\mathbf{R}} \! \frac{1}{\lambda - \sqrt{-1} \omega } q_1(\omega ) \overline{q_2 (\omega )}g(\omega )d\omega. 
\label{4-7}
\end{equation}
We define the function $D(\lambda )$ to be
\begin{equation}
D(\lambda ) = ( (\lambda - \sqrt{-1}\mathcal{M})^{-1}P_0, P_0 )
 = \int_{\mathbf{R}} \! \frac{1}{\lambda - \sqrt{-1} \omega }g(\omega )d\omega
\label{4-8}
\end{equation}
(recall that $P_0(\omega ) \equiv 1$).
It is holomorphic in $\mathbf{C}\backslash \sigma (\sqrt{-1}\mathcal{M})$ and will be used in later calculations.


\subsection{Eigenvalues of the operator $T_1$ and the transition point formula}

The domain of $T_1 =\sqrt{-1} \mathcal{M} + \frac{K}{2} \mathcal{P}$ is given by 
$\mathsf{D}(\mathcal{M}) \cap \mathsf{D}(\mathcal{P}) = \mathsf{D}(\mathcal{M})$, which is dense in $L^2 (\mathbf{R}, g(\omega )d\omega)$.
Since $\mathcal{M}$ is selfadjoint and $\mathcal{P}$ is bounded, $T_1$ is a closed operator \cite{Kato}.
Let $\mathfrak{\varrho} (T_1)$ be the resolvent set of $T_1$ and 
$\sigma (T_1) = \mathbf{C}\backslash \mathfrak{\varrho} (T_1)$ the spectrum.
Let $\sigma _p(T_1)$ and $\sigma _c(T_1)$ be the point spectrum (the set of eigenvalues)
and the continuous spectrum of $T_1$, respectively.
\\[0.2cm]
\textbf{Proposition 3.2.}\, (i) Eigenvalues $\lambda $ of $T_1$, if they exist, are given as roots of
\begin{equation}
D(\lambda ) = \frac{2}{K},\,\,\, \lambda \in \mathbf{C}\backslash \sigma (\sqrt{-1}\mathcal{M}).
\label{4-10}
\end{equation} 
Furthermore, there are no eigenvalues on the imaginary axis.
\\
(ii) $T_1$ has no residual spectrum.
The continuous spectrum of $T_1$ is given by
\begin{equation}
\sigma _c(T_1) = \sigma (\sqrt{-1} \mathcal{M}) = \sqrt{-1}\cdot \mathrm{supp} (g).
\end{equation}
\\[0.2cm]
\textbf{Proof.} \, (i) Suppose that $\lambda \in \sigma _p(T_1) \backslash  \sigma (\sqrt{-1} \mathcal{M})$.
Then, there exists $x\in L^2 (\mathbf{R}, g(\omega )d\omega)$ such that
\begin{eqnarray*}
\lambda x = (\sqrt{-1}\mathcal{M} + \frac{K}{2} \mathcal{P}) x,\,\,\, x\neq 0.
\end{eqnarray*}
Since $\lambda \notin \sigma (\sqrt{-1} \mathcal{M})$, $(\lambda - \sqrt{-1} \mathcal{M})^{-1}$
is defined and the above is rewritten as
\begin{eqnarray*}
x &=& (\lambda - \sqrt{-1} \mathcal{M})^{-1} \frac{K}{2}\mathcal{P}x  \\
&=& \frac{K}{2}(x, P_0 ) (\lambda - \sqrt{-1} \mathcal{M})^{-1} P_0(\omega ).
\end{eqnarray*}
By taking the inner product with $P_0(\omega )$, we obtain
\begin{equation}
1 = \frac{K}{2} ( (\lambda - \sqrt{-1} \mathcal{M})^{-1}P_0, P_0 ) = \frac{K}{2} D(\lambda ).
\end{equation}
This proves that roots of Eq.(\ref{4-10}) are in $\sigma _p(T_1) \backslash  \sigma (\sqrt{-1} \mathcal{M})$.
The corresponding eigenvector is given by $x = (\lambda - \sqrt{-1} \mathcal{M})^{-1} P_0(\omega )
 = 1/(\lambda - \sqrt{-1}\omega )$.
If $\lambda \in \sqrt{-1}\mathbf{R}$, $x\notin L^2 (\mathbf{R}, g(\omega )d\omega)$.
Thus there are no eigenvalues on the imaginary axis.
In particular, there are no eigenvalues on $\sigma (\sqrt{-1} \mathcal{M})$.
\\
(ii) Since $\mathcal{M}$ is selfadjoint, $\sqrt{-1}\mathcal{M}$ is a Fredholm operator without the residual spectrum.
Since $K$ is $\mathcal{M}$-compact, $T_1$ also has no residual spectrum due to the stability theorem of 
Fredholm operators (see Kato \cite{Kato}).
The latter statement follows from the fact that the essential spectrum is stable under the 
bounded perturbation \cite{Kato}:
the essential spectrum of $T_1$ is the same as $\sigma (\sqrt{-1} \mathcal{M})$.
Since there are no eigenvalues on $\sigma (\sqrt{-1} \mathcal{M})$, it coincides with the continuous spectrum.
\hfill $\blacksquare$
\\[-0.2cm]

Our next task is to calculate roots of Eq.(\ref{4-10}) to obtain eigenvalues of 
$T_1 = \sqrt{-1} \mathcal{M} + \frac{K}{2}\mathcal{P}$.
By putting $\lambda = x + \sqrt{-1}y$ with $ x,y\in \mathbf{R}$,
Eq.(\ref{4-10}) is rewritten as
\begin{equation}
\left\{ \begin{array}{l}
\displaystyle 
\int_{\mathbf{R}} \! \frac{x}{x^2 + (\omega -y)^2}g(\omega )d\omega  = \frac{2}{K},\\[0.4cm]
\displaystyle 
\int_{\mathbf{R}} \! \frac{\omega -y}{x^2 + (\omega -y)^2}g(\omega )d\omega  = 0. \\
\end{array} \right.
\label{4-14}
\end{equation}
The next lemma is easily obtained.
\\[0.2cm]
\textbf{Lemma 3.3.}
\\
(i) If an eigenvalue $\lambda $ exists, it satisfies $\mathrm{Re} (\lambda ) > 0$ for any $K > 0$.
\\
(ii) If $K>0$ is sufficiently large, there exists at least one eigenvalue $\lambda $ near infinity.
\\
(iii) If $K>0$ is sufficiently small, there are no eigenvalues.
\\[0.2cm]
\textbf{Proof.} \, Part (i) of the lemma immediately follows from the first equation of Eq.(\ref{4-14}):
Since the right hand side is positive, $x$ in the left had side has to be positive.
To prove part (ii) of the lemma, note that if $|\lambda |$ is large, Eq.(\ref{4-10}) is expanded as
\begin{eqnarray*}
\frac{1}{\lambda } + O(\frac{1}{\lambda ^2}) = \frac{2}{K}.
\end{eqnarray*}
Thus Rouch\'{e}'s theorem proves that Eq.(\ref{4-10}) has a root $\lambda \sim K/2$ if $K>0$
is sufficiently large.
To prove part (iii) of the lemma, we see that the left hand side of the first equation of 
Eq.(\ref{4-14}) is bounded for any $x, y \in \mathbf{R}$.
To do so, let $G(\omega )$ be the primitive function of $g(\omega )$
and fix $\delta >0$ small. The left hand side of the first equation of Eq.(\ref{4-14}) is calculated as
\begin{eqnarray*}
& & \int_{\mathbf{R}} \! \frac{xg(\omega )d\omega}{x^2 + (\omega -y)^2} \\
 &=& \int^\infty_{y + \delta } \! \frac{xg(\omega )d\omega}{x^2 + (\omega -y)^2}
 + \int^{y-\delta }_{-\infty} \! \frac{xg(\omega )d\omega}{x^2 + (\omega -y)^2}
 + \int^{y + \delta }_{y - \delta } \! \frac{xg(\omega )d\omega}{x^2 + (\omega -y)^2} \\
&=& \int^\infty_{y + \delta } \! \frac{xg(\omega )d\omega}{x^2 + (\omega -y)^2}
 + \int^{y-\delta }_{-\infty} \! \frac{xg(\omega )d\omega}{x^2 + (\omega -y)^2} \\
& & \quad + \frac{x}{x^2 + \delta ^2} \left( G(y + \delta ) - G(y - \delta )\right)
 + \int^{y + \delta }_{y - \delta } \! \frac{2x (\omega -y)}{(x^2 + (\omega -y)^2)^2} G(\omega )d\omega .
\end{eqnarray*}
The first three terms in the right hand side above are bounded for any $x, y\in \mathbf{R}$.
By the mean value theorem, there exists a number $\xi$ such that the last term is estimated as
\begin{eqnarray}
& & \int^{y + \delta }_{y - \delta } \! \frac{2x (\omega -y)}{(x^2 + (\omega -y)^2)^2} G(\omega )d\omega \nonumber \\
&=& \int^{\delta }_{0} \! \frac{2x \omega }{(x^2 + \omega ^2)^2} (G(y+ \omega) - G(y-\omega ))d\omega \nonumber \\
&=& (G(y+ 0) - G(y-0))\int^{\xi}_{0} \! \frac{2x \omega }{(x^2 + \omega ^2)^2}d\omega
 + (G(y+ \delta) - G(y-\delta))\int^{\delta}_{\xi} \! \frac{2x \omega }{(x^2 + \omega ^2)^2} d\omega. \nonumber \\
\label{3-16}
\end{eqnarray}
Since $G$ is continuous, the above is calculated as
\begin{eqnarray*}
(G(y+ \delta) - G(y-\delta)) \left( \frac{x}{x^2 + \xi^2} - \frac{x}{x^2 + \delta ^2 } \right).
\end{eqnarray*}
If $\xi \neq 0$, this is bounded for any $x, y\in \mathbf{R}$.
If $\xi = 0$, Eq.(\ref{3-16}) yields
\begin{eqnarray*}
\int^{\delta }_{0} \! \frac{2x \omega }{(x^2 + \omega ^2)^2} (G(y+ \omega) - G(y-\omega ))d\omega
 = (G(y+ \delta) - G(y-\delta))\int^{\delta}_{0} \! \frac{2x \omega }{(x^2 + \omega ^2)^2} d\omega.
\end{eqnarray*}
Since $G(\omega )$ is monotonically increasing, we obtain
\begin{eqnarray*}
G(y+ \omega) - G(y-\omega ) = G(y+ \delta) - G(y-\delta)
\end{eqnarray*}
for $0 \leq \omega  \leq \delta $. In particular, putting $\omega =0$ gives $G(y+ \delta) - G(y-\delta) = 0$.
Thus $G(y+ \omega) - G(y-\omega ) = 0$ for $0 \leq \omega  \leq \delta $.
This proves that the quantity (\ref{3-16}) is zero.
Now we have proved that the left hand side of the first equation of Eq.(\ref{4-14})
is bounded for any $x, y\in \mathbf{R}$, although the right hand side diverges as $K\to +0$.
Thus Eq.(\ref{4-10}) has no roots if $K>0$ is sufficiently small. 
 \hfill $\blacksquare$
\\[-0.2cm]

\begin{figure}
\begin{center}
\includegraphics[]{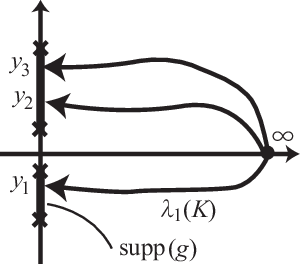}
\caption[]{A schematic view of behavior of roots $\lambda $ of Eq.(\ref{4-10}) when $K$ decreases.
Thick lines denote the continuous spectrum.
As $K$ decreases, eigenvalues $\lambda _1, \lambda _2, \cdots $ converge to $\sqrt{-1}y_1, \sqrt{-1}y_2, \cdots $
and disappear at some $K = K_1, K_2,\cdots $, respectively.}
\label{fig3}
\end{center}
\end{figure}

Lemma 3.3 shows that if $K>0$ is sufficiently large, the trivial solution $Z_1 = 0$ of the equation $dZ_1/dt = T_1Z_1$
is unstable because of eigenvalues with positive real parts.
Our purpose in this section is to determine the bifurcation point $K_c$
such that if $K < K_c$, the operator $T_1$ has no eigenvalues,
while if $K$ exceeds $K_c$, eigenvalues appear on the right half plane
($K_c$ should be positive because of Lemma 3.3 (iii)).
To calculate eigenvalues $\lambda = \lambda (K)$ explicitly is difficult in general.
However, since zeros of the holomorphic function $D(\lambda )-2/K$ do not vanish because of the argument principle,
$\lambda (K)$ disappears if and only if it is absorbed into the continuous spectrum $\sigma (\sqrt{-1}\mathcal{M})$,
on which $D(\lambda )$ is not holomorphic, as $K$ decreases.
This fact suggests that to determine $K_c$, it is sufficient to investigate Eq.(\ref{4-10}) or Eq.(\ref{4-14})
near the imaginary axis.
Thus consider the limit $x \to +0$ in Eq.(\ref{4-14}):
\begin{equation}
\left\{ \begin{array}{l}
\displaystyle 
\lim_{x\to +0}\int_{\mathbf{R}} \! 
\frac{x}{x^2 + (\omega -y)^2}g(\omega )d\omega  = \frac{2}{K},\\[0.4cm]
\displaystyle 
\lim_{x\to +0}\int_{\mathbf{R}} \! 
\frac{\omega -y}{x^2 + (\omega -y)^2}g(\omega )d\omega = 0. \\
\end{array} \right.
\label{4-15}
\end{equation}
These equations determine $K_j$ and $y_j$ such that one of the eigenvalues $\lambda = \lambda _j(K)$
converges to $\sqrt{-1} y_j$
as $K\to K _j + 0$ (see Fig.\ref{fig3}). To calculate them, we need the next lemma.
\\

\noindent \textbf{Lemma 3.4.}\, If $g(\omega )$ is continuous at $\omega = y$, then
\begin{equation}
\lim_{x\to +0}\int_{\mathbf{R}} \! 
\frac{x}{x^2 + (\omega -y)^2}g(\omega )d\omega = \pi g(y ).
\label{4-16}
\end{equation}
\\
\textbf{Proof.} \, This formula is famous and given in Ahlfors~\cite{Ahl}. \hfill $\blacksquare$
\\[-0.2cm]

In what follows, we suppose that $g(\omega )$ is continuous.
Recall that the second equation of Eq.(\ref{4-15})
determines an imaginary part to which $\lambda (K)$ converges as $\mathrm{Re} (\lambda (K)) \to +0$.
Suppose that the number of roots $y_1 , y_2, \cdots$ of the second equation of Eq.(\ref{4-15})
is at most countable for simplicity.
Substituting it into the first equation of Eq.(\ref{4-15}) yields 
\begin{equation}
K_j = \frac{2}{\pi g(y_j)},\,\, j = 1,2,  \cdots ,
\label{4-18}
\end{equation}
which gives the value such that $\mathrm{Re} (\lambda (K)) \to 0$ as $K \to K_j + 0$.
Now we obtain the next theorem.
\\[0.2cm]
\textbf{Theorem 3.5.} Suppose that $g$ is continuous and the number of roots $y_1 , y_2, \cdots$ of 
the second equation of Eq.(\ref{4-15}) is at most countable.
Put
\begin{equation}
K_c:= \inf_{j} K_j = \frac{2}{\pi \sup_{j} g(y_j)}.
\label{4-19}
\end{equation}
If $0 < K \leq K_c$, the operator $T_1$ has no eigenvalues, while if $K$ exceeds $K_c$,
eigenvalues of $T_1$ appear on the right half plane.
In this case, the trivial solution $Z_1 = 0$ of Eq.(\ref{4-1}) is unstable.
\\[-0.2cm]

In general, there exists $K_c^{(2)}$ such that $T_1$ has eigenvalues when $K_c < K< K_c^{(2)}$ but 
they disappear again at $K= K_c^{(2)}$; i.e.
the stability of the trivial solution $Z_1 = 0$ may change many times.
Such $K_c^{(2)}$ is one of the values $K_j$'s.
However, if $g (\omega )$ is an even and unimodal function,
it is easy to prove that $T_1$ has an eigenvalue on the right half plane for any $K > K_c$,
and it is real as is shown in Mirollo and Strogatz \cite{Mir0}.
Indeed, the second equation of Eq.(\ref{4-14}) is calculated as
\begin{eqnarray*}
0 = \int_{\mathbf{R}} \! \frac{\omega -y}{x^2 + (\omega -y)^2} g(\omega )d\omega 
 = \int^\infty_{0} \! \frac{\omega }{x^2 + \omega ^2} (g(y + \omega ) - g(y - \omega )) d\omega . 
\end{eqnarray*}
If $g$ is even, $y = 0$ is a root of this equation.
If $g$ is unimodal, $g(y + \omega ) - g(y - \omega ) > 0$ when $y < 0, \omega >0$
and $g(y + \omega ) - g(y - \omega ) <0$ when $y > 0, \omega >0$. Hence, $y=0$ is a unique root.
This implies that an eigenvalue should be on the real axis, 
and $(K, y) = (K_c, 0)$ is a unique solution of Eq.(\ref{4-15}).
As a corollary, we obtain the transition point (bifurcation point to the synchronous state)
conjectured by Kuramoto~\cite{Kura2}:
\\[0.2cm]
\textbf{Corollary 3.6 (Kuramoto's transition point).}
Suppose that the probability density function $g(\omega )$ is even, unimodal and continuous.
Then, $K_c$ defined as above is given by
\begin{equation}
K_c = \frac{2}{\pi g(0)}.
\label{4-20}
\end{equation}
When $K > K_c$, the solution $Z_1 = 0$ of Eq.(\ref{4-1}) is unstable.
In particular, the order parameter $\eta (t) = (Z_1, P_0)$ is linearly unstable.

\section{Linear stability theory}

Theorem 3.5 shows that $K_c$ is the least bifurcation point and the trivial solution 
$Z_1 = 0$ of Eq.(\ref{4-1}) is unstable if $K$ is larger than $K_c$.
If $0 < K \leq  K_c$, there are no eigenvalues and the continuous spectrum of $T_1$ lies on the imaginary axis: 
$\sigma (T_1) = \sigma (\sqrt{-1}\mathcal{M})$.
In this section, we investigate the dynamics of Eq.(\ref{4-1}) for $0 < K <  K_c$.
We will see that the order parameter $\eta (t)$ may decay exponentially even if the spectrum lies on the imaginary axis
because of the existence of resonance poles.


\subsection{Resonance poles}

Since $\sqrt{-1}\mathcal{M}$ has the semigroup $e^{\sqrt{-1}\mathcal{M}t}$ and since $\mathcal{P}$
is bounded, the operator $\displaystyle T_1 = \sqrt{-1} \mathcal{M} + \frac{K}{2} \mathcal{P}$ also generates
the semigroup $e^{T_1t}$ (Kato~\cite{Kato}) on $L^2(\mathbf{R}, g(\omega )d\omega )$.
A solution of Eq.(\ref{4-1}) with an initial value $\phi(\omega ) \in L^2 (\mathbf{R}, g (\omega )d\omega )$
is given by $e^{T_1t}\phi(\omega )$.
The semigroup $e^{T_1t}$ is calculated by using the Laplace inversion formula
\begin{equation}
e^{T_1t} = \lim_{y\to \infty} \frac{1}{2\pi \sqrt{-1}} \int^{x + \sqrt{-1}y}_{x-\sqrt{-1}y} \!
e^{\lambda t} (\lambda -T_1)^{-1} d\lambda , 
\label{semi1}
\end{equation}
for $t>0$,
where $x > 0$ is chosen so that the contour (see Fig.\ref{fig6} (a)) is to the right of the spectrum of $T_1$ 
(Hille and Phillips~\cite{Hil}, Yosida~\cite{Yos}).
The resolvent $(\lambda -T_1)^{-1}$ is given as follows.
\\[0.2cm]
\textbf{Lemma 4.1.} \, For any $\phi(\omega ), \psi(\omega ) \in L^2 (\mathbf{R}, g (\omega )d\omega )$, the equality
\begin{eqnarray}
& & (  (\lambda -T_1)^{-1} \phi, \psi) \nonumber \\
& =& \!\!\!\! ((\lambda - \! \sqrt{-1} \mathcal{M})^{-1} \phi, \psi )
 + \frac{K/2}{1 - KD(\lambda )/2} ( (\lambda - \! \sqrt{-1} \mathcal{M})^{-1} \phi, P_0)((\lambda - \! \sqrt{-1} \mathcal{M})^{-1}P_0, \psi)
\,\,\,\,\,\,\,\,\,\,\,\,
\label{semi2}
\end{eqnarray}
holds.
\\[0.2cm]
\textbf{Proof.} \, Put
\begin{eqnarray*}
R(\lambda )\phi := (\lambda -T_1)^{-1} \phi = (\lambda - \sqrt{-1}\mathcal{M}
 - \frac{K}{2}\mathcal{P})^{-1}\phi,
\end{eqnarray*}
which yields
\begin{eqnarray*}
(\lambda - \sqrt{-1}\mathcal{M})R(\lambda )\phi
 &=& \phi + \frac{K}{2}\mathcal{P}R(\lambda )\phi 
=  \phi + \frac{K}{2} (R(\lambda )\phi, P_0) P_0.
\end{eqnarray*}
This is rearranged as
\begin{eqnarray}
R(\lambda )\phi =  (\lambda - \sqrt{-1}\mathcal{M})^{-1}\phi 
    + \frac{K}{2} (R(\lambda )\phi, P_0) (\lambda - \sqrt{-1}\mathcal{M})^{-1}P_0.
\label{semi2-2}
\end{eqnarray}
By taking the inner product with $P_0$, we obtain
\begin{eqnarray*}
(R(\lambda )\phi, P_0) =  ((\lambda - \sqrt{-1}\mathcal{M})^{-1}\phi, P_0)
               + \frac{K}{2} (R(\lambda )\phi, P_0) D(\lambda ).
\end{eqnarray*}
This provides
\begin{eqnarray*}
(R(\lambda )\phi, P_0) = \frac{1}{1 - KD(\lambda )/2}( (\lambda - \sqrt{-1} \mathcal{M})^{-1} \phi, P_0).
\end{eqnarray*}
Substituting it into Eq.(\ref{semi2-2}), we obtain Lemma 4.1. \hfill $\blacksquare$
\\[-0.2cm]

Eq.(\ref{semi1}) and Lemma 4.1 show that $(e^{T_1t}\phi, \psi)$ is given by
\begin{eqnarray}
(e^{T_1t}\phi, \psi) &=& \lim_{y\to \infty} \frac{1}{2\pi \sqrt{-1}} \int^{x + \sqrt{-1}y}_{x-\sqrt{-1}y} \!
e^{\lambda t} \Bigl( (\lambda - \! \sqrt{-1} \mathcal{M})^{-1} \phi, \psi ) \nonumber \\
& & + \frac{K/2}{1 - KD(\lambda )/2} ( (\lambda - \! \sqrt{-1} \mathcal{M})^{-1} \phi, P_0)
     ((\lambda - \! \sqrt{-1} \mathcal{M})^{-1}P_0, \psi) \Bigr) d\lambda.
\label{semi3}
\end{eqnarray}
In particular, the order parameter $\eta (t) = (Z_1, P_0)$ for the linearized system (\ref{4-1})
with the initial condition $Z_1(0, \omega ) = \phi (\omega )$ is given by $\eta (t) = (e^{T_1t}\phi, P_0)$.

One of the effective ways to calculate the integral above is to use the residue theorem.
Recall that the resolvent $(\lambda - T_1)^{-1}$ is holomorphic on $\mathbf{C} \backslash \sigma (T_1)$.
When $0 < K \leq K_c$, $T_1$ has no eigenvalues and the continuous 
spectrum lies on the imaginary axis : 
$\sigma (T_1) = \sigma (\sqrt{-1} \mathcal{M}) = \sqrt{-1} \cdot \mathrm{supp} (g)$.
Thus the integrand $e^{\lambda t} ((\lambda -T_1)^{-1}\phi, \psi)$ in Eq.(\ref{semi1}) is holomorphic
on the right half plane and may not be holomorphic on $\sigma (T_1)$.
However, under assumptions below, we can show that the integrand has an analytic continuation 
through the line $\sigma (T_1)$ from the right to the left.
Then, the analytic continuation may have poles on the left half plane (on the second Riemann sheet
of the resolvent), which are called \textit{resonance poles}~\cite{Reed}.
The resonance pole $\lambda $ affects the integral in Eq.(\ref{semi3}) through the residue theorem
(see Fig.\ref{fig6} (b)). In this manner, the order parameter $\eta (t)$ can decay with the exponential rate $e^{\mathrm{Re} (\lambda )t}$.
Such an exponential decay caused by resonance poles is well known in the theory of Schr\"{o}dinger operators \cite{Reed},
and for the Kuramoto model, it is investigated numerically by Strogatz \textit{et al.}~\cite{Str2}
and Balmforth \textit{et al.}~\cite{Bal}.

For an analytic function $\psi (z)$, the function $\psi^*(z)$ is defined by $\psi^*(z) = \overline{\psi (\overline{z})}$.
At first, we construct an analytic continuation of the function 
$F_0(\lambda ) := ( (\lambda -T_1)^{-1} \phi, \psi^*)$
(the function $\psi^*$ instead of $\psi$ is used to avoid the complex conjugate in the inner product).
\\[0.2cm]
\textbf{Lemma 4.2.}\, Suppose that the probability density function $g(\omega )$ and functions
$\phi(\omega ), \psi (\omega )$ are real analytic on $\mathbf{R}$ and they have
meromorphic continuations to the upper half plane.
Then the function $F_0(\lambda ):=( (\lambda -T_1)^{-1} \phi, \psi^*)$ defined on the right half plane has 
the meromorphic continuation $F_1(\lambda )$ to the left half plane, which is given by
\begin{eqnarray}
F_1(\lambda ) &=& ((\lambda - \! \sqrt{-1} \mathcal{M})^{-1} \phi, \psi^* )
  + 2\pi \phi (-\sqrt{-1}\lambda )\psi (-\sqrt{-1}\lambda ) g(-\sqrt{-1}\lambda ) \nonumber \\
& & \quad + \frac{K/2}{1 - KD(\lambda )/2 - \pi K g(-\sqrt{-1}\lambda )}Q[\lambda , \phi] Q[\lambda , \psi],
\label{semi4}
\end{eqnarray}
where $Q[\lambda , \phi]$ is defined to be
\begin{equation}
Q[\lambda , \phi] =( (\lambda - \! \sqrt{-1} \mathcal{M})^{-1} \phi, P_0) + 2\pi g(-\sqrt{-1}\lambda) \phi (-\sqrt{-1}\lambda ).
\label{semi4-2}
\end{equation}

Note that $Q[\lambda , \,\cdot\, ]$ defines a linear functional for each $\lambda \in \mathbf{C}$.
Actually, we will define a suitable function space in Sec.5 so that $Q[\lambda , \,\cdot\, ]$ becomes 
a continuous linear functional (generalized function).
\\[0.2cm]
\textbf{Proof.} \, Define a function $\widetilde{F}(\lambda )$ to be
\begin{equation}
\widetilde{F}(\lambda ) = 
\left\{ \begin{array}{ll}
((\lambda - \sqrt{-1} \mathcal{M})^{-1} \phi , \psi^*)  & (\mathrm{Re} (\lambda ) > 0),  \\
\lim_{\mathrm{Re}(\lambda ) \to +0} ((\lambda - \sqrt{-1} \mathcal{M})^{-1} \phi , \psi^*) 
   & (\mathrm{Re} (\lambda ) = 0),  \\
( (\lambda - \sqrt{-1} \mathcal{M})^{-1} \phi ,\psi^*)  
     + 2\pi \phi (-\sqrt{-1}\lambda ) \psi (-\sqrt{-1}\lambda )g (-\sqrt{-1}\lambda )
   & (\mathrm{Re} (\lambda ) < 0). 
\end{array} \right.
\label{semi4-3}
\end{equation}
By the formula (\ref{4-16}), we obtain
\begin{eqnarray}
& & \lim_{\mathrm{Re}(\lambda ) \to +0} ((\lambda - \sqrt{-1} \mathcal{M})^{-1}\phi ,\psi^*) 
  - \lim_{\mathrm{Re}(\lambda ) \to -0} ((\lambda - \sqrt{-1} \mathcal{M})^{-1} \phi ,\psi^*) \nonumber \\
& &  = 2\pi \phi (\mathrm{Im} (\lambda )) \cdot \psi(\mathrm{Im} (\lambda )) \cdot g(\mathrm{Im} (\lambda )),
\end{eqnarray}
which proves that $\lim_{\mathrm{Re}(\lambda ) \to +0} \widetilde{F}(\lambda )
 = \lim_{\mathrm{Re}(\lambda ) \to -0} \widetilde{F}(\lambda )$.
Therefore, if we show that $\widetilde{F}(\lambda )$ is continuous on the imaginary axis,
then $\widetilde{F}(\lambda )$ is meromorphic on $\mathbf{C}$ by Schwarz's principle of reflection.
To see this, put $\phi (\omega )\psi (\omega )g(\omega ) = q(\omega )$.
By the formula (\ref{4-16}),
\begin{eqnarray*}
& & \lim_{x\to +0} \int^\infty_{-\infty} \! \frac{1}{\lambda - \sqrt{-1} \omega }q(\omega ) d\omega \\
&=&  \lim_{x\to +0} \int^\infty_{-\infty} \! \frac{x}{x^2 + (\omega -y)^2} q(\omega )d\omega 
 + \sqrt{-1} \lim_{x\to +0} \int^\infty_{-\infty} \! \frac{\omega -y}{x^2 + (\omega -y)^2} q(\omega )d\omega  \\
&=& \pi q(y ) - \pi \sqrt{-1}V(y),
\end{eqnarray*}
where $\lambda  = x + \sqrt{-1}y$ and $V(y)$ is the Hilbert transform of $q$ defined by
\begin{equation}
V(y) = \mathrm{p.v.}\frac{1}{\pi} \int^\infty_{-\infty} \! \frac{1}{t} q(y - t) dt, 
\label{hilbert}
\end{equation}
see Chap.VI of Stein and Weiss \cite{Ste}.
Since $q(y)$ is Lipschitz continuous, so is $V(y)$ (Thm.106 of Titchmarsh \cite{Tit}).
This proves that $\lim_{x\to +0} \int^\infty_{-\infty} \! (\lambda - \sqrt{-1}\omega )^{-1}q(\omega ) d\omega$
is continuous in $y$.
Therefore, $\widetilde{F}(\lambda )$ is meromorphic on $\mathbf{C}$.
Now we have obtained the meromorphic continuation of $( (\lambda - \sqrt{-1} \mathcal{M})^{-1} \phi,\psi^*)$
from the right to the left. 
Applying this to Eq.(\ref{semi2}), we obtain the meromorphic continuation of $F_0(\lambda )$ as Eq.(\ref{semi4}).
\hfill $\blacksquare$
\\[-0.2cm]

Eq.(\ref{semi4}) is rewritten as
\begin{eqnarray}
F_1(\lambda ) &=&
\frac{K/2}{1-KD(\lambda )/2 - \pi K g(-\sqrt{-1}\lambda )} \Biggl( \nonumber \\[0.2cm]
& & \!\!\!\!\! (2/K-D(\lambda ))((\lambda - \! \sqrt{-1} \mathcal{M})^{-1} \phi,\psi^* ) 
      + ((\lambda - \! \sqrt{-1} \mathcal{M})^{-1} \phi, P_0 ) \cdot ((\lambda - \! \sqrt{-1} \mathcal{M})^{-1} \psi , P_0) \nonumber \\[0.2cm]
& + &  2\pi g(-\sqrt{-1}\lambda )\Bigl( \frac{2}{K} \phi (-\sqrt{-1}\lambda ) \psi (-\sqrt{-1}\lambda )
   - D(\lambda ) \phi (-\sqrt{-1}\lambda ) \psi (-\sqrt{-1}\lambda ) \nonumber  \\
& &  \qquad \qquad \qquad   -((\lambda - \! \sqrt{-1} \mathcal{M})^{-1} \phi,\psi ^*) 
          + ((\lambda - \! \sqrt{-1} \mathcal{M})^{-1} \phi, P_0) \cdot \psi (-\sqrt{-1}\lambda ) \nonumber \\
& &  \qquad \qquad \qquad  + ((\lambda - \! \sqrt{-1} \mathcal{M})^{-1} \psi ,P_0 ) \cdot \phi (-\sqrt{-1}\lambda ) \Bigr) \Biggr).
\label{F1}
\end{eqnarray}
This expression shows that poles of $g$ are removable.
Therefore, poles of $F_1(\lambda )$ on the left half plane and the imaginary axis are given as roots of the equation
\begin{equation}
\left\{ \begin{array}{ll}
\displaystyle D(\lambda ) + 2\pi g (-\sqrt{-1} \lambda ) = \frac{2}{K}, & \mathrm{Re}(\lambda ) < 0, \\
\displaystyle \lim_{\mathrm{Re}(\lambda )\to +0}D(\lambda ) =
\lim_{\mathrm{Re}(\lambda )\to -0}\left( D(\lambda ) + 2\pi g(-\sqrt{-1}\lambda ) \right) = \frac{2}{K},
&  \mathrm{Re}(\lambda ) =0,\\
\end{array} \right.
\label{semi7}
\end{equation}
and poles of the functions $\phi (-\sqrt{-1}\lambda )$ and $\psi (-\sqrt{-1}\lambda )$.
To avoid dynamics caused by a special choice of $\phi$ and $\psi$, in what follows,
we will assume that continuations of $\phi$ and $\psi$ have no poles.
\\[0.2cm]
\textbf{Definition 4.3.} \, Roots of Eq.(\ref{semi7}) on the left half plane and the imaginary axis
are called \textit{resonance poles} of the operator $T_1$.
\\[-0.2cm]

\begin{figure}
\begin{center}
\includegraphics[]{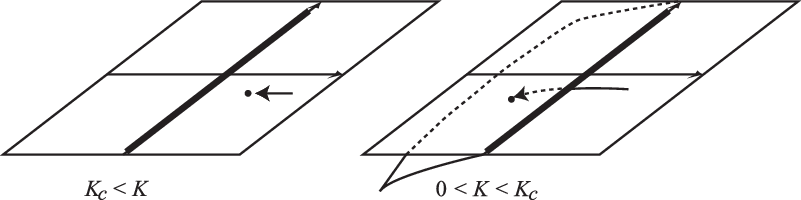}
\caption[]{As is discussed in Sec.3, an eigenvalue $\lambda (K)$ disappears from the original complex plane
at $K = K_c$. But it still exists as a resonance pole on the second Riemann sheet of the resolvent.}
\label{fig5b}
\end{center}
\end{figure}

Since the left hand side of Eq.(\ref{semi7}) is an analytic continuation of that of Eq.(\ref{4-10}), 
at least one of the resonance poles is obtained as a 
continuation of an eigenvalue $\lambda (K)$ coming from the right half plane when $K$ decreases from $K_c$ 
(see Fig.~\ref{fig5b}).
However, if $g(\lambda )$ has an essential singularity, there exist infinitely many resonance poles in general,
which are not obtained as continuations of eigenvalues.

We want to calculate the Laplace inversion formula (\ref{semi3}) by deforming the contour as Fig.\ref{fig6} (b),
and pick up the residues at resonance poles.
We should show that the integral along the arc $C_4$ converges to zero as the radius tends to infinity.
For this purpose, we have to make some assumptions for growth rates of $\phi (\lambda )$ and $\psi (\lambda )$
as $|\lambda | \to \infty$.
Since suitable assumptions depend on the growth rate of $g(\lambda )$,
we calculate the Laplace inversion formula by dividing into two cases:
In Sec.4.2, $g(\omega )$ is assumed to be the Gaussian distribution. 
In Sec.4.3, we consider the case that $g(\omega )$ is a rational function.

\subsection{Gaussian case}

In this subsection, we suppose that $\displaystyle g(\omega ) = e^{-\omega ^2/2}/\sqrt{2\pi}$,
although the results are true for a certain class of density functions.
In this case, the transition point is given by $K_c = 2\sqrt{2/\pi}$.
When $K>K_c$, there exists a unique eigenvalue of $T_1$ on the positive real axis.
When $0<K\leq K_c$, there are no eigenvalues, while resonance poles exist.
The equation (\ref{semi7}) for obtaining the resonance poles is reduced to
\begin{equation}
\int_{\mathbf{R}} \! \frac{1}{\lambda -\sqrt{-1}\omega } g(\omega )d\omega 
 + 2\pi g(-\sqrt{-1}\lambda )
= e^{\lambda ^2/2} \left( \sqrt{\frac{\pi}{2}} - \int^\lambda _{0} \! e^{-x^2/2} dx \right)  = \frac{2}{K}.
\label{A1}
\end{equation}
Let $\lambda _0, \lambda _1, \cdots $ be roots of this equation with 
$0 \geq \mathrm{Re}(\lambda _0) \geq \mathrm{Re}(\lambda _1) \geq \cdots $.
The following properties are easily obtained.
\\[0.2cm]
\textbf{(i)} \, If $\lambda _n$ is a resonance pole, so is its complex conjugate $\overline{\lambda }_n$.
\\
\textbf{(ii)} \, There exist infinitely many resonance poles. As $n\to \infty$, 
$\mathrm{Re}(\lambda _n) \to -\infty$ and they approach to the rays $\mathrm{arg}(z) = 3\pi/4, 5\pi/4$.
\\
\textbf{(iii)} \, When $K=K_c$, there exists a unique resonance pole $\lambda _0 = 0$ on the imaginary axis. 
When $0<K<K_c$, all resonance poles lie on the left half plane.
\\
\textbf{(iv)} All roots of Eq.(\ref{A1}) are simple roots. 

\begin{figure}
\begin{center}
\includegraphics[]{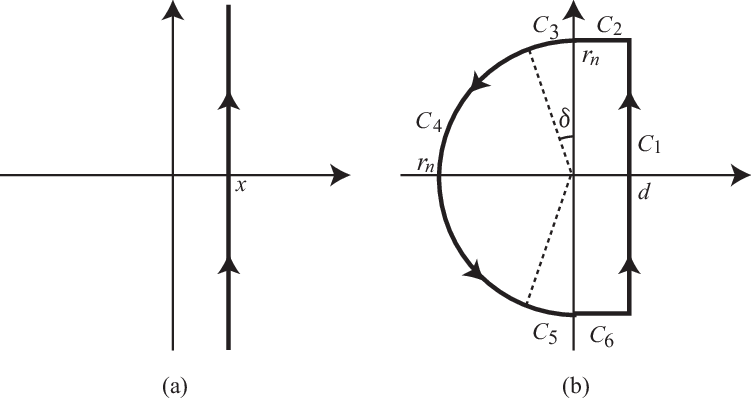}
\caption[]{The contour for the Laplace inversion formula.}
\label{fig6}
\end{center}
\end{figure}

To make assumptions for $\phi$ and $\psi$, we prepare a certain function space.
Let $\mathrm{Exp}_+(\beta, n)$ be the set of holomorphic functions on the region 
$\mathbf{C}_n := \{ z\in \mathbf{C} \, | \, \mathrm{Im} (z) \geq -1/n\}$ such that the norm
\begin{equation}
|| \phi ||_{\beta, n} := \sup_{\mathrm{Im} (z) \geq -1/n} e^{-\beta |z|}|\phi (z)|
\label{exp1}
\end{equation}
is finite. With this norm, $\mathrm{Exp}_+(\beta, n)$ is a Banach space.
Let $\mathrm{Exp}_+ (\beta)$ be their inductive limit with respect to $n=1,2,\cdots $
\begin{equation}
\mathrm{Exp}_+ (\beta) = \varinjlim_{n \geq 1} \mathrm{Exp}_+ (\beta, n) = \bigcup_{n \geq 1} \mathrm{Exp}_+ (\beta, n).
\label{exp2}
\end{equation}
Thus $\mathrm{Exp}_+(\beta)$ is the set of holomorphic functions near the upper half plane that can grow at most the rate $e^{\beta |z|}$.
Next, define $\mathrm{Exp}_+ $ to be their inductive limit with respect to $\beta = 0,1,2,\cdots $
\begin{equation}
\mathrm{Exp}_+ = \varinjlim_{\beta \geq 0} \mathrm{Exp}_+ (\beta) = \bigcup_{\beta \geq 0} \mathrm{Exp}_+ (\beta).
\label{exp3}
\end{equation}
Thus $\mathrm{Exp}_+$ is the set of holomorphic functions near the upper half plane that can grow at most exponentially
; $\phi (z)$ in $\mathrm{Exp}_+$ satisfies $|| \phi ||_{\beta,n} < \infty$ for some $\beta, n$,
and such $\beta$ and $n$ can depend on $\phi$.
Topological properties of $\mathrm{Exp}_+$ will be discussed in Sec.5.2.
In this section, the topology on $\mathrm{Exp}_+$ is not used.
Note that when $\phi \in \mathrm{Exp}_+$, $\phi^*(z) = \overline{\phi (\overline{z})}$ is
holomorphic near the lower half plane and $\phi (- \sqrt{-1}z)$ is holomorphic near the left half plane.
The main theorem in this section is stated as follows.
\\[0.2cm]
\textbf{Theorem 4.4.}\, For any $\phi, \psi \in \mathrm{Exp}_+$, there exists a positive number $t_0$ such that
the semigroup $e^{T_1t}$ satisfies the equality
\begin{equation}
(e^{T_1t}\phi,\psi^*) = S_0[\phi, \psi] e^{\xi_0 t} + \sum^\infty_{n=0} R_n[\phi, \psi] e^{\lambda _n t},
\label{semi8}
\end{equation}
for $t >t_0$, where $\xi_0$ is the eigenvalue of $T_1$ on the right half plane (which exists only when
$K>K_c$), $S_0[\phi, \psi] e^{\xi_0t}$ is a corresponding residue of $F_0(\lambda )e^{\lambda t}$, and
where $\lambda _0, \lambda _1, \cdots $ are resonance poles of $T_1$ 
such that $|\lambda  _0| \leq |\lambda _1| \leq \cdots $,
and $R_n[\phi, \psi] e^{\lambda _nt}$ are corresponding residues of $F_1(\lambda )e^{\lambda t}$.
When $0 < K < K_c$, it is written as
\begin{equation}
(e^{T_1t}\phi,\psi^*) = \sum^\infty_{n=0} R_n[\phi, \psi] e^{\lambda _n t}, \quad \mathrm{Re}(\lambda _n) < 0.
\label{semi8-2}
\end{equation}
In particular, the order parameter $\eta (t) = (e^{T_1t}\phi, P_0)$ for the linearized system (\ref{4-1})
decays to zero exponentially as $t\to \infty$.
\\[0.2cm]
\textbf{Proof.} 
Let $\delta >0$ be a sufficiently small number.
There exist a positive constant $A$ and a sequence $\{r_n\}_{n=0}^\infty$ of positive numbers with 
$r_n \to \infty$ such that
\begin{equation}
\left| 1 - \pi K g(-\sqrt{-1}\lambda ) - \frac{K}{2} D(\lambda )  \right| \geq A, 
\label{4-18b}
\end{equation} 
for $\lambda  = r_ne^{\sqrt{-1}\theta }, \, \pi/2 + \delta < \theta  < 3\pi/2 - \delta $.
Take a positive number $d>0$ so that $\mathrm{Re}(\xi_0) < d$.
With these $d$ and $r_n$, take paths $C_1$ to $C_6$ as are shown in Fig.\ref{fig6} (b):
\begin{eqnarray*}
& & C_1 = \{ d + \sqrt{-1} y \, | \, -r_n \leq y \leq r_n\}, \\
& & C_2 = \{ x + \sqrt{-1}r_n \, | \, 0 \leq x \leq d\}, \\
& & C_3 = \{ r_n e^{\sqrt{-1}\theta } \, | \, \pi /2 \leq \theta \leq \pi / 2 + \delta  \}, \\
& & C_4 = \{ r_n e^{\sqrt{-1}\theta } \, | \, \pi / 2 + \delta \leq \theta \leq 3\pi / 2 - \delta  \},
\end{eqnarray*}
and $C_5$ and $C_6$ are defined in a similar way to $C_3$ and $C_2$, respectively.
We put $C^{(n)} = \sum^6_{j=1}C_j$.

Let $\lambda _0 , \lambda _1 ,\cdots ,\lambda _{\#(n)}$ be resonance poles 
inside the closed curve $C^{(n)}$.
By the definition of $r_n$, there are no resonance poles on the curve $C^{(n)}$.
By the residue theorem, we have
\begin{eqnarray*}
2\pi \sqrt{-1} \left( S_0[\phi, \psi]e^{\xi_0 t}+
 \sum^{\#(n)}_{j=0} R_j[\phi, \psi]e^{\lambda _j t} \right)
 = \!\! \int_{C_6 + C_1 + C_2} \!\!\!\! e^{\lambda t} F_0(\lambda ) d\lambda 
  + \!\!\int_{C_3 + C_4 + C_5} \!\!\!\! e^{\lambda t} F_1(\lambda ) d\lambda,
\end{eqnarray*}
when $r_n$ is sufficiently large so that $C^{(n)}$ encloses the eigenvalue $\xi_0$.
Since the eigenvalue and resonance poles are \textit{simple} roots of Eq.(\ref{4-10}) and Eq.(\ref{semi7}) when $g$
is the Gaussian, $S_0$ and $R_j$ are independent of $t$ (otherwise, they are polynomials in $t$).
Since the integral $\int_{C_1} \! e^{\lambda t} F_0(\lambda ) d\lambda / (2\pi \sqrt{-1})$ 
converges to $(e^{T_1t}\phi,\psi^*)$ as $n \to \infty$, we obtain
\begin{eqnarray}
& & (e^{T_1t}\phi,\psi^*) = S_0[\phi, \psi]e^{\xi_0 t} + 
\lim_{n\to \infty}\sum^{\#(n)}_{j=0} R_j[\phi, \psi] e^{\lambda _jt} \nonumber \\
& & \quad - \frac{1}{2\pi \sqrt{-1}} \lim_{n\to \infty} \int_{C_6 + C_2} \! e^{\lambda t} F_0(\lambda ) d\lambda 
- \frac{1}{2\pi \sqrt{-1}} \lim_{n\to \infty}\int_{C_3 + C_4 + C_5} \! e^{\lambda t} F_1(\lambda ) d\lambda.
\label{semi10}
\end{eqnarray}
It is easy to verify that the integrals along $C_2, C_3, C_5, C_6$ tend to zero as $n\to \infty$.
For example, the integral along $C_2$ is estimated as
\begin{eqnarray*}
\Bigl| \int_{C_2} \! e^{\lambda t}F_0(\lambda ) d\lambda \Bigr| 
 &=&  \Bigl| \int^0_{d} \! e^{(x + \sqrt{-1}r_n)t} F_0(x + \sqrt{-1}r_n) dx \Bigr| \\
&\leq & e^{dt} \int^d_{0} \! |F_0(x + \sqrt{-1}r_n)| dx,
\end{eqnarray*}
where $F_0$ is given as (\ref{semi2}).
Since $|(\lambda - \sqrt{-1}\mathcal{M})^{-1}| \to 0$ as $|\lambda | \to \infty$, 
the integral along $C_2$ proves to be zero as $n\to \infty$ ($r_n \to \infty$).
The integrals along $C_3, C_5, C_6$ are estimated in a similar manner.
The integral along $C_4$ is estimated as
\begin{eqnarray}
\left| \int_{C_4} \! e^{\lambda t} F_1(\lambda ) d\lambda \right| 
&\leq & \int^{3\pi /2 - \delta }_{\pi /2 + \delta }
  \! r_n e^{r_n t \cos \theta } \, |F_1 (r_n e^{\sqrt{-1} \theta })| d\theta  \nonumber \\
&\leq & \max_{\pi / 2 + \delta \leq \theta  \leq 3\pi /2 - \delta } |F_1 (r_n e^{\sqrt{-1}\theta })|
\int^{\pi /2 }_{ \delta } \! 2r_n e^{-r_n t \sin \theta } d\theta  \nonumber \\
&\leq & \max_{\pi / 2 + \delta \leq \theta  \leq 3\pi /2 - \delta } |F_1 (r_n e^{\sqrt{-1}\theta })|
\int^{\pi /2 }_{ \delta } \! 2r_n e^{-2r_n t \theta /\pi} d\theta \nonumber \\
& \leq &  \max_{\pi / 2 + \delta \leq \theta  \leq 3\pi /2 - \delta } |F_1 (r_n e^{\sqrt{-1}\theta })|
   \cdot \frac{\pi }{t} \left( e^{-2r_n t \delta / \pi} - e^{-r_n t}\right) .
\label{semi9}
\end{eqnarray}
Since $|(\lambda - \sqrt{-1}\mathcal{M})^{-1}| \to 0$ as $|\lambda | \to \infty$,
$F_1(\lambda )$ given by Eq.(\ref{F1}) is estimated as
\begin{eqnarray*}
|F_1(\lambda )| \leq \frac{D_0 + |g (-\!\sqrt{-1}\lambda )|\cdot | D_1  + \! D_2 \phi (-\!\sqrt{-1}\lambda)  
+  \!D_3 \psi (-\!\sqrt{-1}\lambda )
 + \!D_4 \phi (-\!\sqrt{-1}\lambda ) \psi (-\!\sqrt{-1}\lambda )|}
{|1 - \pi K g(-\sqrt{-1}\lambda ) - \frac{K}{2} D(\lambda )|} ,
\end{eqnarray*}
where $D_0$ to $D_4$ are some positive constants.
Since $\phi, \psi \in \mathrm{Exp}_+$, there exist $C_1, C_2, \beta_1, \beta_2 \geq 0$ such that
$|\phi (z)| \leq C_1 e^{\beta_1 |z|},\, |\psi (z)| \leq C_2 e^{\beta_2 |z|}$:
\begin{eqnarray*}
|F_1(\lambda )| \leq \frac{D_0 + \left( D_1 + D_2C_1 e^{\beta_1 |\lambda |} + D_3C_2 e^{\beta_2 |\lambda |}
 + D_4C_1C_2 e^{(\beta_1 + \beta_2) |\lambda |} \right)\cdot |g (-\!\sqrt{-1}\lambda )| }
{|1 - \pi K g(-\sqrt{-1}\lambda ) - \frac{K}{2} D(\lambda )|}.
\end{eqnarray*}
Suppose that $g(-\sqrt{-1}r_ne^{\sqrt{-1}\theta })$ diverges as $n\to \infty$ ($r_n \to \infty$).
Then, there exists $D_5 > 0$, which has an upper bound determined by constants $C_1, C_2, \beta_1$ and $\beta_2$,
such that $|F_1 (r_n e^{\sqrt{-1}\theta })|$ is estimated as
\begin{equation}
|F_1 (r_n e^{\sqrt{-1}\theta })| \leq D_5 e^{(\beta_1 + \beta_2) r_n}.
\end{equation}
If $g(-\sqrt{-1}r_ne^{\sqrt{-1}\theta })$ is bounded, Eq.(\ref{4-18b}) shows that 
there exists $D_6 > 0$ , which has an upper bound determined by constants $C_1, C_2, \beta_1$ and $\beta_2$, such that 
\begin{eqnarray}
|F_1(r_n e^{\sqrt{-1}\theta } )| \leq \frac{D_6}{A} e^{(\beta_1 + \beta_2) r_n} \cdot |g (-\!\sqrt{-1}r_n e^{\sqrt{-1}\theta } )|.
\end{eqnarray}
Therefore, we obtain
\begin{equation}
\left| \int_{C_4} \! e^{\lambda t} F_1(\lambda ) d\lambda \right| 
\leq  \frac{D_7 }{t} \left( e^{(\beta_1 + \beta_2 - 2\delta t/\pi) r_n} - e^{(\beta_1 + \beta_2 - t)r_n} \right),
\label{4-22}
\end{equation}
with some $D_7 >0$.
Thus if $t > t_0 := \pi (\beta_1 + \beta_2)/(2\delta )$, 
this integral tends to zero as $n\to \infty$, which proves Eq.(\ref{semi8}).

In particular when $0 < K < K_c$, there are no eigenvalues on the right half plane (Thm.3.5).
Thus Eq.(\ref{semi8}) is reduced to Eq.(\ref{semi8-2}).
\hfill $\blacksquare$
\\[-0.2cm]

Note that $\mathrm{Exp}_+(0)$ is the set of bounded holomorphic functions near the upper half plane.
From the proof above, we immediately obtain the following.
\\[0.2cm]
\textbf{Corollary 4.5.}\, If $\phi, \psi \in \mathrm{Exp}_+(0)$, then Eq.(\ref{semi8}) is true for any $t>0$.


\subsection{Rational case}

In this subsection, we suppose that $g(\omega )$ is a rational function.
Since $g(\omega )$ does not decay so fast as $|\omega | \to \infty$,
we should choose moderate functions for $\phi$ and $\psi$.
Let $\mathbf{C}_+ = \{ z\in \mathbf{C}\, | \, \mathrm{Im}(z) \geq 0\}$ be the real axis and the upper half plane.
Let $H_+$ be the set of bounded holomorphic functions on $\mathbf{C}_+$.
With the norm
\begin{equation}
|| \phi || := \sup_{\mathrm{Im}(z) \geq 0} |\phi (z)|,
\end{equation}
$H_+$ is a Banach space.

It is remarkable that if $g(\omega )$ is a rational function, Eq.(\ref{semi7}) is reduced to an 
algebraic equation. Thus the number of resonance poles is finite.
The proof of the following theorem is similar to that of Thm.4.4 and omitted here.
\\[0.2cm]
\textbf{Theorem 4.6.}\, Suppose that $0 < K \leq K_c$ and $g(\omega )$ is a rational function.
For any $\phi, \psi \in H_+$, the semigroup $e^{T_1t}$ satisfies the equality
\begin{equation}
(e^{T_1t}\phi,\psi^*) = \sum^M_{n=0} R_n[t, \phi, \psi] e^{\lambda _n t},
\label{thm4.6}
\end{equation}
for $t >0$, where $\lambda _0,  \cdots , \lambda _M$ are resonance poles of $T_1$ 
and $R_0[t, \phi, \psi]e^{\lambda _0t}, \cdots , R_M[t, \phi, \psi]e^{\lambda _Mt} $ are 
corresponding residues of $F_1(\lambda )e^{\lambda t}$.
In particular, when $0<K<K_c$, $\lambda _0, \cdots  , \lambda _M$ are on the left half plane and 
the order parameter $\eta (t) = (e^{T_1t}\phi, P_0)$ for the linearized system (\ref{4-1}) decays to zero exponentially as $t\to \infty$.
\\[-0.2cm]

Since the right hand side of Eq.(\ref{thm4.6}) is a finite sum,
the semigroup $e^{T_1t}$ looks like an exponential of a matrix.
The reason of this fact will be revealed in Sec.5.3 by means of the theory of rigged Hilbert spaces.
\\[0.2cm]
\textbf{Example 4.7.} \, If $g(\omega ) = 1/(\pi (1 + \omega ^2))$ is the Lorentzian distribution,
a resonance pole is given by $\lambda  = K/2 -1$ (a root of Eq.(\ref{semi7})).
Therefore $\eta (t)$ decays with the exponential rates $e^{(K/2 -1)t}$.
Note that the transition point is $K_c = 2/\pi/g(0) = 2$.


\section{Spectral theory}

We have proved that when $K>K_c$, the de-synchronous state ($\eta (t) \equiv 0$) is linearly unstable
because of eigenvalues on the right half plane, while when $0 < K < K_c$, it is linearly stable
because of resonance poles on the left half plane.
Next, we want to investigate bifurcations at $K=K_c$.
However, a center manifold in the usual sense is of infinite dimensional because the continuous spectrum
lies on the imaginary axis.
To handle this difficulty, we develop a spectral theory of resonance poles based on a rigged Hilbert space.


\subsection{Rigged Hilbert space}

Let $X$ be a locally convex Hausdorff topological vector space over $\mathbf{C}$ and $X'$ its dual space.
$X'$ is a set of continuous anti-linear functionals on $X$.
For $\mu \in X'$ and $\phi \in X$, $\mu (\phi )$ is denoted by
$\langle \mu \,|\, \phi \rangle$.
For any $a,b \in \mathbf{C},\, \phi, \psi \in X$ and $\mu, \xi \in X'$, the equalities
\begin{eqnarray}
& & \langle \mu \,|\,  a \phi + b\psi\rangle 
   = \overline{a} \langle \mu \,|\,  \phi \rangle + \overline{b} \langle \mu \,|\, \psi \rangle, \\
& & \langle a\mu + b\xi \,|\, \phi \rangle
   = a \langle \mu \,|\, \phi \rangle + b \langle \xi \,|\, \phi \rangle,
\end{eqnarray}
hold. 
Several topologies can be defined on the dual space $X'$.
Two of the most usual topologies are the weak dual topology (weak * topology) 
and the strong dual topology (strong * topology).
A sequence $\{ \mu_j\} \subset X'$ is said to be weakly convergent to $\mu \in X'$
if $\langle \mu_j \,|\, \phi \rangle \to \langle \mu  \,|\, \phi \rangle$ for each $\phi \in X$;
a sequence $\{ \mu_j\} \subset X'$ is said to be strongly convergent to $\mu \in X'$
if $\langle \mu_j \,|\,\phi \rangle \to \langle \mu \,|\, \phi \rangle$ uniformly on any bounded subset of $X$.

Let $\mathcal{H}$ be a Hilbert space with the inner product $(\cdot\, , \, \cdot)$ such that $X$ is a dense subspace of
$\mathcal{H}$.
Since a Hilbert space is isomorphic to its dual space, we obtain $\mathcal{H} \subset X'$ through $\mathcal{H} \simeq \mathcal{H}'$.
\\[0.2cm]
\textbf{Definition 5.1.} If a locally convex Hausdorff topological vector space $X$ is a dense subspace of 
a Hilbert space $\mathcal{H}$ and a topology of $X$ is stronger than that of $\mathcal{H}$, the triplet
\begin{equation}
X \subset \mathcal{H} \subset X'
\end{equation}
is called the \textit{rigged Hilbert space} or the \textit{Gelfand triplet}.
The \textit{canonical inclusion} $i: \mathcal{H} \to X'$ is defined as follows; for $\psi\in \mathcal{H}$,
we denote $i(\psi)$ by $\langle \psi |$, which is defined to be
\begin{equation}
i(\psi)(\phi) = \langle \psi \,|\, \phi \rangle = (\psi, \phi),
\label{cano}
\end{equation}
for any $\phi \in X$. 
Thus if $\mathcal{H} = L^2 (\mathbf{R}, g(\omega )d\omega )$, then
\begin{eqnarray*}
i(\psi)(\phi) = \int_{\mathbf{R}} \! \overline{\phi (\omega )} \psi (\omega ) g(\omega ) d\omega . 
\end{eqnarray*}
We will usually substitute $\phi^*$ instead of $\phi$ to avoid the complex conjugate in the right hand side.
\\[-0.2cm]

Let $A : X \to X$ be a linear operator on $X$.
The (Hilbert) adjoint $A^*$ of $A$ is defined through $(A\phi, \psi) = (\phi , A^* \psi)$ as usual.
If $A^*$ is continuous on $X$, the dual operator $A^\times : X' \to X'$ of $A^*$ defined through
\begin{equation}
\langle A^\times \mu \,|\, \phi \rangle = \langle \mu  \,|\, A^* \phi \rangle,\quad
\phi \in X,\, \mu\in X'
\label{dualop}
\end{equation}
is also continuous on $X'$ for both of the weak dual topology and the strong dual topology.
We can show the equality
\begin{equation}
A^\times i(\psi) = i(A\psi),
\label{dual}
\end{equation}
for any $\psi \in X$, which implies that $A^\times$ is an extension of $A$.

It is easy to show that the canonical inclusion is injective if and only if 
$X$ is a dense subspace of $H$, and the canonical inclusion is continuous 
(for both of the weak dual topology and the strong dual topology) if and only if
a topology of $X$ is stronger than that of $H$ (see Tr\'{e}ves \cite{Tre}).
If $X$ is not dense in $H$, two functionals on $H$ may not be distinguished as functionals on $X$.
As a result, $H' \nsubset X'$ in general.
\\[0.2cm]
\textbf{Definition 5.2.} \,When $X \subset H$ is not a dense subspace of $H$, the triplet $(X, H, X')$
is called the \textit{degenerate} rigged Hilbert space.
\\[-0.2cm]

For applications to the Kuramoto model, we investigate two triplets,
$\mathrm{Exp}_- \subset L^2(\mathbf{R}, g(\omega )d\omega ) \subset \mathrm{Exp}_-'$,
and a degenerate one $(H_-, L^2(\mathbf{R}, g(\omega )d\omega ), H_-')$.


\subsection{Spectral theory on $\mathrm{Exp}_- \subset L^2(\mathbf{R}, g(\omega )d\omega ) \subset \mathrm{Exp}_-'$}

In this subsection, we suppose that $g(\omega )$ is the Gaussian.
Since $g$ decays faster than any exponential functions $e^{-\beta |\omega |}$, we have 
$\mathrm{Exp}_+ \subset L^2(\mathbf{R}, g(\omega )d\omega )$,
and indeed, $\mathrm{Exp}_+$ is dense in $L^2(\mathbf{R}, g(\omega )d\omega )$ and 
the topology of $\mathrm{Exp}_+$ is stronger than that of $L^2(\mathbf{R}, g(\omega )d\omega )$ (see Prop.5.3 below).
Thus the rigged Hilbert space 
$\mathrm{Exp}_+ \subset L^2(\mathbf{R}, g(\omega )d\omega ) \subset \mathrm{Exp}_+'$ is well defined.
Recall that $\mathrm{Exp}_+(\beta, n)$ is a Banach space of holomorphic functions on 
$\mathbf{C}_n=\{ z \, | \, \mathrm{Im}(z)\geq -1/n\}$ with 
the norm $|| \cdot ||_{\beta, n}$, and $\mathrm{Exp}_+(\beta)$ is their inductive limit
with respect to $n \geq 1$.
By the definition of the inductive limit, the topology of $\mathrm{Exp}_+(\beta)$ is defined as follows:
a set $U \subset \mathrm{Exp}_+(\beta)$ is open if and only if $U \cap \mathrm{Exp}_+(\beta, n)$
is open for every $n \geq 1$. 
This implies that the inclusions $\mathrm{Exp}_+(\beta, n) \to \mathrm{Exp}_+(\beta)$ are continuous for every $n \geq 1$.
Similarly, $\mathrm{Exp}_+$ is an inductive limit of $\mathrm{Exp}_+(\beta)$, and its topology is induced
from that of $\mathrm{Exp}_+(\beta)$: a set $U \subset \mathrm{Exp}_+$ is open if and only if $U \cap \mathrm{Exp}_+(\beta)$
is open for every $\beta = 0,1, \cdots $.
The inclusions $\mathrm{Exp}_+(\beta) \to \mathrm{Exp}_+$ are continuous for every $\beta = 0,1,\cdots $.
On the dual space $\mathrm{Exp}_+'$, both of the weak dual topology and the strong dual topology 
can be introduced as was explained. 
The space $\mathrm{Exp}_-$ is defined by $\mathrm{Exp}_- = \{ \phi^* \, | \, \phi \in \mathrm{Exp}_+ \}$,
on which the topology of $\mathrm{Exp}_+$ is introduced by the mapping $\phi \mapsto \phi^*$
(recall that $\phi^*(z):=\overline{\phi(\overline{z})}$).
Then, $\mathrm{Exp}_-$ is an inductive limit of Banach spaces $\mathrm{Exp}_-(\beta, n)$,
which are defined in a similar manner to $\mathrm{Exp}_+ (\beta, n)$.
A Gelfand triplet
\begin{eqnarray}
\mathrm{Exp}_- \subset L^2(\mathbf{R}, g(\omega )d\omega ) \subset \mathrm{Exp}_-'
\label{triplet}
\end{eqnarray}
has the same topological properties as 
$\mathrm{Exp}_+ \subset L^2(\mathbf{R}, g(\omega )d\omega ) \subset \mathrm{Exp}_+'$.
We will use the triplet $\mathrm{Exp}_- \subset L^2(\mathbf{R}, g(\omega )d\omega ) \subset \mathrm{Exp}_-'$,
however, functions in $\mathrm{Exp}_+$ also play an important role.

A topological vector space $X$ is called Montel if every bounded set of $X$ is relatively compact.
A Montel space has a convenient property that on a bounded set $A$ of a dual space of a Montel space,
the weak dual topology coincides with the strong dual topology.
In particular, a weakly convergent series in a dual of a Montel space also converges with respect to 
the strong dual topology (see Tr\'{e}ves \cite{Tre}).
This property is very important to develop a theory of generalized functions.

The topology of $\mathrm{Exp}_+$ has following properties.
Obviously the space $\mathrm{Exp}_-$ has the same properties.
\\[0.2cm]
\textbf{Proposition 5.3.}\, $\mathrm{Exp}_+$ is a topological vector space satisfying
\\[0.2cm]
(i) $\mathrm{Exp}_+$ is a complete Montel space.
\\
(ii) if $\{ \phi_j \}_{j=1}^\infty$ is a convergent series in $\mathrm{Exp}_+$, there exist $n \geq 1$
and $\beta \geq 0$ such that $\{ \phi_j \}_{j=1}^\infty \subset \mathrm{Exp}_+(\beta, n)$ and $\{ \phi_j \}_{j=1}^\infty$
converges with respect to the norm $|| \cdot ||_{\beta, n}$.
\\
(iii) $\mathrm{Exp}_+$ is a dense subspace of $L^2(\mathbf{R}, g(\omega )d\omega )$.
\\
(iv) the topology of $\mathrm{Exp}_+$ is stronger than that of $L^2(\mathbf{R}, g(\omega )d\omega )$.
\\[0.2cm]
\textbf{Proof.} (i)  At first, we prove that $\mathrm{Exp}_+ (\beta)$ is Montel.
To do so, it is sufficient to show that 
the inclusion $\mathrm{Exp}_+(\beta, n) \to \mathrm{Exp}_+(\beta, n+1)$ is a compact operator for every $n$
(see Grothendieck \cite{Gro}, Chap.4.3.3).
To prove it, let $A$ be a bounded set of $\mathrm{Exp}_+(\beta, n)$.
There exists a constant $C$ such that 
$|| \phi ||_{\beta, n} = \sup_{z\in \mathbf{C}_n}e^{-\beta |z|} |\phi (z)| < C$ for any $\phi \in A$.
This means that the set $A$ is locally bounded in the interior of $\mathbf{C}_n$.
Therefore, for any sequence $\{ \phi_j\}^\infty_{j=1} \subset A$, there exists a subsequence
converging to some holomorphic function
$\psi$ uniformly on compact subsets in $\mathbf{C}_n$ (Montel's theorem).
In particular, the subsequence converges to $\psi$ on $\mathbf{C}_{n+1}$, and it satisfies
$|| \psi ||_{\beta, n+1} < C$ and $\psi \in \mathrm{Exp}_+(\beta, n+1)$.
This proves that the inclusion $\mathrm{Exp}_+(\beta, n) \to \mathrm{Exp}_+(\beta, n+1)$ is compact
and thus $\mathrm{Exp}_+ (\beta)$ is Montel.
In a similar manner, we can prove by using Montel's theorem that the inclusion $\mathrm{Exp}_+(\beta) \to \mathrm{Exp}_+(\beta+1)$
is a compact operator for every $\beta = 0,1,\cdots $, which proves that $\mathrm{Exp}_+$ is also Montel.
Next, we show that $\mathrm{Exp}_+$ is complete.
Since $\mathrm{Exp}_+ (\beta,n)$ is a Banach space, in particular it is a DF space,
their inductive limit $\mathrm{Exp}_+(\beta)$ is a DF space by virtue of Prop.5 in Chap.4.3.3 of \cite{Gro},
in which it is shown that an inductive limit of DF spaces is DF.
The same proposition also shows that the inductive limit $\mathrm{Exp}_+$ of DF spaces $\mathrm{Exp}_+(\beta)$
is a DF space. Since $\mathrm{Exp}_+$ is Montel and DF, it is complete because of Cor.2 in Chap.4.3.3 of \cite{Gro}.

(ii) It is known that if the inclusion $\mathrm{Exp}_+(\beta) \to \mathrm{Exp}_+(\beta+1)$ is a compact operator for
every $\beta = 0,1,\cdots $, then, for any bounded set $A \subset \mathrm{Exp}_+$, there exists $\beta$ such that
$A \subset \mathrm{Exp}_+(\beta)$ and $A$ is bounded on $\mathrm{Exp}_+ (\beta)$ (see Komatsu \cite{Kom} and references therein).
By using the same fact again, it turns out that for any bounded set $A \subset \mathrm{Exp}_+$, there exist $\beta$ and $n$ such that
$A \subset \mathrm{Exp}_+(\beta, n)$.
In particular, since a convergent series $\{ \phi_j \}_{j=1}^\infty$ is bounded, there exists $\beta$ and $n$ such that
$\{ \phi_j \}_{j=1}^\infty \subset \mathrm{Exp}_+(\beta, n)$ and it converges with respect to the topology of $\mathrm{Exp}_+(\beta, n)$.

To prove (iii), note that $L^2(\mathbf{R}, g(\omega )d\omega )$ is obtained by the completion of the set of polynomials
because the Gaussian has all moments.
Obviously $\mathrm{Exp}_+$ includes all polynomials, and thus $\mathrm{Exp}_+$ is dense in $L^2(\mathbf{R}, g(\omega )d\omega )$.

For (iv), note that $\mathrm{Exp}_+$ satisfies the first axiom of countability because it is defined through the 
inductive limits of Banach spaces. Therefore, to prove (iv), it is sufficient to show that the inclusion
$\mathrm{Exp}_+ \to L^2(\mathbf{R}, g(\omega )d\omega )$ is sequentially continuous.
Let $\{\phi_j\}_{j=1}^\infty$ be a sequence in $\mathrm{Exp}_+$ which converges to zero.
By (ii), there exist $\beta$ and $n$ such that $\{\phi_j\}_{j=1}^\infty$ converges in the topology of $\mathrm{Exp}_+ (\beta, n)$:
$|| \phi_j ||_{\beta, n} \to 0$. Then,
\begin{eqnarray*}
|| \phi_j ||^2_{L^2(\mathbf{R}, g(\omega )d\omega )} &=&
  \int^\infty_{-\infty} \! |\phi_j (\omega )|^2 g(\omega )d\omega \\
&\leq & \sup_{\omega \in \mathbf{R}} e^{-2\beta |\omega |} |\phi_j (\omega )|^2
          \int^\infty_{-\infty} \! e^{2\beta |\omega |} g(\omega )d\omega \\
&\leq & || \phi_j ||^2_{\beta,n}
          \int^\infty_{-\infty} \! e^{2\beta |\omega |} g(\omega )d\omega.
\end{eqnarray*}
The right hand side exists and tends to zero as $j\to \infty$.
This means that the inclusion $\mathrm{Exp}_+ \to L^2(\mathbf{R}, g(\omega )d\omega )$ is continuous.
\hfill $\blacksquare$
\\[-0.2cm]

The topology of the dual space $\mathrm{Exp}_+'$ has following properties, and so is $\mathrm{Exp}_-'$
\\[0.2cm]
\textbf{Proposition 5.4.}
\\[0.2cm]
(i) $\mathrm{Exp}_+'$ is a complete Montel space with respect to the strong dual topology.
\\
(ii) $\mathrm{Exp}_+'$ is sequentially complete with respect to the weak dual topology; that is,
for a sequence $\{ \mu _j\}^\infty_{j=1}$, if $\langle \mu_j  \,|\, \phi \rangle$ converges to some complex number
$C_\phi$ for every $\phi\in \mathrm{Exp}_+$ as $j\to \infty$, then there exists $\mu \in \mathrm{Exp}_+'$
such that $C_\phi = \langle \mu \,|\, \phi\rangle$ and $\mu_j \to \mu$ with respect to the strong dual topology.
\\[0.2cm]
\textbf{Proof.} (i) It is known that the strong dual of a Montel space is Montel and complete, see Tr\'{e}ves \cite{Tre}.
(ii) Suppose that $\langle \mu_j\,|\, \phi  \rangle$ converges to some complex number
$C_\phi$ for every $\phi\in \mathrm{Exp}_+$. This means that the set $\{ \mu _j\}^\infty_{j=1}$
is weakly bounded and is a Cauchy sequence with respect to the weak dual topology.
As was explained before, on a bounded set of a dual space of a Montel space, the weak dual topology and the strong dual topology
coincide with one another. Thus $\{ \mu _j\}^\infty_{j=1}$ is a Cauchy sequence with respect to the strong dual topology.
Since $\mathrm{Exp}_+'$ is complete with respect to the strong dual topology, $\mu_j$ converges to some element 
$\mu \in \mathrm{Exp}_+'$. In particular, $\langle \mu_j \,|\, \phi \rangle$ converges to $\langle \mu \,|\, \phi \rangle = C_\phi$.
\hfill $\blacksquare$
\\[-0.2cm]

Next, we restrict the domain of the operator $T_1 = \sqrt{-1}\mathcal{M} + \frac{K}{2}\mathcal{P}$
to $\mathrm{Exp}_{+}$. We simply denote $T_1|_{\mathrm{Exp}_{+}}$ by $T_1$.
We will see that $T_1$ is quite moderate if restricted to $\mathrm{Exp}_{+}$.
The next proposition also holds for $\mathrm{Exp}_-$.
\\[0.2cm]
\textbf{Proposition 5.5.}\, 
\\
(i) The operator $T_1 : \mathrm{Exp}_{+} \to \mathrm{Exp}_{+}$ is continuous (note that it is not continuous on 
$L^2(\mathbf{R}, g(\omega )d\omega )$).
\\
(ii) The operator $T_1 : \mathrm{Exp}_{+} \to \mathrm{Exp}_{+}$ generates a holomorphic semigroup
$e^{T_1t} : \mathrm{Exp}_{+} \to \mathrm{Exp}_{+}$ on the positive $t$ axis (note that it is not holomorphic on 
$L^2(\mathbf{R}, g(\omega )d\omega )$).
\\[0.2cm]
\textbf{Proof.} 
(i) It is easy to see by the definition that if $\phi \in \mathrm{Exp}_+$, then $T_1 \phi \in \mathbf{\mathrm{Exp}_+}$.
Let $\{\phi_j\}^\infty_{j=1}$ be a sequence in $\mathrm{Exp}_+$ converging to zero.
By Prop.5.3 (ii), there exist $\beta \geq 0$ and $n\geq 1$ such that $|| \phi_j ||_{\beta, n} \to 0$.
For any $\varepsilon >0$, $|| T_1 \phi_j ||_{\beta + \varepsilon ,n}$ is estimated as 
\begin{eqnarray*}
|| T_1 \phi_j ||_{\beta + \varepsilon ,n} 
 &\leq & || \sqrt{-1}\omega \phi_j  ||_{\beta + \varepsilon ,n}
               + \frac{K}{2}|(\phi_j, P_0)| \cdot || P_0 ||_{\beta + \varepsilon ,n} \\
& \leq & \sup_{\omega \in \mathbf{C}_n} e^{-(\beta + \varepsilon )|\omega |}|\omega \phi_j(\omega )|
         + \frac{K}{2} || \phi_j ||_{L^2(\mathbf{R}, g(\omega )d\omega )} \\
& \leq & || \phi_j ||_{\beta, n} \cdot \sup_{\omega \in \mathbf{C}_n} e^{-\varepsilon |\omega |}|\omega |
         + \frac{K}{2} || \phi_j ||_{L^2(\mathbf{R}, g(\omega )d\omega )}, 
\end{eqnarray*}
which tends to zero as $j\to \infty$.
This proves that $T_1 \phi_j$ tends to zero as $j\to \infty$ with respect to the topology of $\mathrm{Exp}_+$,
and thus $T_1 : \mathrm{Exp}_+ \to \mathrm{Exp}_+$ is continuous.
\\
(ii) We know that the operator $T_1$ generates the semigroup $e^{T_1t}$ as an operator on $L^2(\mathbf{R}, g(\omega )d\omega )$
(see Sec.4.1). In other words, the differential equation
\begin{equation}
\frac{d}{dt}x(t, \omega ) = T_1x(t, \omega )
 = \sqrt{-1}\omega x_1(t,\omega ) + \frac{K}{2}(x(t, \cdot), P_0)
\end{equation}
has a unique solution $x(t, \omega ) = e^{T_1t}\phi (\omega )$ in $L^2(\mathbf{R}, g(\omega )d\omega )$
if an initial condition $\phi$ is in $L^2(\mathbf{R}, g(\omega )d\omega )$.
We have to prove that if $\phi \in \mathrm{Exp}_+$, then $x(t, \cdot ) \in \mathrm{Exp}_+$.
For this purpose, we integrate the above equation as
\begin{equation}
e^{T_1t}\phi (\omega ) = e^{\sqrt{-1}\omega t}\phi(\omega )
 + \frac{K}{2}\int^t_{0} \! e^{\sqrt{-1}\omega (t-s)} (e^{T_1s}\phi, P_0) ds. 
\label{5-6}
\end{equation}
From this expression, it is obvious that if $\phi \in \mathrm{Exp}_+(\beta, n)$, 
then $e^{T_1t}\phi (\omega ) \in \mathrm{Exp}_+$.
By the same way as the standard proof of the existence of holomorphic semigroups \cite{Kato},
we can show that $e^{T_1t}$ is a holomorphic semigroup near the positive real axis.
\hfill $\blacksquare$
\\[-0.2cm]

Eigenvalues of $T_1$ are given as roots of the equation $D(\lambda ) = 2/K,\, \mathrm{Re}(\lambda ) > 0$,
and corresponding eigenvectors are 
\begin{equation}
v_\lambda (\omega ) = \frac{1}{\lambda - \sqrt{-1} \omega } \in L^2(\mathbf{R}, g(\omega )d\omega ).
\label{5-11}
\end{equation}
If we regard it as a functional on $\mathrm{Exp}_-$ through the canonical inclusion 
$i: L^2(\mathbf{R}, g(\omega )d\omega ) \to \mathrm{Exp}_-'$, it acts on $\mathrm{Exp}_-$ as
\begin{equation}
i(v_\lambda )(\phi^*) = \langle v_\lambda  \,|\, \phi^* \rangle = 
 (v_\lambda , \phi^*) = \int_{\mathbf{R}} \! \phi (\omega ) v_\lambda (\omega ) g(\omega ) d\omega 
 = \int_{\mathbf{R}} \! \frac{1}{\lambda - \sqrt{-1}\omega } \phi (\omega ) g(\omega ) d\omega ,
\end{equation}
for $\phi \in \mathrm{Exp}_+$ (i.e. for $\phi^*\in \mathrm{Exp}_-$).
Due to Eq.(\ref{semi4-3}), the analytic continuation of this value from the right to the left is given as
\begin{equation}
\int_{\mathbf{R}} \! \frac{1}{\lambda - \sqrt{-1}\omega } \phi (\omega ) g(\omega ) d\omega
 + 2\pi \phi (-\sqrt{-1}\lambda )g(-\sqrt{-1}\lambda ).
\end{equation}
Motivated by this observation, let us define a linear functional $\mu (\lambda ) \in \mathrm{Exp}_-'$ to be
\begin{equation}
\langle \mu (\lambda ) \,|\,  \phi^* \rangle
 = \int_{\mathbf{R}} \! \frac{1}{\lambda -\sqrt{-1}\omega }\phi (\omega )g(\omega )d\omega 
 +  2\pi \phi (-\sqrt{-1}\lambda )g(-\sqrt{-1}\lambda ),
\label{5-15b}
\end{equation}
when $\mathrm{Re}(\lambda ) < 0$, and 
\begin{equation}
\langle \mu (\lambda ) \,|\, \phi^* \rangle
 = \lim_{x \to +0} 
\int_{\mathbf{R}} \! \frac{1}{(x+ \sqrt{-1}y) -\sqrt{-1}\omega }\phi (\omega )g(\omega )d\omega ,
\label{5-15}
\end{equation}
when $\lambda = \sqrt{-1}y, \, y\in \mathbf{R}$.
It is easy to verify that $\mu (\lambda )$ is continuous and thus an element of $\mathrm{Exp}_-'$.
We expect that $\mu (\lambda )$ plays a similar role to eigenvectors.
Indeed, we can prove the following theorem.
\\[0.2cm]
\textbf{Theorem 5.6.}\, Let $\lambda _0, \lambda _1, \cdots $ be resonance poles of the operator $T_1$
and $T_1^\times$ the dual operator of $T_1$ defined through
\begin{equation}
\langle T_1^\times \xi \,|\, \phi^* \rangle = \langle \xi \,|\, T_1^*\phi^*  \rangle,
\quad \phi \in \mathrm{Exp}_+,\, \xi \in \mathrm{Exp}_-'.
\end{equation}
Then, the equality
\begin{equation}
T_1^\times \mu (\lambda _n) = \lambda _n \mu (\lambda _n)
\label{5-17}
\end{equation}
holds for $n = 0,1,2,\cdots $. In this sense, $\lambda _n$ is an eigenvalue of $T_1^\times$,
and $\mu (\lambda _n)$ is an eigenvector. In what follows, $\mu (\lambda _n)$ is denoted by $\mu_n$
and we call it the \textit{generalized eigenfunction} associated with the resonance pole $\lambda _n$.
\\[0.2cm]
\textbf{Proof.} The proof is straightforward. Suppose that $\mathrm{Re}(\lambda _n) <0$.
For any $\phi \in \mathrm{Exp}_+$,
\begin{eqnarray}
\langle T_1^\times \mu_n \,|\, \phi^* \rangle
&=& \langle \mu_n  \,|\, T_1^* \phi^* \rangle \nonumber \\
&=& \int_{\mathbf{R}} \! \frac{1}{\lambda _n - \sqrt{-1}\omega }(T_1\phi )(\omega )g(\omega )d\omega 
       + 2\pi (T_1\phi )(-\sqrt{-1}\lambda _n)g(-\sqrt{-1}\lambda _n) \nonumber\\
&=&  \int_{\mathbf{R}} \! \frac{\sqrt{-1}\omega }{\lambda _n - \sqrt{-1}\omega }\phi (\omega )g(\omega )d\omega
     + \frac{K}{2}\int_{\mathbf{R}} \! \frac{1}{\lambda _n - \sqrt{-1}\omega }g(\omega )d\omega  \cdot
          \int_{\mathbf{R}} \! \phi (\omega ) g(\omega )d\omega \nonumber\\
& & + 2\pi \left( \lambda _n \phi (-\sqrt{-1}\lambda _n) 
       + \frac{K}{2} \int_{\mathbf{R}} \! \phi (\omega )g(\omega )d\omega  \right) g(-\sqrt{-1}\lambda _n)   \nonumber\\
&=& \lambda _n \left( \int_{\mathbf{R}} \! \frac{1}{\lambda _n - \sqrt{-1}\omega } \phi (\omega ) g(\omega )d\omega
      +  2\pi \phi (-\sqrt{-1}\lambda _n) g(-\sqrt{-1}\lambda _n) \right) \nonumber\\
& & + \frac{K}{2} \int_{\mathbf{R}} \! \phi (\omega ) g(\omega )d\omega 
        \left(  D(\lambda _n) + 2\pi g(-\sqrt{-1}\lambda _n) - \frac{2}{K}\right) \nonumber \\
&=& \lambda _n \langle \mu_n  \,|\, \phi^* \rangle
       + \frac{K}{2} \int_{\mathbf{R}} \! \phi (\omega ) g(\omega )d\omega 
        \left(  D(\lambda _n) + 2\pi g(-\sqrt{-1}\lambda _n) - \frac{2}{K}\right).
\label{thm5.6}
\end{eqnarray}
Since $\lambda _n$ is a resonance pole, it is a root of Eq.(\ref{semi7}).
Thus we obtain
\begin{eqnarray*}
\langle T_1^\times \mu  \,|\, \phi^* \rangle
 = \lambda _n \langle \mu_n \,|\, \phi^* \rangle  = \langle \lambda _n \mu_n \,|\, \phi^* \rangle ,
\end{eqnarray*}
which proves the theorem. The proof for the case $\mathrm{Re}(\lambda _n) = 0$ is done in the same way. \hfill $\blacksquare$
\\[-0.2cm]

Define a dual semigroup $(e^{T_1t})^\times$ through
\begin{equation}
\langle (e^{T_1t})^\times \mu \,|\, \phi ^* \rangle = \langle \mu  \,|\, (e^{T_1t})^* \phi^* \rangle.
\label{5-18}
\end{equation}
for any $\phi \in \mathrm{Exp}_+$ and $\mu \in \mathrm{Exp}_-'$,
where $(e^{T_1t})^*$ is the (Hilbert) adjoint of $e^{T_1t}$.
\\[0.2cm]
\textbf{Proposition 5.7.} \, (i) A solution of the initial value problem
\begin{equation}
\frac{d}{dt} \xi = T_1^\times \xi , \quad \xi (0) = \mu \in \mathrm{Exp}_-'
\label{5-19}
\end{equation}
in $\mathrm{Exp}_-'$ is uniquely given by $\xi (t) = (e^{T_1t})^\times \mu$.
\\
(ii) $(e^{T_1t})^\times$ has eigenvalues $e^{\lambda _0 t}, e^{\lambda _1 t},\cdots $,
where $\lambda _0, \lambda _1,\cdots $ are resonance poles of $T_1$.
\\[0.2cm]
\textbf{Proof.}
This follows from the standard (dual) semigroup theory \cite{Yos}.
\hfill $\blacksquare$
\\[-0.2cm]

If we define a semigroup $e^{T_1^\times t}$ generated by $T_1^\times$ to be the flow of (\ref{5-19}),
then Prop.5.7 (i) means $e^{T_1^\times t} = (e^{T_1t})^\times$.
Prop.5.7 (i) also implies that a solution of the inhomogeneous equation
\begin{equation}
\frac{d}{dt}\xi = T_1^\times \xi + f(t),\quad f(t) \in \mathrm{Exp}_-',
\end{equation}
is uniquely given by
\begin{equation}
\xi (t) = (e^{T_1t})^\times\mu + \int^{t}_{0}\! (e^{T_1(t-s)})^\times f(s) ds.
\label{vari}
\end{equation}
This formula will be used so often when analyzing the nonlinear system (\ref{4-0}),(\ref{4-0b}).

In what follows, we suppose that $0 < K \leq K_c$ so that Eq.(\ref{semi8-2}) is applicable.
Since $R_n[\phi, \psi]$ is the residue of $F_1(\lambda )$ given as Eq.(\ref{semi4}), it is calculated as
\begin{equation}
R_n[\phi, \psi] = \frac{K}{2D_n} \langle \mu_n \,|\, \phi^* \rangle  \langle \mu_n \,|\, \psi^* \rangle ,
\end{equation}
where
\begin{equation}
D_n := \lim_{\lambda \to \lambda _n} \frac{1}{\lambda - \lambda _n}
\left( 1 - \frac{K}{2}D(\lambda ) - \pi K g(-\sqrt{-1}\lambda )\right)
\label{D_n}
\end{equation}
is a constant which is independent of $\phi, \psi$.
Note that $Q[\lambda , \phi ]$ given by Eq.(\ref{semi4-2}) is just the definition of the functional $\mu (\lambda )$.
Thus Eq.(\ref{semi8-2}) is rewritten as
\begin{equation}
(e^{T_1t}\psi , \phi^*) = \sum^\infty_{n=0} \frac{K}{2D_n}e^{\lambda _nt}
\langle \mu_n \,|\, \phi^* \rangle  \langle \mu_n \,|\, \psi^* \rangle,
\label{5-23}
\end{equation}
for $t > t_0$. 
Let $i : L^2(\mathbf{R}, g(\omega )d\omega ) \to \mathrm{Exp}_-'$ be the canonical inclusion with respect to
the triplet (\ref{triplet}).
Since $\mathrm{Exp}_+ \subset L^2(\mathbf{R}, g(\omega )d\omega )$,  $i(\psi) \in \mathrm{Exp}_-'$
is well defined for $\psi \in \mathrm{Exp}_+ $.
Sometimes we will denote $i(\psi)$ by $\psi$ for simplicity.
Thus the left hand side above is rewritten as
\begin{eqnarray}
(e^{T_1t}\psi , \phi^*) = i(e^{T_1t}\psi) (\phi^*) = (e^{T_1t})^\times i(\psi) (\phi^*)
= \langle (e^{T_1t})^\times \psi  \,|\,  \phi^* \rangle.
\label{5-23-2}
\end{eqnarray}
Therefore, we obtain
\begin{equation}
(e^{T_1t})^\times \psi
 = \sum^\infty_{n=0} \frac{K}{2D_n}e^{\lambda _nt} \langle \mu_n \,|\, \psi^* \rangle \mu_n,
\label{5-24}
\end{equation}
for $t > t_0$ and $\psi \in \mathrm{Exp}_+$.
Since Eq.(\ref{5-24}) comes from Eq.(\ref{5-23}), the infinite series in the right hand side of Eq.(\ref{5-24}) 
converges with respect to the weak dual topology on $\mathrm{Exp}_-'$.
However, since $\mathrm{Exp}_-$ is Montel, it also converges with respect to the strong dual topology.
We divide the infinite sum in Eq.(\ref{5-23}) into two parts as
\begin{equation}
(e^{T_1t}\psi , \phi^*) = \sum^M_{n=0} \frac{K}{2D_n}e^{\lambda _nt}
\langle \mu_n \,|\, \phi^* \rangle \langle \mu_n \,|\, \psi^* \rangle
 + \sum^\infty_{n=M+1} \frac{K}{2D_n}e^{\lambda _nt}
\langle \mu_n \,|\, \phi^* \rangle \langle \mu_n \,|\, \psi^* \rangle,
\label{5-25}
\end{equation}
where $M\in \mathbf{N}$ is any natural number.
The second part $\sum^\infty_{n= M+1}[\cdots ]$ does not converge when $0<t < t_0$ in general.
However, since $(e^{T_1t}\psi , \phi^*)$ is holomorphic in $t > 0$ and continuous at $t=0$, we obtain
\begin{equation}
(\psi , \phi^*) = \sum^M_{n=0} \frac{K}{2D_n}
\langle \mu_n \,|\, \phi^* \rangle \langle \mu_n \,|\, \psi^* \rangle
 + \lim_{t\to +0} \sum^\infty_{n=M+1} \frac{K}{2D_n}e^{\lambda _nt}
\langle \mu_n \,|\, \phi^* \rangle \langle \mu_n \,|\, \psi^* \rangle,
\end{equation}
where the second term has a meaning in the sense of an analytic continuation in $t$.
Through the canonical inclusion, the above yields
\begin{eqnarray}
& & i(\psi) = \sum^M_{n=0} \frac{K}{2D_n}
\langle \mu_n \,|\, \psi^* \rangle \cdot  \mu_n 
 + \mathcal{R}_M[\psi], \label{5-27} \\
& & \mathcal{R}_M[\psi] := \lim_{t\to +0} \sum^\infty_{n=M+1} \frac{K}{2D_n}e^{\lambda _nt}
\langle \mu_n \,|\, \psi^* \rangle \cdot  \mu_n , \nonumber 
\end{eqnarray}
which gives the spectral decomposition of $i(\psi) \in i(\mathrm{Exp}_+)$ in $\mathrm{Exp}_-'$.
\\[0.2cm]
\textbf{Theorem 5.8 (Spectral decomposition).}
Suppose that $0<K\leq K_c$.
\\
(i) A system of generalized eigenfunctions $\{  \mu_n \}^\infty_{n=0}$
is complete in the sense that if $\langle \mu_n \,|\, \psi^* \rangle = 0$ for $n= 0,1,\cdots $,
then $\psi =0$.
\\
(ii) $\mu_0, \mu_1, \cdots$ are linearly independent of each other:
if $\sum^\infty_{n=0} a_n \mu_n = 0$ with $a_n\in \mathbf{C}$, then $a_n = 0$ for every $n$.
\\
(iii) Let $V_M$ be a complementary subspace of $\mathrm{span}\{ \mu_0, \cdots , \mu_M \}$ in $\mathrm{Exp}_-'$,
which includes $\mu_j$ for every $j = M+1, M+2 ,\cdots $.
Then, any $i( \psi)\in i(\mathrm{Exp}_+)$ is uniquely decomposed with respect to the direct sum
$\mathrm{span}\{ \mu_0, \cdots , \mu_M \} \oplus V_M$ as Eq.(\ref{5-27}), and this decomposition is independent of 
the choice of the complementary subspace $V_M$ including $\mu_{M+1}, \mu_{M+2}, \cdots $.
\\[0.2cm]
\textbf{Proof.} (i) If $\langle \mu_n \,|\, \psi^* \rangle = 0$ for all $n$, Eq.(\ref{5-23}) provides
$(e^{T_1t}\psi , \phi^*) =0$ for any $\phi \in \mathrm{Exp}_+$.
Since $\mathrm{Exp}_+$ is dense in $L^2(\mathbf{R}, g(\omega )d\omega )$, we obtain $e^{T_1t}\psi = 0$
for $t > t_0$. Since $e^{T_1t}$ is holomorphic in $t > 0$ and strongly continuous at $t = 0$,
we obtain $\psi = 0$ by taking the limit $t\to +0$. 

(ii) Suppose that $\sum^\infty_{n=0} a_n  \mu_n = 0$.
Since $(e^{T_1t})^\times$ is continuous,
\begin{eqnarray*}
0 = (e^{T_1t})^\times \sum^\infty_{n=0} a_n  \mu_n 
 = \sum^\infty_{n=0} a_n (e^{T_1t})^\times  \mu_n 
 = \sum^\infty_{n=0} a_n e^{\lambda _nt} \mu_n.
\end{eqnarray*}
We can assume that
\begin{eqnarray*}
0 \geq \mathrm{Re}(\lambda _0) \geq \mathrm{Re}(\lambda _1) \geq \mathrm{Re}(\lambda _2) \geq \cdots,
\end{eqnarray*}
without loss of generality.
Further, on each vertical line
$\{ \lambda \, | \, \mathrm{Re}(\lambda ) = a\leq 0 \}$, there are only finitely many
resonance poles (see Sec.4.2).
Suppose that $\mathrm{Re}(\lambda _0) = \cdots = \mathrm{Re}(\lambda _k)$
and $\mathrm{Re}(\lambda _k) > \mathrm{Re}(\lambda _{k+1})$.
Then, the above equality provides
\begin{eqnarray*}
0 = \sum^k_{n=0} a_n e^{\sqrt{-1}\mathrm{Im}(\lambda _n)t} \mu_n 
 + \sum^\infty_{n=k+1} a_n e^{(\lambda _n - \mathrm{Re}(\lambda _0))t} \mu_n.
\end{eqnarray*} 
Taking the limit $t\to \infty$ yields
\begin{eqnarray*}
0 = \lim_{t\to \infty}\sum^k_{n=0} a_n e^{\sqrt{-1}\mathrm{Im}(\lambda _n)t} \mu_n.
\end{eqnarray*}
Since the finite set $\mu_0, \cdots , \mu_k$ of generalized eigenfunctions are linearly independent
as in a finite-dimensional case, we obtain $a_n = 0$ for $n = 0, \cdots , k$.
The same procedure is repeated to prove $a_n = 0$ for every $n$.

(iii) Let $V_M$ and $V_M'$ be two complementary subspaces of $\mathrm{span}\{ \mu_0, \cdots , \mu_M \}$,
both of which include $\mu_{M+1}, \mu_{M+2}, \cdots $.
Let 
\begin{eqnarray*}
\mathrm{Exp}_-' = \mathrm{span}\{ \mu_0, \cdots , \mu_M \}\oplus V_M
 = \mathrm{span}\{ \mu_0, \cdots , \mu_M \} \oplus V_M'
\end{eqnarray*}
be two direct sums and let
\begin{eqnarray*}
i(\psi) = \sum^M_{n=0}a_n  \mu_n + v
 = \sum^M_{n=0}a_n' \mu_n + v',\quad v\in V_M,\, v'\in V_M'
\end{eqnarray*}
be corresponding decompositions.
We will use the fact that the decomposition of $(e^{T_1t})^\times \psi$ is uniquely given by
(\ref{5-24}) because of part (ii) of this theorem.
Then, $(e^{T_1t})^\times \psi$ is given by
\begin{eqnarray*}
(e^{T_1t})^\times \psi = 
\sum^M_{n=0}a_n e^{\lambda _n t} \mu_n + (e^{T_1t})^\times v
 = \sum^M_{n=0}a_n' e^{\lambda _n t}\mu_n + (e^{T_1t})^\times v'.
\end{eqnarray*}
They give decompositions of $(e^{T_1t})^\times  \psi$ with respect to two direct sums
\begin{eqnarray*}
\mathrm{span} \{ \mu_0, \cdots  , \mu_n\} \oplus (e^{T_1t})^\times V_M ,
\quad \mathrm{span} \{ \mu_0, \cdots  , \mu_n\} \oplus (e^{T_1t})^\times V'_M,
\end{eqnarray*}
respectively.
Since $(e^{T_1t})^\times \mu_n = e^{\lambda _n t}\mu_n$,
the sets $(e^{T_1t})^\times V_M$ and $(e^{T_1t})^\times V_M'$ also include $\mu_{M+1}, \mu_{M+2}, \cdots $.
Because of part (ii) of the theorem, the decomposition of $(e^{T_1t})^\times  \psi$ with respect to
above direct sums is uniquely given as Eq.(\ref{5-24}) for $t > t_0$.
This implies
\begin{eqnarray*}
a_n = a_n' \,\, ( = \frac{K}{2D_n} \langle  \mu_n  \,|\, \psi^* \rangle ),
\quad (e^{T_1t})^\times v = (e^{T_1t})^\times v'\, \,
( = \sum^\infty_{n=M+1} \frac{K}{2D_n}e^{\lambda _nt}
\langle \mu_n \,|\, \psi^* \rangle \cdot \mu_n ).
\end{eqnarray*}
Since $(e^{T_1t})^\times $ is continuous in $t\geq 0$, we obtain $v = v'$ by the limit $t\to +0$.
\hfill $\blacksquare$
\\[-0.2cm]

When $K>K_c$, there exists an eigenvalue of $T_1$ in the usual sense and
the spectral decomposition involves the eigenvalue.
Eq.(\ref{semi8}) proves
\begin{equation}
(e^{T_1t})^\times \psi
 = \frac{K}{2E_0}e^{\xi_0t}
\langle v_0  \,|\, \psi^* \rangle \cdot  v_0 + 
\sum^\infty_{n=0} \frac{K}{2D_n}e^{\lambda _nt}
\langle \mu_n  \,|\, \psi^* \rangle \cdot \mu_n ,
\end{equation}
for $t> t_0$, and 
\begin{eqnarray}
i(\psi )= \frac{K}{2E_0} \langle v_0 \,|\, \psi^* \rangle \cdot v_0
 + \sum^{M}_{n=0} \frac{K}{2D_n}
\langle \mu_n \,|\, \psi^* \rangle \cdot  \mu_n  + \mathcal{R}_{M}[\psi],
\label{5-27b}
\end{eqnarray}
where $\xi_0$ is an eigenvalue of $T_1$ on the right half plane, 
$v_0$ is a corresponding eigenvector defined by Eq.(\ref{5-11}), and where $E_0$ is defined to be
\begin{equation}
E_0 = \lim_{\xi \to \xi_0} \frac{1}{\xi - \xi_0}\left( 1 - \frac{K}{2}D(\xi)\right).
\label{en}
\end{equation}

Theorem 5.8 suggests the expression of the projection to the generalized eigenspace.
\\[0.2cm]
\textbf{Definition 5.9.}\, 
Denote by $\Pi_n : \mathrm{Exp}_-' \to \mathrm{span}\{ \mu_n \} \subset \mathrm{Exp}_-'$ the projection 
to the generalized eigenspace with respect to the direct sum given in Thm.5.8.
For $ \psi \in i(\mathrm{Exp}_+)$, it is given as
\begin{equation}
\Pi_n \psi = \frac{K}{2D_n} \langle \mu_n \,|\,  \psi^* \rangle \cdot \mu_n.
\label{5-28}
\end{equation}

\begin{figure}
\begin{center}
\includegraphics[]{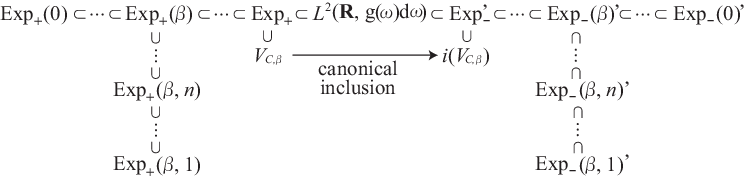}
\caption[]{A diagram for the rigged Hilbert space $\mathrm{Exp}_+ \subset L^2(\mathbf{R}, g(\omega )d\omega ) \subset
\mathrm{Exp}_-'$.}
\label{fig7b}
\end{center}
\end{figure}

Unfortunately, the projection $\Pi_n$ is not a continuous operator.
For example, put $\psi_m (\omega ) = e^{-\sqrt{-1}m\omega }$.
Then, $i(\psi_m)$ converges to zero as $m\to \pm \infty$
with respect to the weak dual topology of $\mathrm{Exp}_-'$ by virtue of the Riemann-Lebesgue lemma.
However,
\begin{eqnarray}
\langle \mu_n \,|\, \psi_m^* \rangle
 = \int_{\mathbf{R}} \! \frac{1}{\lambda_n - \sqrt{-1}\omega }e^{-\sqrt{-1}m\omega }g(\omega )d\omega 
 + 2\pi e^{-m\lambda_n } g(-\sqrt{-1}\lambda_n ) 
\label{5-28b}
\end{eqnarray}
does not tend to zero.
It diverges as $m\to \infty$ when $\mathrm{Re}(\lambda _n) < 0$.
This means that $\Pi_n : \mathrm{Exp}_-' \to \mathrm{Exp}_-'$ is not continuous
with respect to the weak dual topology.
To avoid such a difficulty caused by the weakness of the topology of the domain, we will restrict the domain of $\Pi_n$.
To discuss the continuity, let us introduce the projective topology on $\mathrm{Exp}_-'$
(see also Fig.~\ref{fig7b} and Table 1).
In the dual space $\mathrm{Exp}_-'$, the weak dual topology and the strong dual topology are defined.
Another topology called the projective topology is defined as follows:
Recall that $\mathrm{Exp}_+(\beta,n)$ is a Banach space with the norm $|| \cdot ||_{\beta, n}$, 
and the strong dual $\mathrm{Exp}_- (\beta, n)'$ of $\mathrm{Exp}_- (\beta, n)$ is a
Banach space with the norm
\begin{equation}
|| \mu ||^*_{\beta, n} := \sup_{|| \phi ||_{\beta, n} = 1} |\langle \mu  \,|\, \phi^* \rangle|.
\end{equation}
We introduce the projective topology on $\mathrm{Exp}_-(\beta)' = \bigcap_{n\geq 1} \mathrm{Exp}_-(\beta, n)'$
as follows: $U \subset \mathrm{Exp}_- (\beta)'$ is open if and only if there exist open sets $U_{n} \subset \mathrm{Exp}_- (\beta, n)'$
such that $U_n \cap \mathrm{Exp}_- (\beta)' = U$ for every $n\geq 1$.
It is known that the projective topology is equivalent to that induced by the metric
\begin{equation}
d_\beta (\mu_1, \mu_2) := \sum^\infty_{n=1} \frac{1}{2^n} \frac{|| \mu_1 - \mu_2 ||^*_{\beta, n}}{1 + || \mu_1 - \mu_2 ||^*_{\beta, n}},
\label{db}
\end{equation}
see Gelfand and Shilov \cite{Gel0}.
When the projective topology is equipped, $\mathrm{Exp}_- (\beta)'$ is called the projective limit of $\mathrm{Exp}_-(\beta, n)'$
and denoted by $\mathrm{Exp}_- (\beta)' = \varprojlim \mathrm{Exp}_-(\beta, n)'$.
In a similar manner, the projective topology on $\mathrm{Exp}_-' = \bigcap_{\beta \geq 0} \mathrm{Exp}_-(\beta)'$
is introduced so that $U \subset \mathrm{Exp}_- '$ is open if and only if there exist open sets $U_{\beta} \subset \mathrm{Exp}_- (\beta)'$
such that $U_\beta \cap \mathrm{Exp}_- ' = U$ for every $\beta\geq 0$.
This topology coincides with the topology induced by the metric
\begin{equation}
d(\mu_1, \mu_2) := \sum^\infty_{\beta =0} \frac{1}{2^{\beta}} \frac{d_{\beta} (\mu_1, \mu_2)}{1+ d_{\beta} (\mu_1, \mu_2)}.
\label{5-32}
\end{equation}
In this manner, $\mathrm{Exp}_-'$ equipped with the projective topology is a complete metric vector space.
When the projective topology is equipped, $\mathrm{Exp}_-'$ is called the projective limit of $\mathrm{Exp}_-(\beta)'$
and denoted by $\mathrm{Exp}_-' = \varprojlim \mathrm{Exp}_-(\beta)'$.

By the definition, the projective topology on $\mathrm{Exp}_-'$ is weaker than the strong dual topology
and stronger than the weak dual topology.
Since $\mathrm{Exp}_-$ is a Montel space, the weak dual topology coincides with the strong dual topology
on any bounded set of $\mathrm{Exp}_-'$.
This implies that the projective topology also coincides with the weak dual topology and the strong dual topology
on any bounded set of $\mathrm{Exp}_-'$.
In particular, a weakly convergent series in $\mathrm{Exp}_-'$ also converges with respect to the metric $d$ 
and the strong dual topology.

\begin{table}[h]
\begin{center}
\begin{tabular}{|c|c|}
\hline 
& \\[-0.3cm]
$\displaystyle \mathrm{Exp}_+ (\beta, n)$ & 
    Banach space: $\displaystyle || \phi ||_{\beta, n} = \sup_{\mathrm{Im}(z) \geq -1/n}|\phi (z)| e^{-\beta |z|}$ \\ \hline
& \\[-0.3cm]
$\displaystyle \mathrm{Exp}_- (\beta, n)'$ & 
   Banach space:  $\displaystyle || \xi ||^*_{\beta, n} = \sup_{|| \phi ||_{\beta,n} =1 }
         |\langle \xi  \,|\, \phi^* \rangle| $ \\ \hline 
& \\[-0.3cm]
$\displaystyle \mathrm{Exp}_-(\beta)' = \varprojlim \mathrm{Exp}_- (\beta, n)'$ &  
    $\displaystyle d_\beta (\xi ,\zeta) 
               = \sum^\infty_{n=1} \frac{1}{2^n} \frac{|| \xi -\zeta ||^*_{\beta, n}}{1 + || \xi -\zeta ||^*_{\beta, n}}$ \\ \hline
$\displaystyle \mathrm{Exp}_-' = \varprojlim \mathrm{Exp}_- (\beta)'$ &  
    $\displaystyle d(\xi ,\zeta) = \sum^\infty_{\beta=0} \frac{1}{2^\beta} \frac{d_\beta (\xi ,\zeta)}{1 + d_\beta (\xi ,\zeta)}$ \\ \hline
\end{tabular}
\caption{Metric vector spaces used in Section 5.}
\end{center} 
\end{table}

For constants $C\geq 1, \alpha  \geq 0$ and $0 < p \leq \infty$, 
define a subset $V_{C, \alpha, p} \subset \mathrm{Exp}_+$ to be
\begin{equation}
V_{C, \alpha , p}  = \{ \phi \in \mathrm{Exp}_+ \, | \, 
|\phi (z)|e^{-\alpha  |z|} \leq C \,\, \mathrm{when} \,\, 0\leq \mathrm{Im}(z) \leq 2p \}.
\label{Vca}
\end{equation}
When the choice of a number $0 < p \leq \infty$ is not important, we denote it as $V_{C, \alpha }$ for simplicity.
Note that the set $\{e^{-\sqrt{-1}m\omega }\}^\infty_{m=0}$ above is \textit{not} included in $V_{C,\alpha }$
for any $C$ and $\alpha $.
Let $i(V_{C, \alpha })$ be an inclusion into $\mathrm{Exp}_-'$.

\noindent \textbf{Theorem 5.10.}\, For any $C\geq 1, \alpha  \geq 0$ and $p>0$, the following holds.
\\
(i) On $i(V_{C, \alpha }) \subset \mathrm{Exp}_-'$, 
the weak dual topology, the projective topology and the strong dual topology coincide with one another.
\\
(ii) The closure $\overline{i(V_{C, \alpha })}$ of $i(V_{C, \alpha })$ is a connected, compact metric space.
\\
(iii) For the system (\ref{conti}), give an initial condition $\rho_0(\theta , \omega ) = h(\theta )$,
where $h$ is an arbitrary probability measure on $S^1$.
Then, corresponding solutions of (\ref{4-0}),(\ref{4-0b}) satisfy $Z_j(t, \, \cdot \,) \in V_{1,0}$ for 
any $t\geq 0$ and $j=1,2,\cdots $
(In particular, $Z_j(t, \,\cdot \,)\in V_{1,0, \infty}$).
\\[0.2cm]
\textbf{Proof.}\, (i) At first, we show that the set $i(V_{C, \alpha }) \subset \mathrm{Exp}_-'$ is equicontinuous.
For any small $\varepsilon >0$, we define a neighborhood $U = U(\varepsilon ) \subset \mathrm{Exp}_+$
of the origin so that if $\phi \in U \cap \mathrm{Exp}_+ (\gamma )$, then
\begin{eqnarray*}
\sup_{\mathrm{Im}(z) \geq 0} e^{-\gamma |z|} |\phi (z)| < \varepsilon D(\alpha , \gamma ),
\end{eqnarray*}
where $D(\alpha , \gamma )$ is a positive number to be determined.
Then, for any $\phi \in U \cap \mathrm{Exp}_+(\gamma )$ and $\psi \in i(V_{C, \alpha })$,
\begin{eqnarray*}
|\langle \psi \,|\, \phi^* \rangle| 
&\leq & \int_{\mathbf{R}} \! |\phi (\omega )| \cdot |\psi (\omega )| g(\omega ) d\omega  \\
&\leq & \int_{\mathbf{R}} \! e^{-\gamma |\omega |} |\phi (\omega )| \cdot e^{-\alpha  |\omega |} |\psi (\omega )| 
          \cdot e^{(\alpha  + \gamma )|\omega |} g(\omega ) d\omega \\
&\leq & \varepsilon CD(\alpha , \gamma ) \int_{\mathbf{R}} \! e^{(\alpha  + \gamma )|\omega |} g(\omega ) d\omega.
\end{eqnarray*}
Since $g(\omega )$ is the Gaussian, the integral $E(\alpha , \gamma )
:= \int_{\mathbf{R}} \! e^{(\alpha  + \gamma )|\omega |} g(\omega ) d\omega$ exists.
If we put $D(\alpha , \gamma ) = 1/CE(\alpha , \gamma )$, we obtain $|\langle \phi \,|\, \psi \rangle| < \varepsilon $
for any $\phi\in U$ and $\psi \in i(V_{C, \alpha })$.
This proves that $i(V_{C, \alpha })$ is an equicontinuous set.
In particular, $i(V_{C, \alpha })$ is a bounded set of $\mathrm{Exp}_-'$ for both of the weak dual topology 
and the strong dual topology due to the Banach-Steinhaus theorem (see Prop.32.5 of Tr\'{e}ves \cite{Tre}).
Since $\mathrm{Exp}_-$ is Montel, the weak dual topology, the projective topology and the strong dual topology coincide on the 
bounded set $i(V_{C, \alpha })$. Thus it is sufficient to prove (ii) for one of these topologies.

(ii) Obviously $V_{C, \alpha } \subset \mathrm{Exp}_+$ is connected (actually it is a convex set).
Since the canonical inclusion $i$ is continuous, $i(V_{C, \alpha })$ and $\overline{i(V_{C, \alpha })}$
are connected.
Since the strong dual $\mathrm{Exp}_-'$ is Montel (Prop.5.4), every bounded set of $\mathrm{Exp}_-'$ is relatively
compact, which proves that $\overline{i(V_{C, \alpha })}$ is compact.
By the projective topology, $\overline{i(V_{C, \alpha })}$ is a metrizable space with the metric (\ref{5-32}). 

(iii) To prove $Z_j\in V_{1, 0}$, recall that $Z_j$ is defined by Eq.(\ref{4+1}).
We want to estimate the analytic continuation of $Z_j(t, \omega )$ with respect to $\omega $.
Put $X(t) = e^{\sqrt{-1}x(t,0;\theta , \omega )}$. From Eq.(\ref{cha0}), it turns out that $X$ satisfies the equation
\begin{eqnarray*}
\left\{ \begin{array}{ll}
\displaystyle \frac{dX}{dt} = \sqrt{-1}\omega X + \frac{K}{2}\eta (t) - \frac{K}{2}\overline{\eta (t)}X^2,  \\
X(0) = e^{\sqrt{-1}\theta }.  \\
\end{array} \right.
\end{eqnarray*}
Put $X(t) = \xi (t)e^{\sqrt{-1}p(t)},\, \eta (t) = \zeta (t) e^{\sqrt{-1} q(t)}$ with $\xi, \zeta, p, q\in \mathbf{R}$.
Then, the above equation is rewritten as 
\begin{eqnarray*}
\frac{d\xi}{dt} + \sqrt{-1}\frac{dp}{dt} \xi 
 = (\sqrt{-1} \mathrm{Re}(\omega ) - \mathrm{Im}(\omega )) \xi
    + \frac{K}{2}\zeta e^{\sqrt{-1}(q-p )} - \frac{K}{2} \zeta \xi^2 e^{-\sqrt{-1}(q-p)},
\end{eqnarray*}
which yields
\begin{eqnarray}
\frac{d\xi}{dt} = - \mathrm{Im}(\omega ) \xi + \frac{K}{2}\zeta (1 - \xi^2) \cos (q- p).
\label{xi}
\end{eqnarray}
This equation shows that if $\mathrm{Im}(\omega ) \geq 0$ and $|\xi| = 1$, then $d\xi/dt \leq 0$.
Therefore, if the initial condition satisfies $|X(0)| \leq 1$, then $|X(t)| \leq 1$ for any $t >0$ and any $\mathrm{Im}(\omega ) \geq 0$.
Thus Eq.(\ref{4+1}) shows that 
the analytic continuation of $Z_j(t, \omega )$ to the upper half plane is estimated as
\begin{eqnarray}
|Z_j(t, \omega )| \leq \int^{2\pi}_{0} \! |X(t)|^j h(\theta ) d\theta  \leq 1, \quad j=1,2,\cdots ,
\label{xi2}
\end{eqnarray}
which means that $Z_j(t, \omega ) \in V_{1,0}$ for every $t \geq 0$.
\hfill $\blacksquare$
\\[-0.2cm]

Although solutions of the system (\ref{4-0}),(\ref{4-0b}) are included in the set $V_{1,0}$,
this set is inconvenient because it is not closed under the multiplication by a scalar.
Let us introduce a new set $W_{C, \alpha }$.
For $C\geq 1, \alpha \geq 0$ and $p>0$, we define a subset $W_{C, \alpha }$ of $\mathrm{Exp}_+$ to be
\begin{equation}
W_{C, \alpha }:=\{\psi \in \mathrm{Exp}_+ \, | \, | \psi (z +\sqrt{-1}p)/\psi (z) |e^{-\alpha |z|}\leq C
,\,\, \mathrm{when}\,\, 0\leq \mathrm{Im}(z) \leq p \}.
\end{equation}
The choice of a number $p>0$ is not important.
If $\xi \in i(W_{C, \alpha })$, then $k\xi \in i(W_{C, \alpha })$ for any $k\in \mathbf{C}$.
For elements in $i(W_{C, \alpha })$, let us estimate the norm of the projection $\Pi_j$.
\\[0.2cm]
\textbf{Lemma 5.11.}
(i) For each $\xi \in i(\mathrm{Exp}_+)$, $|| \xi ||^*_{\beta, n}$ is bounded as $n\to \infty$.
\\
(ii) For every $\beta = 0,1,\cdots $ and $n = 1,2, \cdots $, there exists a positive number $Q(\beta)$
such that the inequalities
\begin{equation}
|| \xi ||^*_{0, 1} \leq || \xi ||^*_{\beta, n}, \quad  || \xi ||^*_{\beta, n} \leq Q(\beta)|| \xi ||^*_{0, 1}
\label{compara}
\end{equation}
hold for $\xi \in i(\mathrm{Exp}_+)$ (this means that norms $|| \cdot ||^*_{\beta ,n}$ are comparable \cite{Gel0}).
\\
(iii) For $\mu (\lambda )$ defined by Eq.(\ref{5-15b}), the linear mapping 
$i(\psi) \mapsto \langle \mu (\lambda ) \,|\, \psi^* \rangle$
from $i(W_{C, \alpha })$ into $\mathbf{C}$ is continuous with respect to the projective topology
when $-p<\mathrm{Re}(\lambda ) \leq 0$. In particular, if a resonance pole $\lambda _j$ satisfies $-p < \mathrm{Re}(\lambda _j) \leq 0$,
the corresponding projection $\Pi_j$ is continuous on $i(W_{C, \alpha })$.
\\
(iv) For every $\beta = 0,1,\cdots $ and $n = 1,2, \cdots $, there exists a positive number $D_{C, \alpha , \beta, j}$
such that the inequality
\begin{equation}
|| \Pi_j \xi ||^*_{\beta, n} \leq D_{C, \alpha ,\beta, j} || \xi ||^*_{\beta, n} 
\label{thm5-17}
\end{equation}
holds for $\xi \in i(W_{C, \alpha })$.
\\[0.2cm]
\textbf{Proof.} 
(i) $|| \xi ||^*_{\beta, n}$ has an upper bound
\begin{eqnarray*}
|| \xi ||^*_{\beta, n}:= \sup_{|| \phi ||_{\beta, n} = 1} \left| \int^{}_{\mathbf{R}}\!
\phi (\omega )\xi (\omega )g(\omega )d\omega \right| \leq \int_{\mathbf{R}}\! e^{\beta |\omega |} |\xi (\omega )|g(\omega )d\omega, 
\end{eqnarray*}
which is independent of $n = 1,2,\cdots $.

(ii) The inequality $|| \xi ||^*_{0, 1} \leq || \xi ||^*_{\beta, n}$ follows from the definition.
It is easy to verify that the inclusion $\mathrm{Exp}_-(0,1) \to \mathrm{Exp}_-(\beta, n)$ is continuous.
Thus its dual operator from $\mathrm{Exp}_-(\beta, n)'$ into $\mathrm{Exp}_-(0,1)'$ is continuous.
This shows that there exists a positive number $Q(\beta)$ such that $|| \xi ||^*_{\beta, n} \leq Q(\beta)|| \xi ||^*_{0, 1}$.
Since the norm $|| \cdot ||^*_{\beta, n}$ is bounded as $n\to \infty$, we can take $Q(\beta)$ not to depend on $n=1,2,\cdots $.

(iii) Let $\{ \psi _m \}^\infty_{m=1} \subset i(W_{C, \alpha })$ be a sequence converging
to zero as $m\to \infty$ with respect to the projective topology.
By the definition of the projective topology, we have $|| \psi _m ||^*_{\beta, n} \to 0$
for every $\beta $ and $n$.
This means that $\langle \psi_m \,|\, f^* \rangle \to 0$ uniformly in $f\in \mathrm{Exp}_+(\beta, n)$
satisfying $|| f ||_{\beta, n} \leq C$ for each $C>0$ and $\beta \geq 0$.
Due to the part (i) of the lemma, $|| \psi _m ||^*_{\beta, n} \to 0$ uniformly in $n=1,2,\cdots $,
which shows that $\langle \psi_m \,|\, f^* \rangle \to 0$ uniformly in $f\in V_{C, \beta}$
for each $C>0$ and $\beta \geq 0$.
For a positive number $p>0$ satisfying $-p < \mathrm{Re}(\lambda )$, $\langle \mu (\lambda ) \,|\, \psi _m^* \rangle$ is given by
\begin{eqnarray*}
\langle \mu (\lambda ) \,|\, \psi _m^* \rangle 
&=& \int^{}_{\mathbf{R}}\! \frac{1}{\lambda - \sqrt{-1}\omega } \psi_m(\omega )g(\omega )d\omega 
      + 2\pi \psi_m(-\sqrt{-1}\lambda )g(-\sqrt{-1}\lambda ) \\
&=& \int_{\mathbf{R}}\! \frac{1}{\lambda  - \sqrt{-1}(\omega + \sqrt{-1}p)}\psi_m(\omega + \sqrt{-1}p) g(\omega + \sqrt{-1}p)d\omega  \\
&=& \int_{\mathbf{R}}\! \frac{1}{\lambda  - \sqrt{-1}(\omega + \sqrt{-1}p)}
\frac{\psi_m(\omega + \sqrt{-1}p) g(\omega + \sqrt{-1}p)}{\psi_m(\omega )g(\omega )}\psi_m(\omega )g(\omega )d\omega ,
\end{eqnarray*}
where we used the residue theorem. Putting
\begin{eqnarray*}
f_m(\omega ) := \frac{1}{\lambda  - \sqrt{-1}(\omega + \sqrt{-1}p)}
\frac{\psi_m(\omega + \sqrt{-1}p) g(\omega + \sqrt{-1}p)}{\psi_m(\omega )g(\omega )}
\end{eqnarray*}
provides $\langle \mu (\lambda ) \,|\, \psi_m^* \rangle = \langle \psi_m \,|\, f_m^* \rangle$.
Since $\psi_m \in W_{C, \alpha }$, there exist $C'\geq 1$ and $\alpha ' \geq 0$ such that
$f_m\in V_{C', \alpha '}$ for every $m$. 
Therefore, $\langle \mu(\lambda )\,|\, \psi^*_m \rangle \to 0$ as $m\to 0$.
Since the projective topology is metrizable, this implies that the mapping 
$i(\psi ) \mapsto \langle \mu(\lambda )\,|\, \psi^* \rangle$ is continuous.

(iv) Since $\Pi_j$ is continuous on $i(W_{C, \alpha })$ with respect to the metric $d$, for any $\varepsilon >0$,
there exists a number $\delta _{C, \alpha , j} > 0$ such that if $d(\xi, 0)\leq \delta _{C, \alpha , j}$, then 
$d(\Pi_j \xi, 0) \leq \varepsilon $ for $\xi \in i(W_{C, \alpha })$.
For $\varepsilon >0$, take $\eta \in i(W_{C, \alpha })$ and numbers $\hat{\delta }_{C, \alpha , \beta ,j}$ such that
$|| \eta ||^*_{\beta, n} \leq \hat{\delta }_{C, \alpha , \beta ,j}$.
We can take $\hat{\delta }_{C, \alpha , \beta ,j}$ sufficiently small so that
$d(\eta ,0) \leq \delta _{C, \alpha , j}$ holds, which implies $d(\Pi_j \eta, 0) \leq \varepsilon$.
By the definition of $d$, this yields
\begin{eqnarray*}
||\Pi_j \eta ||^*_{\beta, n}\leq \frac{2^n\kappa}{1-2^n \kappa},\quad \kappa := \frac{2^\beta \varepsilon }{1-2^\beta \varepsilon }.
\end{eqnarray*}
If $\xi \in i(W_{C, \alpha })$, then $\eta := \hat{\delta }_{C, \alpha , \beta ,j}\xi / || \xi ||^*_{\beta ,n}$
is included in $i(W_{C, \alpha })$ and satisfies 
$|| \eta ||^*_{\beta, n} \leq \hat{\delta }_{C, \alpha , \beta ,j}$.
Thus we obtain
\begin{equation}
\Bigl| \Bigl| \Pi_j \frac{\hat{\delta }_{C, \alpha , \beta ,j} \xi }{|| \xi ||^*_{\beta ,n}}\Bigl| \Bigl| ^*_{\beta ,n} 
\leq \frac{2^n \kappa}{1-2^n \kappa},
\end{equation}
for $\xi \in i(W_{C, \alpha })$, which yields Eq.(\ref{thm5-17}) by putting 
$D_{C,\alpha ,\beta, j} = 2^n \kappa /(1-2^n \kappa ) /\hat{\delta }_{C, \alpha , \beta, j}$.
Since the norm $|| \cdot ||^*_{\beta, n}$ is bounded as $n\to \infty$, we can take $D_{C,\alpha ,\beta, j}$ not to depend on $n=1,2,\cdots $.
\hfill $\blacksquare$
\\[-0.2cm]

Define the \textit{generalized center subspace} of $T_1$ to be 
\begin{equation}
\mathbf{E}_c = \mathrm{span}\{ \, \mu_n \, | \, \lambda _n\in \sqrt{-1}\mathbf{R}\} \subset \mathrm{Exp}_-'.
\label{center}
\end{equation}
When $g$ is the Gaussian distribution and $K=K_c$, $\mathbf{E}_c$ is a one dimensional vector space
because  Eq.(\ref{semi7}) has a unique root $\lambda _0 = 0$ on the imaginary axis when $K=K_c$.
Let $\mu_0$ be the corresponding generalized eigenfunction; 
that is, $\mathbf{E}_c = \mathrm{span}\{ \mu_0 \}$.
Let $\mathbf{E}_c^{\bot}$ be a complementary subspace of $\mathbf{E}_c$ in $\mathrm{Exp}_-'$ 
including $\mu_{1}, \mu_{2}, \cdots $.
Let $\Pi_c : \mathrm{Exp}_-' \to \mathbf{E}_c$ be the projection to $\mathbf{E}_c$ with respect to the direct sum
$\mathrm{Exp}_-' = \mathbf{E}_c \oplus \mathbf{E}_c^{\bot}$.
Although $\mathbf{E}_c^{\bot}$ may not be unique, $\Pi_c \psi $ is uniquely determined
for $\psi \in i(\mathrm{Exp}_+)$ because of Thm.5.8 (iii).
The complementary subspace $\mathbf{E}_c^{\bot}$ including $\mu_{1}, \mu_{2}, \cdots $ 
is called the \textit{stable subspace}.
Eq.(\ref{5-24}) shows that $\Pi_s (e^{T_1t})^\times \psi$ decays exponentially as $t\to \infty$,
because $\mathrm{Re} (\lambda _j) < 0$ for $j=1,2,\cdots $,
where $\Pi_s = id - \Pi_c$ is the projection to $\mathbf{E}_c^{\bot}$.
\\[0.2cm]
\textbf{Theorem 5.12.}\, For any $ \psi \in i(\mathrm{Exp}_+)$, 
the projection to the center subspace $\Pi_c $ satisfies
\begin{eqnarray}
& & \Pi_c  \psi = \Pi_0 \psi = \frac{K_c}{2D_0}\langle \mu_0  \,|\, \psi^* \rangle \cdot \mu_0. \\
& & \Pi_c T_1^\times \psi = T_1^\times \Pi_c  \psi, \\
& & \Pi_c (e^{T_1t})^\times \psi = (e^{T_1t})^\times \Pi_c  \psi.
\label{proj}
\end{eqnarray}
\textbf{Proof.}\, The first equality follows from the definition.
The second one is verified by using Eq.(\ref{dual}) and $T_1^\times \mu_n = \lambda _n \mu_n$ as
\begin{eqnarray*}
\Pi_c T_1^\times \psi &=& \Pi_c T_1^\times i(\psi) = \Pi_c i(T_1 \psi) \\
& & \quad  = \frac{K_c}{2D_0}\langle \mu_0 \,|\, T_1^* \psi^* \rangle \mu_0
 = \frac{K_c}{2D_0}\langle T_1^\times \mu_0 \,|\, \psi^* \rangle \mu_0
 = \frac{K_c}{2D_0}\langle \lambda _0 \mu_0 \,|\, \psi^* \rangle \mu_0, \\
T_1^\times \Pi_c \psi &=& \frac{K_c}{2D_0}\langle \mu_0 \,|\, \psi^* \rangle T_1^\times \mu_0
 = \frac{K_c}{2D_0}\langle \mu_0 \,|\, \psi^* \rangle \lambda _0 \mu_0.
\end{eqnarray*}
The third one is proved in the same way. \hfill $\blacksquare$
\\

Let $\overline{i(W_{C, \alpha })} \subset \mathrm{Exp}_-'$ be a closure of $i(W_{C, \alpha })$
with respect to the projective topology.
The next proposition shows that solutions of the system (\ref{4-0}),(\ref{4-0b}) are included in the 
closure of the set $W_{3,0}$.
\\[0.2cm]
\textbf{Proposition 5.13.} \, 
(i) For any $C\geq 1$, $i(V_{C, 0}) \subset \overline{i(W_{3, 0})}$.
\\
(ii) Put $V = \bigcup_{C\geq 1}V_{C, 0}$. 
Then, the generalized center subspace $\mathbf{E}_c$ is included in $\overline{i(V)}$;
\begin{equation}
\mathbf{E}_c \subset \overline{i(V)} \subset \overline{i(W_{3, 0})} \subset \mathrm{Exp}_-'.
\end{equation}
(iii) $\Pi_c : i(W_{3,0}) \to \mathrm{Exp}_-'$ is continuous with respect to the projective topology.
The continuous extension $\tilde{\Pi}_c : \overline{i(W_{3, 0})}\to \mathrm{Exp}_-'$ satisfies
$\tilde{\Pi}_c \circ \tilde{\Pi}_c = \tilde{\Pi}_c$.
\\[0.2cm]
\textbf{Proof.} \,
(i) If a function $\psi\in V_{C, 0}$ has zeros on the region $0\leq \mathrm{Im}(z) \leq p$,
$\psi \notin W_{C, \alpha }$ for any $C$ and $\alpha $.
To prove that $i(\psi) \in \overline{i(W_{3,0})}$, let us perturb the function $\psi \in V_{C, 0}$.
For $n=1,2,\cdots $, put
\begin{equation}
\tilde{\psi}(\omega ) = \psi (\omega ) + 2Ce^{\sqrt{-1}n\omega  + np}.
\end{equation}
For $0\leq \mathrm{Im}(\omega ) \leq p$, we have
\begin{eqnarray*}
\left| \frac{\tilde{\psi} (\omega +\sqrt{-1}p)}{\tilde{\psi}(\omega )} \right|
 = \left| \frac{\psi (\omega + \sqrt{-1}p) + 2C e^{\sqrt{-1}n\omega }}{\psi (\omega ) + 2Ce^{\sqrt{-1}n\omega + np}} \right|
\leq \frac{C + 2Ce^{-n \mathrm{Im}(\omega )}}{2Ce^{n(p-\mathrm{Im}(\omega ))} - C} \leq 3,
\end{eqnarray*}
which implies $\tilde{\psi} \in W_{3,0}$.
It is easy to verify that $2Ce^{\sqrt{-1}n\omega  + np} \to 0$ as $n\to \infty$ with respect to the weak dual topology.
Therefore, $i(\tilde{\psi}) \to i(\psi) \in \overline{i(W_{3,0})}$  as $n\to \infty$ for any $C\geq 1$.

(ii) Put $v_\lambda (\omega ) = 1/(\lambda -\sqrt{-1}\omega )$.
Let $\lambda _0 = 0$ be a resonance pole on the imaginary axis.
By the definition, the corresponding generalized eigenfunction $\mu_0$ is given by
\begin{equation}
\mu_0 = \lim_{x\to +0}i(v_{x }(\omega )),
\label{mu}
\end{equation}
where the limit is taken with respect to the weak dual topology.
It is easy to verify that $v_\lambda (\omega ) \in V$ for $\mathrm{Re}(\lambda ) >0$.
This implies that $\mu_0 \in \overline{i(V)}$ and thus
the generalized center subspace $\mathbf{E}_c$ is included in $\overline{i(V)}$.

(iii) The continuity was proved in Lemma 5.11.
Since $\bm{E}_c = \mathrm{span} \{ \mu_0\}$, it is sufficient to prove that $\tilde{\Pi}_c \mu_0 = \mu_0$.
Since 
\begin{eqnarray*}
\tilde{\Pi}_c \mu_0 = \lim_{x\to +0} \Pi_c i(v_x) 
        = \lim_{x\to +0} \frac{K_c}{2D_0} \langle \mu_0 \,|\, v_x^* \rangle \mu_0,
\end{eqnarray*}
let us show that
\begin{eqnarray*}
\lim_{x\to +0} \frac{K_c}{2D_0} \langle \mu_0 \,|\, v_x^* \rangle = 1.
\end{eqnarray*}
By the definition of $\mu_0$ given in Eq.(\ref{5-15}), we have
\begin{eqnarray*}
\lim_{x\to +0} \langle \mu_0 \,|\, v_x^* \rangle 
 = \lim_{x\to +0} \lim_{x'\to +0} \int_{\mathbf{R}}\! \frac{1}{x' - \sqrt{-1}\omega } \frac{1}{x-\sqrt{-1}\omega }
    g(\omega )d\omega  .
\end{eqnarray*}
By the definition of $D_0$ given in Eq.(\ref{D_n}), we obtain
\begin{eqnarray*}
\frac{2D_0}{K_c} 
&=& \frac{2}{K_c} \lim_{x\to -0}\frac{1}{x} \left( 1-\frac{K_c}{2}D(\lambda ) - \pi K_c g(-\sqrt{-1}\lambda ) \right)
 = \frac{2}{K_c} \lim_{x\to +0}\frac{1}{x} \left( 1-\frac{K_c}{2}D(\lambda )  \right) \\
&=& \frac{2}{K_c} \lim_{x\to +0}\frac{1}{x} 
        \left( 1-\frac{K_c}{2} \int_{\mathbf{R}}\! \frac{1}{x- \sqrt{-1}\omega }g(\omega )d\omega  \right).
\end{eqnarray*}
Since $K_c = 2/\pi g(0)$ (Cor.3.6), 
\begin{eqnarray*}
\frac{2D_0}{K_c}  =  \lim_{x\to +0}\frac{1}{x} 
\left( \pi g(0) -  \int_{\mathbf{R}}\! \frac{1}{x- \sqrt{-1}\omega }g(\omega )d\omega \right).
\end{eqnarray*}
Lemma 3.4 yields
\begin{eqnarray*}
\frac{2D_0}{K_c}  =  \lim_{x\to +0}\frac{1}{x} \left(  \lim_{x'\to +0}
\int_{\mathbf{R}}\! \frac{x'}{(x')^2 + \omega ^2}g(\omega )d\omega 
   - \int_{\mathbf{R}}\! \frac{1}{x- \sqrt{-1}\omega }g(\omega )d\omega \right).
\end{eqnarray*}
Since $g$ is an even function, the above is rearranged as
\begin{eqnarray*}
\frac{2D_0}{K_c} &=& \lim_{x\to +0}\lim_{x'\to +0} \frac{1}{x} 
    \left( \int_{\mathbf{R}}\! \frac{1}{x'- \sqrt{-1}\omega }g(\omega )d\omega 
             - \int_{\mathbf{R}}\! \frac{1}{x- \sqrt{-1}\omega }g(\omega )d\omega\right) \\
&=& \lim_{x\to +0}\lim_{x'\to +0} \frac{1}{x} \int_{\mathbf{R}}\! 
             \frac{x-x'}{(x' - \sqrt{-1}\omega ) (x-\sqrt{-1}\omega )} g(\omega )d\omega \\
&=& \lim_{x\to +0}\lim_{x'\to +0} \int_{\mathbf{R}}\! 
             \frac{1}{(x' - \sqrt{-1}\omega ) (x-\sqrt{-1}\omega )} g(\omega )d\omega.
\end{eqnarray*}
This completes the proof. \hfill $\blacksquare$
\\

In what follows, the extension $\tilde{\Pi}_c$ of $\Pi_c|_{i(W_{3,0})}$ is denoted by $\Pi_c$ for simplicity.
Next purpose is to estimate norms of semigroups.
At first, we suppose that $K<K_c$.
In this case, there are no resonance poles on the imaginary axis and thus $\Pi_c =0$.
\\[0.2cm]
\textbf{Proposition 5.14.}
Suppose that $0<K<K_c$.
For every $\beta = 0,1,\cdots $ and $n = 1,2, \cdots $, there exist positive numbers $M_{C, \alpha , \beta}$ and $a$ such that
the inequality
\begin{equation}
|| (e^{T_jt})^\times \xi ||^*_{\beta, n} \leq M_{C, \alpha ,\beta} e^{-jat} || \xi ||^*_{\beta, n}, \quad j= 1,2,\cdots 
\label{thm5-18}
\end{equation}
holds for $\xi \in \overline{i(W_{C, \alpha })}$.
\\[0.2cm]
\textbf{Proof.}
At first, we show the proposition for $j=1$.
When $g$ is the Gaussian, there exists a positive constant $a = a(K)$ such that all resonance poles satisfy 
$\mathrm{Re}(\lambda _n) < -a$.
Then,
\begin{eqnarray*}
\frac{|| (e^{T_1t})^\times \xi||^*_{\beta ,n}}{|| \xi ||^*_{\beta ,n}}
 = e^{-at} \frac{|| (e^{(T_1 + a)t})^\times \xi||^*_{\beta ,n}}{|| \xi ||^*_{\beta ,n}},
\end{eqnarray*}
and $(e^{(T_1 + a)t})^\times \xi/ || \xi ||^*_{\beta ,n}$ is given by
\begin{eqnarray*}
(e^{(T_1 + a)t})^\times \frac{\xi}{|| \xi ||^*_{\beta ,n}}
&=& \sum^\infty_{n=0} e^{(\lambda _n + a)t} \Pi_n \left( \frac{\xi}{|| \xi ||^*_{\beta ,n}} \right).
\end{eqnarray*}
for $\xi \in i(\mathrm{Exp}_+)$.
Due to Lemma 5.11 (iv), the right hand side above is bounded with respect to the norm $|| \cdot ||^*_{\beta, n}$
uniformly in $\xi \in i(W_{C, \alpha })$ and $t > 0$.
This proves that there exists a positive constant $L_{C, \alpha , \beta}$ such that 
$|| (e^{T_1t})^\times \xi||^*_{\beta ,n} \leq e^{-at} L_{C, \alpha , \beta} || \xi ||^*_{\beta ,n}$
for $\xi \in i(W_{C, \alpha })$.
Since the norm $|| \cdot ||^*_{\beta, n}$ is bounded as $n\to \infty$, we can take $L_{C, \alpha , \beta}$ not to depend on $n=1,2,\cdots $.
The result is continuously extended to the closure $\overline{i(W_{C, \alpha })}$.

Next, let us consider $T_j$ for $j = 2,3,\cdots $.
Cauchy's theorem proves that
\begin{eqnarray*}
\langle (e^{T_jt})^\times \phi \,|\, \psi^*  \rangle
  &=& (e^{j\sqrt{-1}\omega t} \phi, \psi^* ) 
  = \int_{\mathbf{R}} \!  e^{j\sqrt{-1}\omega t} \phi (\omega ) \psi (\omega ) g(\omega )d\omega \\
  &=& \int^{\sqrt{-1}a + \infty}_{\sqrt{-1}a - \infty} \! e^{j\sqrt{-1}\omega t} \phi (\omega ) \psi (\omega ) g(\omega )d\omega \\
  &=& e^{-j a t} \int^\infty_{-\infty} \! e^{j\sqrt{-1}\omega t} 
           \phi (\omega + \sqrt{-1}a ) \psi (\omega + \sqrt{-1}a) g(\omega + \sqrt{-1}a)d\omega,
\end{eqnarray*}
for any $\phi, \psi \in \mathrm{Exp}_+$.
Hence, we obtain
\begin{eqnarray*}
e^{jat} || (e^{T_jt})^\times \phi ||^*_{\beta, n} 
&=& e^{jat} \sup_{|| \psi ||_{\beta, n} = 1} |\langle (e^{T_jt})^\times \phi \,|\, \psi^*  \rangle | \\
&=& \sup_{|| \psi ||_{\beta, n} = 1}\left| \int_{\mathbf{R}}\! f_{\psi}(\omega )\phi (\omega )g(\omega )d\omega \right|
= \sup_{|| \psi ||_{\beta, n} = 1} |\langle \phi \,|\, f_{\psi}^*  \rangle |,
\end{eqnarray*}
where
\begin{eqnarray*}
f_\psi (\omega ):= e^{j\sqrt{-1}\omega t} \psi (\omega +\sqrt{-1}a) 
 \frac{\phi (\omega +\sqrt{-1}a) g(\omega +\sqrt{-1}a)}{\phi (\omega ) g(\omega )}.
\end{eqnarray*}
If $\beta '$ is sufficiently large and $\phi \in i(W_{C, \alpha })$, $|| f_\psi ||_{\beta ', n}$ exists and
\begin{eqnarray*}
e^{jat} || (e^{T_jt})^\times \phi ||^*_{\beta, n}
 = \sup_{|| \psi ||_{\beta, n} = 1}\frac{|\langle \phi \,|\, f_{\psi}^*  \rangle |}{|| f_\psi ||_{\beta ', n}}
         || f_\psi ||_{\beta ', n}
\leq \sup_{|| \psi ||_{\beta, n} = 1}|| f_\psi ||_{\beta ', n} \cdot || \phi ||^*_{\beta ', n}.
\end{eqnarray*}
Due to Lemma 5.11(ii), it turns out that there exists a positive number $N_{C,\alpha , \beta}$
such that $e^{jat}|| (e^{T_jt})^\times \phi ||^*_{\beta ,n} \leq  N_{C, \alpha , \beta} || \phi ||^*_{\beta, n}$
for $\phi \in i(W_{C, \alpha })$.
The result is continuously extended to the closure $\overline{i(W_{C, \alpha })}$.
Then, putting $M_{C, \alpha ,\beta} = \max\{L_{C, \alpha ,\beta}, N_{C, \alpha ,\beta} \}$ yields the desired result. 
\hfill $\blacksquare$
\\[-0.2cm]

Next, we suppose that $K=K_c$. In this case, there exists a resonance pole on the imaginary axis and $\Pi_c \neq 0$.
Then, we can prove the following proposition.
\\
\textbf{Proposition 5.15.}\,
Suppose that $K=K_c$.
Then, for every $\beta = 0,1,\cdots $ and $n = 1,2, \cdots $, 
there exist positive constants $L_{C,\alpha , \beta}, M_{C, \alpha , \beta}, N_{C, \alpha , \beta}$ and $a$ such that the inequalities
\begin{eqnarray}
& & || (e^{T_1t})^\times \Pi_c \xi ||^*_{\beta, n} \leq L_{C,\alpha , \beta} || \xi ||^*_{\beta, n}, \label{7-41} \\
& & || (e^{T_1t})^\times \Pi_s \xi ||^*_{\beta, n} \leq M_{C,\alpha , \beta} e^{-at} || \xi ||^*_{\beta, n},\label{7-42}
\end{eqnarray}
hold for $\xi \in \overline{i(W_{C, \alpha })}$, and the inequality
\begin{eqnarray}
|| (e^{T_jt})^\times \xi ||^*_{\beta, n} \leq N_{C,\alpha , \beta} e^{-jat} || \xi ||^*_{\beta, n},\,\, j= 2,3, \cdots ,\label{7-43}
\end{eqnarray}
holds for $\xi \in \overline{i(W_{C, \alpha })}$.
\\[0.2cm]
\textbf{Proof.}
Let $\lambda _0 ( =0)$ be the resonance pole on the imaginary axis.
For $\xi \in i(W_{C, \alpha })$, $(e^{T_1t})^\times \Pi_c \xi$ is calculated as
\begin{eqnarray*}
(e^{T_1t})^\times \Pi_c \xi
= \frac{K}{2D_0}\langle  \mu_0 \,|\, \xi^* \rangle \cdot (e^{T_1t})^\times \mu_0  
= \frac{K}{2D_0} \langle  \mu_0 \,|\,\xi^* \rangle \cdot e^{\lambda_0t} \mu_0.
\end{eqnarray*}
Since $\lambda _0=0$, we obtain
\begin{eqnarray*}
|| \, (e^{T_1t})^\times \Pi_c \xi \, ||^*_{\beta, n}
\leq  || \, \frac{K}{2D_0}  
      \langle \mu_0 \,|\,\xi^* \rangle \cdot \mu_0 \, ||^*_{\beta, n}
=  ||\, \Pi_0 \xi \,||^*_{\beta, n}
\leq D_{C, \alpha ,\beta, 0} || \xi ||^*_{\beta ,n}.
\end{eqnarray*}
This provides Eq.(\ref{7-41}) for $\xi \in i(W_{C, \alpha })$ with $L_{C, \alpha , \beta}= D_{C, \alpha , \beta, 0}$.
The result is continuously extended for $\xi \in \overline{i(W_{C, \alpha })}$.
The proofs of Eq.(\ref{7-42}) and (\ref{7-43}) are the same as that of Prop 5.14 with the aid of the fact that $\Pi_s = id - \Pi_c$
is continuous on $i(W_{C, \alpha })$. 
\hfill $\blacksquare$
\\

Note that $\overline{i(V)}$ is a closed subspace of $\mathrm{Exp}_-'$.
If $\xi \in V$, $i(\xi)$ satisfies inequalities (\ref{thm5-18}),(\ref{7-41}),(\ref{7-42}),(\ref{7-43}),
in which the constants depend only on $\beta$ because $\overline{i(V)} \subset \overline{i(W_{3,0})}$.
Since $\mathbf{E}_c \subset \overline{i(V)}$, the generalized eigenfunctions in $\mathbf{E}_c$ also satisfy
the inequalities with the same constants.
The space $\overline{i(V)}$ has all properties for developing a bifurcation theory:
it is a metric space including all solutions of the Kuramoto model and the generalized center subspace.
The projection $\Pi_c$ is continuous on $\overline{i(V)}$.
The semigroup $(e^{T_1t})^\times$ admits the spectral decomposition on it,
and norms of the semigroups $(e^{T_jt})^\times$ satisfy the appropriate inequalities.
By using these properties, we will prove the existence of center manifolds in Section 7.


\subsection{Spectral theory on $(H_-, L^2(\mathbf{R}, g(\omega )d\omega ),H_-')$}

In this subsection, we suppose that $g(\omega )$ is a rational function.
Let $H_+$ be a Banach space of bounded holomorphic functions on the real axis and the upper half plane (see Sec.4.3)
and $H_- = \{ \phi^* \, | \, \phi \in H_+\}$.
In this case, $H_-$ is not a dense subspace of $L^2(\mathbf{R}, g(\omega )d\omega )$, and thus the triplet
$(H_-, L^2(\mathbf{R}, g(\omega )d\omega ),H_-')$ is a degenerate rigged Hilbert space.
\\[0.2cm]
\textbf{Proposition 5.16.} \, The canonical inclusion $i : H_+ \to H_-'$ is a finite dimensional operator ;
that is, $i(H_+)\subset H_-'$ is a finite dimensional vector space.
\\[0.2cm]
\textbf{Proof.} \, By the definition,
\begin{eqnarray*}
i(\psi)(\phi^*) = \langle \psi \,|\, \phi^*  \rangle
 = (\psi, \phi^* ) = \int_{\mathbf{R}} \! \phi (\omega )\psi (\omega ) g(\omega )d\omega , 
\end{eqnarray*}
for $\phi, \psi \in H_+$.
Let $z_1, \cdots  ,z_n$ be poles of $g(\omega )$ on the upper half plane.
By the residue theorem, we obtain
\begin{eqnarray*}
& & \int_{\mathbf{R}} \! \phi (\omega )\psi (\omega ) g(\omega )d\omega
   + \lim_{r\to \infty } 
     \int^{\pi}_{0} \! \phi (re^{\sqrt{-1}\theta }) \psi (re^{\sqrt{-1}\theta }) g(re^{\sqrt{-1}\theta })
                  \sqrt{-1}re^{\sqrt{-1}\theta }d\theta  \\
& = & 2\pi \sqrt{-1} \sum^n_{j=1} \mathrm{Res}(z_j), 
\end{eqnarray*}
where, $\mathrm{Res}(z_j)$ denotes the residue of $\phi (\omega )\psi (\omega )g(\omega )$ at $z_j$.
Since $g(\omega )$ is a rational function which is integrable on the real axis, the degree of the denominator
is at least two greater than the degree of the numerator : $g(\omega ) \sim O(1/|\omega |^2)$ as $|\omega | \to \infty$.
Since $\phi, \psi \in H_+$ is bounded on the upper half plane, we obtain
\begin{equation}
\langle \psi \,|\, \phi^*  \rangle =2\pi \sqrt{-1} \sum^n_{j=1}\mathrm{Res}(z_j),
\end{equation}
as $r\to \infty$. This means that the action of $i(\psi ) \in i(H_+)$ on $H_-$ is determined by
the values of $\psi (\omega )$ and its derivatives at $z_1, \cdots  ,z_n$.
In particular, if the denominator of $g$ is of degree $M$, then $i(H_+) \simeq \mathbf{C}^M$. 
\hfill $\blacksquare$
\\[-0.2cm]

Since $i(H_+)$ is of finite dimensional, the semigroup $(e^{T_1t})^\times$ restricted to $i(H_+)$ is a finite dimensional
operator. This is the reason that Eq.(\ref{thm4.6}) consists of a finite sum.
In what follows, we suppose that all resonance poles are simple roots of Eq.(\ref{semi7}).
Then Eq.(\ref{thm4.6}) is rewritten as
\begin{equation}
(e^{T_1t}\phi, \psi^* ) = \langle (e^{T_1t})^\times\phi \,|\,  \psi^* \rangle
 = \sum^M_{n=0} \frac{K}{2D_n}e^{\lambda _nt} \langle \mu_n \,|\, \phi^* \rangle \langle \mu_n \,|\, \psi^* \rangle,
\end{equation}
where definitions of $D_n$ and $\mu_n$ are the same as those in previous sections.
Now we have obtained the following theorem.
\\[0.2cm]
\textbf{Theorem 5.17.} \, For any $\psi \in H_+$, the equalities
\begin{eqnarray}
(e^{T_1t})^\times \psi
  &=& \sum^M_{n=0} \frac{K}{2D_n}e^{\lambda_nt}\langle \mu_n  \,|\, \psi^* \rangle \cdot \mu_n, \\
i(\psi )
  &=& \sum^M_{n=0} \frac{K}{2D_n}\langle \mu_n  \,|\, \psi^* \rangle \cdot \mu_n,
\end{eqnarray}
hold. In particular, a system of generalized eigenfunctions $\{ \mu_n \}^M_{n=0}$ forms a base of $i(H_+)$.
\\[-0.2cm]

The projection $\Pi_n : i(H_+) \to \mathrm{span }\{ \mu_n\} \subset H_-'$ is defined to be
\begin{equation}
\Pi_n \psi = \frac{K}{2D_n}\langle \mu_n \,|\, \psi^* \rangle \cdot \mu_n,
\quad n = 0, \cdots ,M
\end{equation}
as before. Since $i(H_+)$ is a finite dimensional vector space, $\Pi_n$ is continuous on the whole space.
Note that solutions $Z_1, Z_2, \cdots $ of the Kuramoto model are included in $H_+$ 
(we have proved that $Z_j\in V_{1,0, \infty}$ in Thm.5.10 (iii)).
Thus the bifurcation problem of the Kuramoto model is reduced to the bifurcation theory on a finite dimensional space,
and the usual center manifold theory is applicable.


\section{Nonlinear stability}

Before going to the bifurcation theory, let us consider the nonlinear stability of the de-synchronous state.
In Sec.4 and Sec.5, we proved that the order parameter $\eta (t) \equiv 0$ is linearly stable when $0 < K < K_c$;
that is, the asymptotic stability of $\eta (t) \equiv 0$ is proved for the linearized system (\ref{4-1}).
For a system on an infinite dimensional space, in general, the linear stability does not imply the nonlinear stability.
Infinitesimally small nonlinear terms may change the stability of fixed points.
In this section, we show that the de-synchronous state $Z_j(t) \equiv 0\, (j=1,2,\cdots )$ 
(which corresponds to $\rho_t \equiv 1/2\pi$) is locally stable with
respect to a suitable topology when $0 < K < K_c$. In particular, the order parameter proves to decay to zero
as $t\to \infty$ without neglecting the nonlinear terms.

Recall that the continuous model (\ref{conti}) is rewritten as Eqs.(\ref{4-0}),(\ref{4-0b}) by putting
$Z_j(t, \omega ) = \int^{2\pi}_{0} \! e^{\sqrt{-1}j \theta } \rho_t(\theta , \omega ) d\theta $ with the 
initial condition
\begin{equation}
Z_j(0,\omega ) = \int^{2\pi}_{0} \! e^{\sqrt{-1}j \theta }h(\theta ) d\theta  := h_j\in \mathbf{C}.
\label{6-1}
\end{equation}
We need not suppose that $h(\theta )$ is a usual function.
It may be a probability measure on $S^1$.
We have proved that solutions $Z_j(t, \cdot)$ are included in the set $V_{1,0} \subset \mathrm{Exp}_+$. 
By using the canonical inclusion $i : L^2(\mathbf{R}, g(\omega )d\omega ) \to \mathrm{Exp}_-'$,
, Eqs.(\ref{4-0}),(\ref{4-0b}) are rewritten as a system of evolution equations
on $\prod^\infty_{j=1}\mathrm{Exp}_-'$ of the form
\begin{eqnarray}
\left\{ \begin{array}{ll}
\displaystyle \frac{d}{dt} Z_1
  = T_1^\times  Z_1
      - \frac{K}{2} \langle \overline{Z_1 \,|\, P_0} \rangle Z_2, \\[0.4cm]
\displaystyle \frac{d}{dt} Z_j  
  = T_j^\times Z_j + \frac{jK}{2} \left( \,
     \langle Z_1 \,|\, P_0 \rangle  Z_{j-1} 
      - \langle \overline{Z_1 \,|\, P_0 } \rangle Z_{j+1} \right), \quad j = 2,3,\cdots , \\[0.4cm]
 Z_j(0, \,\cdot \,) = h_j P_0,   
\end{array} \right.
\label{6-2}
\end{eqnarray}
where $Z_j$ is an abbreviation for $i(Z_j) \in i(\mathrm{Exp}_+) \subset \mathrm{Exp}_-'$.
Linear operators $T_j$ are defined to be
\begin{equation}
T_1\phi (\omega ) = (\sqrt{-1}\mathcal{M} + \frac{K}{2}\mathcal{P})\phi (\omega )
   = \sqrt{-1}\omega \phi (\omega ) + \frac{K}{2} \langle  \phi \,|\, P_0 \rangle P_0(\omega ),
\label{t1}
\end{equation}
and
\begin{equation}
T_j \phi (\omega ) = \sqrt{-1}j\mathcal{M}\phi (\omega ) = \sqrt{-1}j\omega \phi (\omega ),
\end{equation}
for $j = 2,3,\cdots $, and $T_j^\times$ are their dual operators.
The main theorem in this section is stated as follows.
\\[0.2cm]
\textbf{Theorem 6.1 (local stability of the de-synchronous state).}\, 
Suppose that $g (\omega )$ is the Gaussian and $0<K<K_c$.
Then, there exists a positive constant $\delta_\beta$ such that if the initial condition $h(\theta )$ of the initial value
problem (\ref{conti}) satisfies
\begin{equation}
|h_j| = \left| \int^{2\pi}_{0} \! e^{j\sqrt{-1}\theta } h(\theta ) d\theta  \right| \leq  \delta_\beta,
 \quad j=1,2,\cdots ,
\label{6-11}
\end{equation}
then the quantities 
\begin{eqnarray*}
(Z_j, \phi) = \int^{2\pi}_{0} \! \int_{\mathbf{R}} \! e^{\sqrt{-1}j \theta } \phi (\omega )
         \rho_t(\theta , \omega ) d\omega d\theta   
\end{eqnarray*}
tend to zero as $t\to \infty$ for every $\phi \in \mathrm{Exp}_+(\beta )$ uniformly in $j = 1,2, \cdots $. 
In particular, the order parameter $\eta (t) = (Z_1, P_0)$ tends to zero as $t\to \infty$.

This theorem means that the trivial solution $Z_j \equiv 0$ of (\ref{6-2}) is locally stable with
respect to the weak dual topology on $\mathrm{Exp}_-'$.
In general, $\delta _\beta \to 0$ as $\beta \to \infty$.
One of the reasons is that the norm $|| \cdot ||_{\beta,n}^*$ goes to infinity as $\beta \to \infty$.
For the case $g(\omega )$ is a rational function,
we can show the same statement : $(Z_j, \phi)$ tends to zero as $t\to \infty$ for every $\phi \in H_+$ 
if the initial condition satisfies (\ref{6-11}), in which $\delta _\beta$ is independent of $\beta$. 
\\[0.2cm]
\textbf{Proof of Thm.6.1.}\, Since we have Prop.5.14, the proof is done in a similar manner to the proof of the local stability 
of fixed points of finite dimensional systems. Eq.(\ref{6-2}) provides
\begin{equation}
\left\{ \begin{array}{l}
\displaystyle Z_1(t, \cdot) = (e^{T_1(t-t_0)})^\times  Z_1(t_0, \cdot)  
   - \frac{K}{2} \int^t_{t_0} \! \langle \overline{Z_1(s, \cdot) \,|\, P_0} \rangle 
             (e^{T_1(t-s)})^\times  Z_2(s, \cdot) ds, \\[0.4cm]
\displaystyle  Z_j(t, \cdot) = (e^{T_j(t-t_0)})^\times Z_j(t_0, \cdot) 
   + \frac{jK}{2} \int^t_{t_0} \! \Bigl( \langle Z_1(s, \cdot) \,|\,P_0  \rangle 
             (e^{T_j(t-s)})^\times Z_{j-1}^*(s, \cdot) \\[0.2cm]
\displaystyle  \qquad \qquad \qquad \qquad \qquad \qquad \qquad
           - \langle \overline{Z_1(s, \cdot) \,|\, P_0 } \rangle 
             (e^{T_j(t-s)})^\times Z_{j+1}(s, \cdot)\Bigr) ds,
\end{array} \right.
\end{equation}
for $0 \leq t_0 < t$.
Since $ Z_j \in i(V_{1,0}) \subset i(W_{3,0})$ for every $t>0$ and $j=1,2,\cdots $, 
Prop.5.14 is applied to show that there exists $M_{3,0,\beta} = M_{\beta} > 0$ such that 
\begin{equation}
\left\{ \begin{array}{l}
\displaystyle 
||\, Z_1(t, \cdot) ||^*_{\beta, n} \leq M_\beta e^{-a(t-t_0)} || \, Z_1(t_0, \cdot) ||^*_{\beta, n} 
   + \frac{K}{2} \int^t_{t_0} \! M_\beta e^{-a(t-s)} ||\, Z_1(s, \cdot) ||^*_{\beta, n} \cdot 
             ||\, Z_2(s, \cdot) ||^*_{\beta, n} ds, \\[0.4cm]
\displaystyle 
   || \, Z_j(t, \cdot) ||^*_{\beta, n} \leq M_\beta e^{-ja(t-t_0)} || \, Z_j(t_0, \cdot) ||^*_{\beta, n} 
   + \frac{jK}{2} \int^t_{t_0} \! M_\beta e^{-ja (t-s)} || \, Z_1(s, \cdot) ||^*_{\beta, n} \cdot \\[0.2cm]
\displaystyle   \qquad \qquad \qquad \qquad \qquad \qquad \qquad \qquad \qquad
             \Bigl( || \, Z_{j-1}(s, \cdot) ||^*_{\beta, n}
                            + || \, Z_{j+1}(s, \cdot) ||^*_{\beta, n} \Bigr) ds.
\end{array} \right.
\label{6-13}
\end{equation}
Take a small constant $\delta _\beta > 0$ such that $h_j = \int^{2\pi}_{0}\! e^{j\sqrt{-1}\theta }h(\theta )d\theta $
satisfies Eq.(\ref{6-11}).
Let us show that there exists $N_\beta > 0$ such that
\begin{equation}
|| \, Z_j(t, \cdot) ||^*_{\beta, n} \leq \delta _\beta N_\beta
\label{6-15}
\end{equation}
holds for any $t>t_0$ and $j=1,2,\cdots $.
Indeed, Eq.(\ref{xi2}) shows that $|Z_j(t, \omega )| \leq \delta _\beta$ holds for any $t, \omega $ and $j$.
Hence, 
\begin{eqnarray*}
|| Z_j ||^*_{\beta, n} 
& = & \sup _{|| \psi ||_{\beta, n} = 1} |\langle Z_j \,|\, \psi^* \rangle| 
 = \sup _{|| \psi ||_{\beta, n} = 1} \left| \int_{\mathbf{R}}\! Z_j(t, \omega )\psi (\omega )g(\omega )d\omega \right| \\
&\leq & \delta _\beta \int_{\mathbf{R}}\! e^{\beta |\omega |} g(\omega )d\omega .
\end{eqnarray*}
Therefore, putting $N_\beta := \int_{\mathbf{R}}\! e^{\beta |\omega |} g(\omega )d\omega $ proves Eq.(\ref{6-15}).
Then, the first equation of (\ref{6-13}) gives
\begin{equation}
||\,  Z_1(t, \cdot) ||^*_{\beta, n}
\leq \delta _\beta M_\beta N_{\beta} e^{-a(t-t_0)} + \frac{\delta _\beta KM_\beta N_\beta }{2} \int^t_{t_0} \! e^{-a(t-s)} 
||\,  Z_1(s, \cdot) ||^*_{\beta, n} ds,
\end{equation}
for $t > t_0$. Now the Gronwall inequality proves
\begin{equation}
||\, Z_1(t, \cdot) ||^*_{\beta, n}
\leq \delta _\beta M_\beta N_{\beta} e^{(\delta _\beta KM_\beta N_\beta /2 - a)(t-t_0)}.
\end{equation}
Since $M_\beta$ and $N_\beta$ are independent of the choice of $\delta _\beta$, by taking $\delta _\beta$ sufficiently small,
this quantity proves to tend to zero as $t\to \infty$.
Substituting it into the second equation of (\ref{6-13}), we obtain
\begin{equation}
||\, Z_j(t, \cdot) ||^*_{\beta, n} \leq 
\delta _\beta M_\beta N_{\beta}e^{-ja(t-t_0)} + jK (\delta _\beta M_\beta N_\beta )^2
 \int^t_{t_0} \! e^{-ja(t-s)} e^{(\delta _\beta KM_\beta N_\beta /2 - a)(s-t_0)}ds
\end{equation}
for $j = 2,3, \cdots $.
It is easy to verify that the right hand side tends to zero as $t\to \infty$ uniformly in $j$.

Now we have proved that if the initial condition satisfies (\ref{6-11}) for each $\beta$, then
$||\, Z_j(t, \cdot) ||^*_{\beta, n}$ decays to zero as $t\to \infty$
for every $j$ and $n$. By the definition of the norm $|| \cdot ||^*_{\beta, n}$, this means that
$(Z_j(t, \cdot), \phi) \to 0$ as $t\to \infty$ for every $\phi \in \mathrm{Exp}_+(\beta)$.
\hfill $\blacksquare$


\section{Bifurcation theory}

Now we are in a position to investigate bifurcation of the Kuramoto model by using the center manifold reduction.
Our strategy to detect bifurcation is that we use the space of functionals $\mathrm{Exp}_-'$ instead of the spaces of usual functions
$\mathrm{Exp}_+$ or $L^2 (\mathbf{R}, g(\omega )d\omega )$ because the linear operator $T_1$ admits the spectral 
decomposition on $\mathrm{Exp}_-'$ consisting of a countable number of eigenfunctions, 
while the spectral decomposition on $L^2 (\mathbf{R}, g(\omega )d\omega )$
involves the continuous spectrum on the imaginary axis;
that is, a center manifold on $L^2 (\mathbf{R}, g(\omega )d\omega )$ is an infinite dimensional manifold.
To avoid such a difficulty, we will seek a center manifold on $\mathrm{Exp}_-'$.
At first, we have to prove the existence of center manifolds.
Standard results of the existence of center manifolds (see \cite{Bat, Che, Kri, Van}) are not applicable
to our system because the space $\mathrm{Exp}_-'$ is not a Banach space and the projection $\Pi_c$ to the 
center subspace is continuous only on a subspace of $\mathrm{Exp}_-'$.
Thus in Sec.7.1, the existence theorem of center manifolds for our system and a strategy for proving it are given.
The proof of the theorem is given in Sec.7.2 to 7.4.
In Sec.7.5, the dynamics on the center manifold is derived and the Kuramoto's conjecture is solved.
Readers who are interested in a practical method for obtaining  a bifurcation structure can skip 
Sec.7.1 to 7.4 and go to Sec.7.5.
Throughout this section, we suppose that $g(\omega )$ is the Gaussian.
Existence of center manifolds for the case that $g(\omega )$ is a rational function is trivial because
the phase space $i(H_+)$ is a finite dimensional vector space.


\subsection{Center manifold theorem}

Let $i(\mathcal{F})$ be a certain metric subspace of the 
product space $\prod^\infty_{k=1} \mathrm{Exp}_-'$ with a distance $d_\infty$, and $\overline{i(\mathcal{F})}$ its closure.
These spaces and the metric $d _\infty$ will be introduced in Sec.7.2 and Sec.7.3. 
Let $\Phi_t$ be the semiflow on $\overline{i(\mathcal{F})}$ generated by the system (\ref{6-2}).
For the generalized center subspace $\mathbf{E}_c$ of $T_1$ defined by (\ref{center}), put
\begin{equation}
\hat{\mathbf{E}}_c = \mathbf{E}_c \times \{0 \} \times \{0\} \times \cdots \subset \prod^\infty_{k=1} \mathrm{Exp}_-'.
\end{equation}
Let $\hat{\mathbf{E}}^\bot_c = \mathbf{E}^\bot_c \times \mathrm{Exp}_-' \times \mathrm{Exp}_-' \times \cdots $ be the complement of $\hat{\mathbf{E}}_c$.
The existence theorem of center manifolds is stated as follows.
\\[0.2cm]
\textbf{Theorem 7.1.} 
There exist a positive number $\varepsilon _0$
and an open set $U$ of the origin in $\overline{i(\mathcal{F})} \subset \prod^\infty_{k=1} \mathrm{Exp}_-'$ such that when 
$| K- K_c |< \varepsilon _0$, the following holds:
\\[0.2cm]
(I) (center manifold). There exists a $C^1$ mapping 
$\hat{q} : \hat{\mathbf{E}}_c \to \hat{\mathbf{E}}^\bot_c \cap \overline{i(\mathcal{F})}$
such that the one dimensional $C^1$ manifold defined to be 
\begin{equation}
W^c_{loc} = \{y + \hat{q}(y) \, | \, y\in \hat{\mathbf{E}}_c \} \cap U
\end{equation}
is $\Phi_t$-invariant (that is, $\Phi_t(W^c_{loc}) \cap U \subset W^c_{loc}$).
This is called the \textit{local center manifold}.
The mapping $\hat{q}$ is also $C^1$ with respect to the parameter $\varepsilon  := K-K_c$, and 
$\hat{q}(y) \sim O(y^2, \varepsilon y, \varepsilon ^2)$ as $y, \varepsilon \to 0$
\\[0.2cm]
(II) (negative semi-orbit). 
For every $\xi_0 \in W^c_{loc}$, there exists a function $u : (-\infty, 0] \to \overline{i(\mathcal{F})}$
such that $u(0) = \xi_0$ and $\Phi _t (u(s)) = u(t+s)$ when $t\geq 0,\, s\leq -t$.
Such a $u(t)$ is called a \textit{negative semi-orbit} of (\ref{6-2}).
As long as $u(t)\in U$, $u(t) \in W^c_{loc}$. In this case, there exist $C_1>0$ and a small number $b >0$ such that
\begin{equation}
d_\infty (u(t), 0) \leq C_1 e^{b t}.
\label{7-3}
\end{equation} 
\\
(III) (invariant foliation).
There exists a family of manifolds $\{M_\xi\}_{\xi \in W^c_{loc}} \subset U$, parameterized by $\xi \in W^c_{loc}$, satisfying that
\\
\quad (i) $ M_\xi \cap W^c_{loc} = \{\xi\},\, \bigcup_{\xi \in W^c_{loc}} M_\xi = U$,
and $ M_\xi \cap M_{\xi'} = \emptyset$ if $\xi \neq \xi'$.
\\
\quad (ii) when $\Phi_t (\xi) \in U$, $\Phi_t(M_\xi) \cap U \subset M_{\Phi_t(\xi)}$.
\\
\quad (iii)
\\[-0.8cm]
\begin{eqnarray*}
M_\xi &=& \{ u\in \overline{i(\mathcal{F})} \cap U \, \Bigl| \, \\
& & \quad \mathrm{there \,\, exist}\,\, a>b \,\, \mathrm{and}\,\, C_2 >0 \,\, \mathrm{such \,\,that} \,\,
d_\infty(\Phi_t(u) , \Phi_t(\xi) ) \leq C_2 e^{-at}\}.
\end{eqnarray*}

Part (III) of the theorem means that $W^c_{loc}$ is attracting with the decay rate $e^{-at}$, where the constant $a$
is the same as that in Prop.5.15.
Further, (III)-(iii) means that the semiflow near $W^c_{loc}$ is eventually well approximated by the semiflow on $W^c_{loc}$ if $t>0$ is large.
In particular, if (\ref{6-2}) has an attractor $N$ near the origin, $N$ is included in $W^c_{loc}$.
Since the topology induced by the metric coincides with the strong and the weak dual topologies on any bounded set,
$N$ is attracting for both of the strong and the weak dual topologies.
Due to the spectral decomposition, any element $ Z_1 \in i(\mathrm{Exp}_+)$ is decomposed as
\begin{equation}
Z_1 = \alpha \mu_0 + Y_1, \quad \alpha \in \mathbf{C},\, Y_1 \in \mathbf{E}_c^\bot,
\label{7-4}
\end{equation}
where $\mu_0 \in \mathbf{E}_c$ is a generalized eigenfunction associated with the resonance pole $\lambda _0=0$.
Then, part (I) of the theorem means that if $(Z_1, Z_2, \cdots  ) \in W^c_{loc}$,
\begin{equation}
Y_1 \sim O(\alpha ^2, \varepsilon \alpha , \varepsilon ^2), \quad 
     Z_k  \sim O(\alpha ^2, \varepsilon \alpha , \varepsilon ^2), \,k= 2,3, \cdots ,
\label{7-5}
\end{equation}
as $\varepsilon , \alpha \to 0$.
Substituting Eq.(\ref{7-4}) into the system (\ref{6-2}) with the condition (\ref{7-5}),
we can obtain the expression of $\hat{q}(y)$ as a function of $\varepsilon , \alpha$.
The dynamics on $W^c_{loc}$ is realized by an ordinary differential equation of $\alpha$:
\begin{equation}
\frac{d\alpha }{dt} = f(\varepsilon ; \alpha ).
\label{7-6}
\end{equation}
If (\ref{6-2}) has an attractor $N$ near the origin, $N$ is an attractor of the system (\ref{7-6}).
In this manner, (\ref{6-2}) is reduced to a finite dimensional problem.
Such a technique to investigate bifurcation is called the \textit{center manifold reduction}.
Part (II) implies that the center manifold is characterized by the property that 
the dynamics on it is sufficiently slow.

Although we prove the existence of $W^c_{loc}$ in $\overline{i(\mathcal{F})}$, from a physical viewpoint, especially
we are interested in an initial condition of the form $Z_j(0)= h_j  P_0$
(which corresponds to the initial condition of the form $\rho_0 (\theta , \omega ) = h(\theta )$ for the system (\ref{conti})).
Then $Z_j(t, \cdot ) \in i(V_{1,0})$ (Thm.5.10 (iii)).
This means that an attractor of (\ref{6-2}) which is reached from the initial condition 
$Z_j(0)= h_j  P_0$ is included in $W^c_{loc} \cap \prod^\infty_{k=1}\overline{i(V_{1,0})}$.

Sec.7.2 to 7.4 are devoted to prove Thm.7.1.
It is well known that a global center manifold uniquely exists only when a Lipschitz constant of nonlinear terms of a system
is sufficiently small.
Thus in Sec.7.2, we consider a perturbed system of (\ref{6-2}) so that its Lipschitz constant becomes sufficiently small,
while it coincides with the original system in the vicinity of the origin.
Because of the perturbation, a solution may fall out of $i(V_{1,0})$ and go into a larger space.
Thus we will introduce the space $\mathcal{F}$, and show that solutions 
$(Z_1, Z_2, \cdots )$ of the perturbed system are included in $\mathcal{F}$.
We will prove in Sec.7.3 that the perturbed system generates a smooth flow
to prove that the center manifold is smooth.
Once we obtain the existence of a proper phase space, a spectral decomposition of the linear operator, 
estimates of norms of the semigroups and a smooth flow whose Lipschitz constant of nonlinear terms is sufficiently small,
then the existence of the center manifold is proved in the usual way
with a slight modification.
We demonstrate it in Sec.7.4.
In Sec.7.5, we perform the center manifold reduction: 
an equation of $\alpha $ is obtained and investigated.
The order parameter $\eta (t)$ is defined as $\eta (t) = (Z_1, P_0) = \langle Z_1 \,|\, P_0 \rangle$.
On the center manifold, it is written as
\begin{equation}
\langle Z_1 \,|\, P_0 \rangle = \alpha (t) \langle \mu_0 \,|\, P_0 \rangle
 + \langle Y_1 \,|\, P_0 \rangle
   = \frac{2}{K_c} \alpha (t) + O(\alpha ^2, \varepsilon \alpha , \varepsilon ^2),
\end{equation}
where we use $\langle \mu_0  \,|\, P_0 \rangle = 2/K_c$, which follows from the definition of resonance poles.
Therefore, if a bifurcation diagram of (\ref{7-6}) is obtained, a bifurcation diagram of the order parameter
is also obtained. In this way, the Kuramoto's conjecture will be proved in Sec.7.5.


\subsection{Phase space of the perturbed system}

Recall that the trivial solution $Z_j(t) \equiv 0\, (j=1,2,\cdots )$ which corresponds to $\rho_t \equiv 1/2\pi$
is called the de-synchronous state (Sec.3). Since we are interested in bifurcations from $\rho_t \equiv 1/2\pi$
at $K = K_c$, put $\rho_t = 1/2\pi + \hat{\rho}_t$ and $K = K_c + \varepsilon $.
Then, Eq.(\ref{conti}) is rewritten as
\begin{eqnarray}
& & \frac{\partial \hat{\rho}_t}{\partial t} +
\frac{\partial }{\partial \theta }
\Bigl( \omega \hat{\rho}_t + \frac{K_c}{4\pi \sqrt{-1}}(\eta (t) e^{-\sqrt{-1}\theta }-\overline{\eta (t)}e^{\sqrt{-1}\theta }) \nonumber \\ 
& &  + \frac{\varepsilon }{4\pi \sqrt{-1}} (\eta (t) e^{-\sqrt{-1}\theta }-\overline{\eta (t)} e^{\sqrt{-1}\theta }) 
    + \frac{K}{2\sqrt{-1}}(\eta (t) e^{-\sqrt{-1}\theta }-\overline{\eta (t)}e^{\sqrt{-1}\theta })\hat{\rho}_t \Bigr) = 0, \quad
\label{7-8}
\end{eqnarray}
where
\begin{eqnarray*}
\eta (t) = \int_{\mathbf{R}} \! \int^{2\pi}_{0} \! e^{\sqrt{-1}\theta } \rho_t(\theta ,\omega )d\theta d\omega 
 = \int_{\mathbf{R}} \! \int^{2\pi}_{0} \! e^{\sqrt{-1}\theta } \hat{\rho}_t(\theta ,\omega )d\theta d\omega.
\end{eqnarray*}
An initial condition $\hat{\rho}_0 = \hat{h}(\theta , \omega )$ is a real-valued measure (signed measure) on $S^1$
parameterized by $\omega \in \mathbf{R}$ satisfying $\int^{2\pi}_{0} \! \hat{h}(\theta , \omega ) d\theta =0$.
The first step to prove the existence of center manifolds is to localize the nonlinear term so that the Lipschitz constant
of the nonlinear term becomes sufficiently small.
For this purpose, let $\hat{\chi} : [0, \infty) \to [0,1]$ be a $C^\infty$ function, and consider 
the perturbed continuous model of the form
\begin{eqnarray}
& & \frac{\partial \hat{\rho}_t}{\partial t} +
\frac{\partial }{\partial \theta }
\Bigl( \omega \hat{\rho}_t + \frac{K_c}{4\pi \sqrt{-1}}(\eta (t) e^{-\sqrt{-1}\theta }-\overline{\eta (t)}e^{\sqrt{-1}\theta }) \nonumber \\ 
& &  + \frac{\varepsilon }{4\pi \sqrt{-1}} (\eta (t) e^{-\sqrt{-1}\theta }-\overline{\eta (t)} e^{\sqrt{-1}\theta }) \hat{\chi}(t)
 + \frac{K}{2\sqrt{-1}}(\eta (t) e^{-\sqrt{-1}\theta }-\overline{\eta (t)}e^{\sqrt{-1}\theta })\hat{\chi}(t)\hat{\rho}_t \Bigr) = 0. 
\quad \quad \quad
\label{7-9}
\end{eqnarray}
If we put 
\begin{equation}
\hat{Z}_j(t,\omega ) = \int^{2\pi}_{0} \! e^{\sqrt{-1}j \theta } \hat{\rho}_t (\theta , \omega ) d\theta , 
\label{7-10}
\end{equation}
Eq.(\ref{7-9}) yields a system of equations 
\begin{equation}
\left\{ \begin{array}{l}
\displaystyle \frac{d\hat{Z}_1}{dt}
   = \sqrt{-1}\omega \hat{Z}_1 + \frac{K_c}{2}\eta (t) + 
 \left( \frac{\varepsilon }{2}\eta (t)- \frac{K}{2}\overline{\eta (t)}\hat{Z}_2 \right) \hat{\chi}(t), \\[0.4cm]
\displaystyle \frac{d\hat{Z}_j}{dt}
   = j\sqrt{-1}\omega \hat{Z}_j
+ \frac{jK}{2}\left( \eta (t)\hat{Z}_{j-1} -\overline{\eta (t)}\hat{Z}_{j+1} \right) \hat{\chi}(t),\,\, j=2, 3,\cdots .
\end{array} \right.
\label{7-11}
\end{equation}
If $\hat{\chi} \equiv 1$, this coincides with the original system (\ref{4-0}),(\ref{4-0b}).
If $\hat{\chi}(t)$ is sufficiently small, this perturbation makes the Lipschitz constant of the nonlinear terms
of (\ref{4-0}),(\ref{4-0b}) sufficiently small
(Note that when proving the existence of center manifolds, the bifurcation parameter $\varepsilon $ is regarded as
a dependent variable. Thus $\varepsilon \eta (t) = \varepsilon (\hat{Z}_1, P_0)$ is regarded as a nonlinear term).
A concrete definition of $\hat{\chi}$ will be specified in Sec.7.3.
Note that $\hat{Z}_0 \equiv 0$ because of $\int^{2\pi}_{0} \! \hat{h}(\theta , \omega ) d\theta =0$.

Eq.(\ref{7-9}) is integrated by using the characteristic curve method.
The characteristic curve $x = x(t,s; \theta ,\omega )$ is defined as a solution of the equation
\begin{equation}
\frac{dx}{dt} = \omega + \frac{K}{2\sqrt{-1}}\left( \eta (t) e^{-\sqrt{-1}x} - \overline{\eta (t)} e^{\sqrt{-1}x}\right) \hat{\chi}(t),
\label{7-13}
\end{equation}
satisfying the initial condition $x(s,s;\theta ,\omega ) = \theta $ at an initial time $s$.
Along the characteristic curve, (\ref{7-9}) is integrated to yield
\begin{eqnarray}
\hat{\rho}_t(\theta , \omega ) &=& 
\hat{h}(x(0,t;\theta ,\omega ),\omega ) \exp \Bigl[ \frac{K}{2} \int^t_{0} \! \left( \eta (s) 
e^{-\sqrt{-1}x(s,t;\theta ,\omega )} + \overline{\eta(s)}e^{\sqrt{-1}x(s,t;\theta ,\omega )} \right) \hat{\chi}(s) ds \Bigr] \nonumber \\
& & + \int^t_{0} \! \exp \Bigl[ \frac{K}{2} \int^t_{s} \! \left( \eta (\tau) e^{-\sqrt{-1}x(\tau,t;\theta ,\omega )}
     + \overline{\eta(\tau)}e^{\sqrt{-1}x(\tau,t;\theta ,\omega )} \right) \hat{\chi}(\tau) d\tau \Bigr] \times \nonumber \\
& & \qquad \qquad \left( \frac{K_c + \varepsilon \hat{\chi}(s)}{4\pi} \right) 
\left( \eta (s) e^{-\sqrt{-1}x(s,t;\theta ,\omega )} + \overline{\eta(s)}e^{\sqrt{-1}x(s,t;\theta ,\omega )}\right) ds.
\label{7-14}
\end{eqnarray}
Once $x(t,s;\theta , \omega )$ and $\eta (t)$ are determined, this $\hat{\rho}_t$ gives a weak solution of (\ref{7-9}).
The existence of solutions of (\ref{7-11}) follows from that of the integro-ODE (\ref{7-13}) and (\ref{7-14}),
which will be proved by the standard iteration method (see also Prop.7.3).

In Sec.5, we have proved that a solution of (\ref{4-0}),(\ref{4-0b}) is included in $V_{1,0, \infty}$.
This property may break down because of the perturbation $\hat{\chi}$.
Thus we define an appropriate phase space for (\ref{7-11}) and prove the existence of the flow on it.
Fix a finite number $p > 0$ and let $V_{1,0,p}$ be a set defined by Eq.(\ref{Vca}).
$V = V_p$ is defined by 
\begin{eqnarray*}
V = \bigcup_{C\geq 1} V_{C,0,p}
 = \{ \phi \in \mathrm{Exp}_+ \, ;\, |\phi (z)| < \infty ,\,\, \mathrm{when} \,\, 0 \leq \mathrm{Im} (z) \leq 2p \}.
\end{eqnarray*}
Recall that $i(V) \subset \overline{i(W_{3,0})}$.
Define a subspace $\mathcal{F}$ of the product space $\prod^\infty_{j=1} \mathrm{Exp}_+$ as follows.
$(Z_1, Z_2 , \cdots ) \in \mathcal{F}$ if and only if 
\\[0.2cm]
\textbf{(F1)} there exists a signed measure
 $\hat{h}(\theta , \omega )$ on $S^1$ parameterized by $\omega \in \mathbf{R}$ such that
\begin{eqnarray}
\int^{2\pi}_{0} \! \hat{h}(\theta , \omega ) d\theta =0,\,\, \quad
Z_j(\omega ) = \int^{2\pi}_{0} \! e^{\sqrt{-1}j \theta } \hat{h}(\theta , \omega ) d\theta.
\label{7-12}
\end{eqnarray}
\textbf{(F2)} Define $Z_{-1}, Z_{-2}, \cdots $ by Eq.(\ref{7-12}).
Then, there exist positive constants $C_1, C_2$ and $\gamma $ such that 
\begin{equation}
\sup_{\omega \in \mathbf{R}} |Z_j(\omega )| \leq C_1,\,\, \quad
\sup_{0\leq \mathrm{Im}(\omega ) \leq p}|Z_j(\omega )| \leq C_2 e^{|j| \gamma }
\label{F3}
\end{equation}
for all $j= \pm 1, \pm 2, \cdots $
In particular, $Z_j \in V = V_p$.
\\[0.2cm]
Hence, $\mathcal{F}$ is a vector space of Fourier coefficients $\{Z_j\}^\infty_{j=1}$ of 
signed measures included in $\prod^\infty_{j=1}V$ whose growth rate in $j$ is not so fast.
From a physical viewpoint, we are interested in an initial condition of the form 
$\hat{h}(\theta , \omega ) = h(\theta )$ (see Eq.(\ref{conti})), which satisfies (F2).
\\[-0.2cm]

The existence of solutions of Eq.(\ref{7-11}) will be prove in Sec.7.3 after $\hat{\chi}$ is specified.
In this section, we show that if solutions exist, they are included in $\mathcal{F}$.
\\[0.2cm]
\textbf{Proposition 7.2.} \, For a given function $\hat{\chi} : [0, \infty) \to [0,1]$ and an initial condition in $\mathcal{F}$,
suppose that a solution of (\ref{7-11}) exists and a function $x(t,s; \theta ,\omega )$ has an analytic continuation with respect to
$\theta $ and $\omega $ (these facts will be verified in Sec.7.3).
Then,
\\
(i) the solution is included in $\mathcal{F}$ for any $t\geq 0$.
\\
(ii) for each $t\geq 0$, $\beta = 0,1,\cdots $ and $n = 1,2,\cdots $, $|| Z_j(t, \cdot) ||^*_{\beta, n}$
is bounded uniformly in $j \in \mathbf{Z}$.
\\[0.2cm]
\textbf{Proof.} 
Note that when $\omega \in \mathbf{R}$, $\hat{Z}_j(t, \omega )$ is bounded uniformly in $j\in \mathbf{Z}$;
$|\hat{Z}_j (t, \omega )| \leq \int^{2\pi}_{0}\! |\hat{\rho}_t (\theta , \omega )| d\theta $,
which verifies the first equality of (\ref{F3}).
This also shows Part (ii); $||Z_j(t, \cdot) ||^*_{\beta, n}$ is bounded uniformly in $j\in \mathbf{Z}$.

To show that $\hat{Z}_j(t, \, \cdot \, ) \in V$ if $\hat{Z}_j(0 ,\, \cdot \,) \in V$, we use the equality
\begin{eqnarray}
& & \int^{2\pi}_{0} \! a(\theta , \omega ) \hat{\rho}_t(\theta , \omega ) d\theta = 
\int^{2\pi}_{0} \! a(x(t,0;\theta ,\omega ), \omega ) \hat{h}(\theta ,\omega ) d\theta \nonumber \\
& & + \int^{2\pi}_{0} \!\!\! \int^t_{0} \! a(x(t,0;\theta ,\omega ), \omega ) 
           \exp \Bigl[-\frac{K}{2} \! \int^s_{0} \! \left( \eta (\tau) e^{-\sqrt{-1}x(\tau,0;\theta ,\omega )}
     + \overline{\eta(\tau)}e^{\sqrt{-1}x(\tau,0;\theta ,\omega )} \right) \hat{\chi}(\tau) d\tau \Bigr] \times \nonumber \\
& & \qquad \qquad \left( \frac{K_c + \varepsilon \hat{\chi}(s)}{4\pi} \right) 
\left( \eta (s) e^{-\sqrt{-1}x(s,0;\theta ,\omega )} + \overline{\eta(s)}e^{\sqrt{-1}x(s,0;\theta ,\omega )}\right) ds d\theta,
\label{7-15}
\end{eqnarray}
for any measurable function $a(\theta , \omega )$, which is proved by substitution of (\ref{7-14}).
Note that if $\hat{\chi}(t) \equiv 1$, it is reduced to Eq.(\ref{cha3}).
From this, it turns out that $\hat{Z}_j$ is expressed as
\begin{eqnarray}
& & \hat{Z}_j(t, \omega ) =
\int^{2\pi}_{0} \! e^{\sqrt{-1}j x(t,0;\theta ,\omega )}\hat{h}(\theta ,\omega ) d\theta \nonumber \\
& & + \int^{2\pi}_{0} \!\!\! \int^t_{0} \! e^{\sqrt{-1}j x(t,0;\theta ,\omega )} 
           \exp \Bigl[-\frac{K}{2} \! \int^s_{0} \! \left( \eta (\tau) e^{-\sqrt{-1}x(\tau,0;\theta ,\omega )}
     + \overline{\eta(\tau)}e^{\sqrt{-1}x(\tau,0;\theta ,\omega )} \right) \hat{\chi}(\tau) d\tau \Bigr] \times \nonumber \\
& & \qquad \qquad \left( \frac{K_c + \varepsilon \hat{\chi}(s)}{4\pi} \right) 
\left( \eta (s) e^{-\sqrt{-1}x(s,0;\theta ,\omega )} + \overline{\eta(s)}e^{\sqrt{-1}x(s,0;\theta ,\omega )}\right) ds d\theta.
\label{7-16}
\end{eqnarray}
At first, let us show that $e^{\pm \sqrt{-1}x(t,0;\theta ,\omega )}\in V$.
This is proved in the same way as Thm.5.10 (iii).
Put $X(t) = e^{\sqrt{-1}x(t,0;\theta ,\omega )}$. Then $X$ satisfies the equation
\begin{equation}
\left\{ \begin{array}{l}
\displaystyle \frac{dX}{dt} = \sqrt{-1}\omega X + \frac{K}{2}\left( \eta (t) - \overline{\eta(t)} X^2 \right) \hat{\chi}(t),  \\[0.1cm]
X(0) = e^{\sqrt{-1}\theta }. \\
\end{array} \right.
\label{7-17}
\end{equation}
Putting $X = \xi e^{\sqrt{-1}p},\, \eta = \zeta e^{\sqrt{-1}q}$ with $\xi, \zeta, p,q \in \mathbf{R}$ yields
\begin{equation}
\frac{d\xi}{dt} = - \mathrm{Im}(\omega ) \xi + \frac{K}{2} \zeta (1- \xi^2) \cos (p-q) \hat{\chi}(t).
\label{7-18}
\end{equation}
This equation implies that if $\mathrm{Im}(\omega )\geq 0$ and $\xi = 1$, then $d\xi/dt \leq 0$.
Since $|X(0)| = 1$, we obtain $|X(t)| \leq 1$ for any $t\geq 0$ and any $\mathrm{Im}(\omega ) \geq 0$.
This proves that $X(t)$ is bounded on the real axis and the upper half plane:
$X(t) = e^{\sqrt{-1}x(t,0;\theta ,\omega )} \in V_{1, 0}\subset V$.

Next thing to do is to investigate $Y(t) = e^{-\sqrt{-1}x(t,0;\theta ,\omega )}$, which satisfies
\begin{equation}
\left\{ \begin{array}{l}
\displaystyle \frac{dY}{dt} = -\sqrt{-1}\omega Y - \frac{K}{2}\left( \eta (t)Y^2 - \overline{\eta(t)} \right) \hat{\chi}(t),  \\[0.1cm]
Y(0) = e^{-\sqrt{-1}\theta }. \\
\end{array} \right.
\label{7-19}
\end{equation}
Putting $Y = \xi e^{\sqrt{-1}p},\, \eta = \zeta e^{\sqrt{-1}q}$ with $\xi, \zeta, p,q \in \mathbf{R}$ yields
\begin{equation}
\left\{ \begin{array}{l}
\displaystyle \frac{d\xi}{dt} = \mathrm{Im}(\omega ) \xi + \frac{K}{2} \zeta (1- \xi^2) \cos (p+q) \hat{\chi}(t),  \\[0.2cm]
\displaystyle \frac{dp}{dt} = - \mathrm{Re}(\omega ) - \frac{K}{2} \zeta (\xi + \frac{1}{\xi}) \sin (p+q) \hat{\chi}(t).   \\
\end{array} \right.
\label{7-20}
\end{equation}
When $|\mathrm{Re}(\omega )|$ is sufficiently large, the averaging method is applicable to construct an approximate solution.
Eq.(\ref{7-20}) is averaged with respect to $p$ to provide the averaging equation $d\xi/dt = \mathrm{Im}(\omega ) \xi $, which is solved as
$\xi(t) = e^{\mathrm{Im}(\omega )t} \xi (0)$. Therefore, a solution of Eq.(\ref{7-20}) is given as
\begin{equation}
\xi(t) = e^{\mathrm{Im}(\omega )t} + O\Bigl( \frac{1}{|\mathrm{Re}(\omega )|} \Bigr),
\label{7-21}
\end{equation}
as $|\mathrm{Re}(\omega )| \to \infty$. See Sanders and Verhulst \cite{San} for the averaging method.
This implies that $Y$ is in $\mathrm{Exp}_+$ for each $t$ and is
bounded as $\mathrm{Re}(\omega ) \to \pm \infty$ for each $0\leq \mathrm{Im}(\omega ) \leq p$ and $t$.
Thus $e^{-\sqrt{-1}x(t,0;\theta ,\omega )} \in V$.
Therefore, the second term in the right hand side of Eq.(\ref{7-16}) is in $V$;
the second term is bounded as $\mathrm{Re}(\omega ) \to \pm \infty$ for each $j, t$ and 
$0\leq \mathrm{Im}(\omega ) \leq p$. 

\begin{figure}
\begin{center}
\includegraphics[]{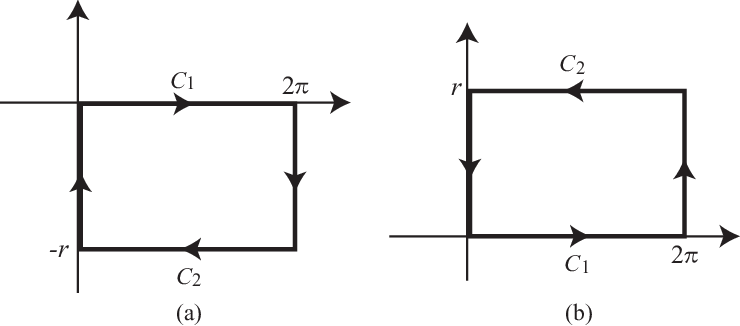}
\caption[]{The contour for obtaining Eq.(\ref{7-25}).}
\label{fig10}
\end{center}
\end{figure}

Next, we show that the first term in the right hand side of Eq.(\ref{7-16}) is in $V$.
Let
\begin{equation}
e^{\sqrt{-1}jx(t,0;\theta ,\omega )} = \sum^\infty_{n= -\infty} a_{jn}(t, \omega ) e^{\sqrt{-1}n \theta }
\label{7-22}
\end{equation}
be a Fourier expansion of $e^{\sqrt{-1}jx(t,0;\theta ,\omega )}$.
Then,
\begin{eqnarray}
\int^{2\pi}_{0} \! e^{\sqrt{-1}jx(t,0;\theta ,\omega )} \hat{h}(\theta ,\omega ) d\theta 
 &=& \sum^\infty_{n=-\infty} a_{jn}(t, \omega ) \int^{2\pi}_{0} \! e^{\sqrt{-1}n \theta } \hat{h}(\theta , \omega )d\theta \nonumber \\
 &=& \sum^\infty_{n=-\infty} a_{jn}(t, \omega )\hat{Z}_n(0, \omega ).
\label{7-23}
\end{eqnarray}
Since the series (\ref{7-22}) converges uniformly in $\theta $, the right hand side of (\ref{7-23}) exists for each $\omega $.
Since $e^{\sqrt{-1}jx(t,0;\theta ,\omega )}$ is holomorphic in $\omega $, so is $a_{jn}(t, \omega )$.
By the assumption (F2), $\hat{Z}_n(0, \omega )$ is also holomorphic.
Therefore, the right hand side of (\ref{7-23}) converges to a holomorphic function on the region $0\leq \mathrm{Im}(\omega ) \leq p$.
By the assumption (F2), there are positive constants $C$ and $\gamma $ such that
\begin{equation}
\sup_{0\leq \mathrm{Im}(\omega ) \leq p} |\hat{Z}_n(0, \omega )| \leq Ce^{|n| \gamma}.
\label{F3-2}
\end{equation}
This provides the inequality
\begin{equation}
\sup_{0\leq \mathrm{Im}(\omega ) \leq p} \Bigl| \sum^\infty_{n=-\infty} a_{jn}(t, \omega )\hat{Z}_n(0, \omega ) \Bigr| 
\leq \sup_{0\leq \mathrm{Im}(\omega ) \leq p} C \sum^\infty_{n=-\infty} e^{|n| \gamma }|a_{jn}(t, \omega )|.
\label{7-24}
\end{equation}
Let us prove that the right hand side exists.
Eq.(\ref{7-13}) shows that $x(t,0; \theta +2\pi, \omega ) = x(t,0; \theta , \omega ) + 2\pi$.
With this property, we use Cauchy's theorem to the function $e^{\sqrt{-1}jx(t,0;\theta ,\omega )}$
along the path represented in Fig.\ref{fig10}(a) to yield
\begin{eqnarray}
a_{jn}(t, \omega ) &=& \frac{1}{2\pi} \int_{C_1} \! e^{-\sqrt{-1}n\theta }e^{\sqrt{-1}j x(t,0;\theta ,\omega )} d\theta \nonumber \\
&=& -\frac{1}{2\pi} \int_{C_2} \! e^{-\sqrt{-1}n\theta }e^{\sqrt{-1}j x(t,0;\theta ,\omega )} d\theta \nonumber \\
&=& \frac{e^{-nr}}{2\pi} \int^{2\pi}_{0} \! e^{-\sqrt{-1}n\theta }e^{\sqrt{-1}j x(t,0;\theta - \sqrt{-1}r,\omega )} d\theta,
\label{7-25}
\end{eqnarray}
for $n=0,1,2,\cdots $ and $j = \pm 1, \pm 2, \cdots $,
where $r>0$ can be taken arbitrarily large because $e^{\sqrt{-1}j x(t,0;\theta ,\omega )}$ is analytic in $\theta \in \mathbf{C}$.
By the same way as above, we can show that $e^{\sqrt{-1}j x(t,0;\theta - \sqrt{-1}r,\omega )}$ is estimated as
\begin{eqnarray*}
|e^{\sqrt{-1}j x(t,0;\theta - \sqrt{-1}r,\omega )}| = e^{(-\mathrm{Im}(\omega )t + r)j} + O\Bigl( \frac{1}{|\mathrm{Re}(\omega )|} \Bigr),
\end{eqnarray*}
as $|\mathrm{Re}(\omega )| \to \infty$ for each $t$ and $\mathrm{Im}(\omega )$.
This provides
\begin{eqnarray*}
|a_{jn}(t, \omega )| \leq e^{-nr}\left( e^{(-\mathrm{Im}(\omega )t + r)j}
     + O\Bigl( \frac{1}{|\mathrm{Re}(\omega )|} \Bigr) \right)  ,
\end{eqnarray*}
for $n= 0,1,2,\cdots $.
When $n<0$, we take a path represented as Fig.\ref{fig10}(b), which yields
\begin{eqnarray*}
|a_{jn}(t, \omega )| \leq e^{nr}\left( e^{(-\mathrm{Im}(\omega )t - r)j} + O\Bigl( \frac{1}{|\mathrm{Re}(\omega )|} \Bigr) \right)
\end{eqnarray*}
as $|\mathrm{Re}(\omega )| \to \infty$ in the same way.
Therefore, we obtain
\begin{eqnarray*}
& & \sup_{0\leq \mathrm{Im}(\omega ) \leq p} \Bigl| \sum^\infty_{n=-\infty} a_{jn}(t, \omega )\hat{Z}_n(0, \omega ) \Bigr| \\
& \leq  & \sup_{0\leq \mathrm{Im}(\omega ) \leq p}  C \sum^\infty_{n=1} e^{|n|(\gamma - r)} 
      \left( e^{(-\mathrm{Im}(\omega )t + r)j} + e^{(-\mathrm{Im}(\omega )t - r)j} + O\Bigl( \frac{1}{|\mathrm{Re}(\omega )|} \Bigr) \right) .
\label{7-26}
\end{eqnarray*}
By taking $r > \gamma$, it turns out that the right hand side of Eq.(\ref{7-24}) exists and bounded as 
$\mathrm{Re}(\omega ) \to \pm \infty$ for each $j, t$ and $0\leq \mathrm{Im}(\omega ) \leq p$.
This proves that $\hat{Z}_j(t, \omega )$ is bounded as $\mathrm{Re}(\omega ) \to \pm \infty$ and
$\hat{Z}_j(t, \omega ) \in V$ for each $j$ and $t$.

To verify the second equality of (\ref{F3}), put
\begin{eqnarray*}
A(s, \theta ) &=&  \exp \Bigl[-\frac{K}{2} \! \int^s_{0} \! \left( \eta (\tau) e^{-\sqrt{-1}x(\tau,0;\theta ,\omega )}
     + \overline{\eta(\tau)}e^{\sqrt{-1}x(\tau,0;\theta ,\omega )} \right) \hat{\chi}(\tau) d\tau \Bigr] \times \nonumber \\
& & \qquad \qquad \left( \frac{K_c + \varepsilon \hat{\chi}(s)}{4\pi} \right) 
\left( \eta (s) e^{-\sqrt{-1}x(s,0;\theta ,\omega )} + \overline{\eta(s)}e^{\sqrt{-1}x(s,0;\theta ,\omega )}\right).
\end{eqnarray*}
Then, $\hat{Z}_j$ is rewritten as
\begin{eqnarray}
\hat{Z}_j(t, \omega ) = \sum^\infty_{n=-\infty}a_{jn}(t, \omega )\hat{Z}_n(0, \omega )
 + \int^{t}_{0}\!\! \int^{2\pi}_{0}\! e^{\sqrt{-1}jx(t,0;\theta , \omega )} A(s, \theta ) d\theta ds.  
\label{pr7-2}
\end{eqnarray}
From the above calculation, we obtain
\begin{eqnarray*}
\sup_{0\leq \mathrm{Im}(\omega )\leq p} |\hat{Z}_j(t, \omega )|
&\leq & \sup_{0\leq \mathrm{Im}(\omega )\leq p} \sup_{0\leq \theta \leq 2\pi} 
 C\sum^\infty_{n=0}e^{|n|(\gamma -r)} 
\left( |e^{\sqrt{-1}x(t,0; \theta -\sqrt{-1}r, \omega )}|^j + |e^{\sqrt{-1}x(t,0; \theta +\sqrt{-1}r, \omega )}|^j\right)\\
& & + \sup_{0\leq \mathrm{Im}(\omega )\leq p} \sup_{0\leq \theta \leq 2\pi} 2\pi \int^{t}_{0}\! |A(s, \theta )|ds\cdot
|e^{\sqrt{-1}x(t,0;\theta , \omega )}|^j , 
\end{eqnarray*}
which proves that $\hat{Z}_j(t, \omega )$ satisfies (\ref{F3}) for some $C$ and $\gamma $
because $e^{\sqrt{-1}x(t,0;\theta ,\omega )}$ is bounded on the region $0\leq \mathrm{Im}(\omega ) \leq p$.

Finally, let us verify (F1). Note that Eq.(\ref{7-13}) provides
\begin{eqnarray*}
\frac{\partial x}{\partial \theta }(s, 0; \theta , \omega )
 = \exp \Bigl[-\frac{K}{2} \! \int^s_{0} \! \left( \eta (\tau) e^{-\sqrt{-1}x(\tau,0;\theta ,\omega )}
     + \overline{\eta(\tau)}e^{\sqrt{-1}x(\tau,0;\theta ,\omega )} \right) \hat{\chi}(\tau) d\tau \Bigr] .
\end{eqnarray*}
This shows that the second term in the right hand side of Eq.(\ref{7-15}) is rewritten as
\begin{eqnarray*}
 \int^t_{0} \! \frac{K_c + \varepsilon \hat{\chi}(s)}{4\pi}
\int^{2\pi}_{0} \! a(x(t,0;\theta ,\omega ), \omega ) \sqrt{-1} \frac{\partial }{\partial \theta }
\left( \eta (s) e^{-\sqrt{-1}x(s,0;\theta ,\omega )}
     - \overline{\eta(s)}e^{\sqrt{-1}x(s,0;\theta ,\omega )}\right) d\theta ds.
\end{eqnarray*}
In particular, when $a(\theta , \omega ) \equiv 1$, this value vanishes because $x(s, 0; \theta , \omega )$ 
is periodic in $\theta $.
This fact and Eq.(\ref{7-15}) yield
\begin{equation}
\int^{2\pi}_{0}\!\hat{\rho}_t(\theta , \omega )d\theta  = \int^{2\pi}_{0}\!\hat{h}(\theta ,\omega )d\theta .  
\end{equation}
Therefore, if an initial condition satisfies (F1), so is $\hat{\rho}_t(\theta , \omega )$ for any $t\geq 0$.
Now the proof of Prop.7.2 (i) is completed.
\hfill $\blacksquare$
\\

\subsection{Localization of the semiflow}

By using the canonical inclusion, we rewrite Eq.(\ref{7-11}) as an evolution equation on 
$\mathbf{R} \times \prod^\infty_{k=1} \mathrm{Exp}_-'$ of the form
\begin{eqnarray}
\left\{ \begin{array}{l}
\displaystyle \frac{d}{dt}\mathcal{\varepsilon } = 0, \\[0.3cm]
\displaystyle \frac{d}{dt}  Z_1
  = T_{10}^\times Z_1
      + \frac{1}{2}\left( \varepsilon \langle Z_1  \,|\, P_0 \rangle  P_0
            - K \langle \overline{Z_1  \,|\, P_0} \rangle Z_2 \right) \hat{\chi}(t), \\[0.4cm]
\displaystyle \frac{d}{dt} Z_j 
  = T_j^\times Z_j + \frac{jK}{2} \left( \,
     \langle Z_1  \,|\, P_0 \rangle Z_{j-1}
      - \langle \overline{Z_1  \,|\, P_0 } \rangle Z_{j+1} \right) \hat{\chi}(t), 
\quad j = 2, 3,\cdots , 
\end{array} \right.
\label{7-44}
\end{eqnarray}
where $Z_j = i(\hat{Z}_j)\in \mathrm{Exp}_-'$.
The trivial equation $d\varepsilon /dt = 0$ is added to regard $\varepsilon =K - K_c$ as a dependent variable.
The operator $T_{10}$ is defined by (\ref{t1}), in which $K$ is replaced by $K_c$.
Note that $T_{10}$ has a resonance pole $\lambda _0 = 0$.

In what follows, we denote an element of the space $\mathbf{R} \times \prod^\infty_{k=1}\mathrm{Exp}_-'$ by
$z = (z_0, z_1, z_2, \cdots )$, where $z_0\in \mathbf{R}$ and 
$(z_1, z_2, \cdots ) \in \prod^\infty_{k=1}\mathrm{Exp}_-'$.
We also denote it as $z=(z_k)^\infty_{k=1}$.
A metric on this space is defined as follows:
$\mathrm{Exp}_+(\beta, n)$ is a Banach space with the norm $|| \cdot ||_{\beta, n}$ as before.
The dual space $\mathrm{Exp}_-(\beta, n)'$ of $\mathrm{Exp}_-(\beta, n)$ is a Banach space with the norm
$|| \xi ||^*_{\beta, n} = \sup_{|| \phi ||_{\beta,n} =1 } |\langle \xi  \,|\, \phi^* \rangle|$.
The projective limit $\mathrm{Exp}_-(\beta)' = \varprojlim \mathrm{Exp}_-(\beta, n)'$ is a complete metric space with the metric $d_\beta$ 
defined by Eq.(\ref{db}).
Next, for the projective limit $\mathrm{Exp}_-' = \varprojlim \mathrm{Exp}_-(\beta)'$ we introduce the metric by
\begin{equation}
d(\xi ,\zeta) = \sum^\infty_{\beta=0} \frac{1}{P(\beta )} \frac{d_\beta (\xi ,\zeta)}{1 + d_\beta (\xi ,\zeta)},
\end{equation}
where $\{ P(\beta)\}^\infty_{\beta = 0}$ is a certain sequence of positive numbers such that $\sum^\infty_{\beta=0} 1/P(\beta )$
converges.
This defines the same projective topology as that induced by the metric (\ref{5-32}).
The constants $P(\beta)$ will be determined in Sec.7.4 so that $d(\xi, 0)$ plays a similar role to a norm.
Define $d _\infty$ to be
\begin{equation}
d_\infty(z, z') = \sup_{k\geq 0} d(z_k, z_k'),
\quad z = (z_0, z_1, \cdots ) \in \mathbf{R} \times \prod^\infty_{k=1}\mathrm{Exp}_-',
\label{7-45}
\end{equation}
where $d(z_k, z_k')$ for $k\geq 1$ is the distance on $\mathrm{Exp}_-'$ defined as above, 
and $d(z_0, z_0') = |z_0 - z_0'|$ for $z_0, z_0' \in \mathbf{R}$.
Let $\mathcal{G}$ be a subspace of $\mathbf{R} \times \prod^\infty_{k=1}\mathrm{Exp}_-'$ consisting of elements $ z = (z_0, z_1, \cdots )$
such that $\sup_{k\geq 0} d(z_k, 0)$ is finite.
With this metric $d_\infty$, $\mathcal{G}$ is a complete metric vector space.
Metric vector spaces and definitions of the metrics used in this section are listed in Table 2.

\begin{table}[h]
\begin{center}
\begin{tabular}{|c|c|}
\hline 
& \\[-0.3cm]
$\displaystyle \mathrm{Exp}_+ (\beta, n)$ & 
    Banach space: $\displaystyle || \phi ||_{\beta, n} = \sup_{\mathrm{Im}(z) \geq -1/n}|\phi (z)| e^{-\beta |z|}$ \\ \hline
& \\[-0.3cm]
$\displaystyle \mathrm{Exp}_-(\beta, n)'$ & 
   Banach space:  $\displaystyle || \xi ||^*_{\beta, n} = \sup_{|| \phi ||_{\beta,n} =1 }
         |\langle \xi  \,|\, \phi^* \rangle| $ \\ \hline 
& \\[-0.3cm]
$\displaystyle \mathrm{Exp}_-(\beta)' = \varprojlim \mathrm{Exp}_- (\beta, n)'$ &  
    $\displaystyle d_\beta (\xi ,\zeta) 
               = \sum^\infty_{n=1} \frac{1}{2^n} \frac{|| \xi -\zeta ||^*_{\beta, n}}{1 + || \xi -\zeta ||^*_{\beta, n}}$ \\ \hline
$\displaystyle \mathrm{Exp}_-' = \varprojlim \mathrm{Exp}_-(\beta)'$ &  
    $\displaystyle d(\xi ,\zeta) = \sum^\infty_{\beta=0} \frac{1}{P(\beta )} \frac{d_\beta (\xi ,\zeta)}{1 + d_\beta (\xi ,\zeta)}$ \\ \hline
$\displaystyle \mathcal{G} \subset \mathbf{R} \times \prod^\infty_{k=1}\mathrm{Exp}_-'$ &  
    $\displaystyle d_\infty(z, z') =\sup_{k\geq 0}d(z_k, z_k')$ \\ \hline
& \\[-0.3cm]
$\mathcal{F},\,\, i(\mathcal{F})$ &  
    $\mathcal{F}$ is a subspace of $\prod^\infty_{k=1}\mathrm{Exp}_+$ satisfying (F1), (F2);
 \\
& 
$i(\mathcal{F}) $ is its inclusion into $\prod^\infty_{k=1}\mathrm{Exp}_-'$.\,\, $\mathbf{R} \times i(\mathcal{F}) \subset \mathcal{G}$.
\\ \hline
& \\[-0.3cm]
$\displaystyle X^{(-m)} \,( = \mathbf{R} \times \overline{i(\mathcal{F})})$ &
   $\displaystyle D_{m}(z, z') = \kappa ^m d_\infty (z, z')$  \\
& \\[-0.3cm]  \hline
$\displaystyle X \subset \prod^{\infty}_{m=0} X^{(-m)}$ &
   $\displaystyle D(u, v) = \sup_{m \geq 0} D_{m}(u^{(-m)}, v^{(-m)})$, \\
&  $u=(u^{(0)}, u^{(-1)}, \cdots ) \in \prod^{\infty}_{m=0} X^{(-m)},\,\,\, 
u^{(-m)}=(u^{(-m)}_{k})^\infty_{k=0} \in X^{(-m)}$ 
\\
\hline
\end{tabular}
\caption{Metric vector spaces used in Section 7.
Definitions of the spaces $X^{(-m)}, X$ and the constant $\kappa$ will be given in Section 7.4.}
\end{center} 
\end{table}
By the definition, it is easy to verify that
\begin{equation}
d_\infty(z, \pm z') = d_\infty(z\mp z',0) \leq d_\infty(z,0) + d_\infty(z',0).
\end{equation}
A sequence $\{z^{(m)} = (z_{k}^{(m)})^\infty_{k=0} \}^\infty_{m=1}$ in $\mathcal{G}$
converges to $x =(x_{k})^\infty_{k=0}$ if and only if $z_0^{(m)} \to x_0$ on $\mathbf{R}$ and
$|| z_k^{(m)} - x_k ||^*_{\beta, n} \to 0$ uniformly in $k\geq 1$ for every $\beta \geq 0$ and $n \geq 1$.
On the other hand, since a weakly convergent series in $\mathrm{Exp}_-'$ also converges with respect to the metric $d$,
a sequence $\{z^{(m)} = (z_{k}^{(m)})^\infty_{k=0} \}^\infty_{m=1}$ in $\mathcal{G}$
converges to $x = (x_{k})^\infty_{k=0}$ if and only if $z_0^{(m)} \to x_0$ on $\mathbf{R}$ and
$\langle z_k^{(m)} \,|\, \phi^* \rangle \to \langle x_k \,|\, \phi^* \rangle$ uniformly in $k\geq 1$ for every $\phi \in \mathrm{Exp}_+$.

Let $i(\mathcal{F})$ be the subspace of $\prod^\infty_{k=1}\mathrm{Exp}_-'$ consisting of elements of the form
\begin{eqnarray*}
(i(Z_1), i(Z_2), \cdots ), \quad \mathrm{where}\,\, (Z_1, Z_2, \cdots ) \in \mathcal{F}.
\end{eqnarray*}
Due to Eq.(\ref{F3}), $\mathbf{R} \times i(\mathcal{F}) \subset \mathcal{G}$.
Thus with the distance $d_\infty$, $\mathbf{R} \times i(\mathcal{F})$ is a metric vector subspace of $\mathcal{G}$,
and the closure $\mathbf{R} \times \overline{i(\mathcal{F})}$ is a complete metric vector space.

Now we specify the function $\hat{\chi}(t)$.
Let $\mathbf{E}_c = \mathrm{span} \{ \mu_0\}$ be the generalized center subspace of $T_{10}$.
Put
\begin{equation}
\hat{\hat{\mathbf{E}}}_c =  
\mathbf{R} \times \mathbf{E}_c \times \{0\} \times \{0\} \times \cdots  \subset \mathbf{R} \times \prod^\infty_{k=1} \mathrm{Exp}_-'.
\label{7-47}
\end{equation}
Let $P_c : \mathbf{R} \times \prod^\infty_{k=1}\mathrm{Exp}_-' \to \hat{\hat{\mathbf{E}}}_c$ 
be the projection to $\hat{\hat{\mathbf{E}}}_c$ defined by
\begin{equation}
P_c = (id_{\mathbf{R}},\, \Pi_c,\, 0,0,\cdots ),
\label{7-48}
\end{equation}
($id_{\mathbf{R}}$ is the identity on $\mathbf{R}$) and $P_s = id - P_c$ the projection to the complement of $\hat{\hat{\mathbf{E}}}_c$.
Because of Lemma 5.11 and Prop.5.13, $P_c$ and $P_s$ are continuous on $\mathbf{R} \times i(\mathcal{F})$
and $\hat{\hat{\mathbf{E}}}_c$ is included in the closure $\mathbf{R} \times \overline{i(\mathcal{F})}$.
Let $\chi(t)$ be a $C^\infty$ function such that $\chi(t) \equiv 1$ when $0 \leq t\leq 1$, $0\leq \chi (t) \leq 1$ when
$1\leq t\leq 2$, and $\chi(t) \equiv 0$ when $t\geq 2$.
Taking a small positive constant $\delta_1 $, we replace $\hat{\chi}(t)$ in (\ref{7-44}) by
\begin{equation}
\hat{\chi}(t) := \chi \left(  \frac{|| P_cz ||_{E_c}}{\delta _1}\right) \cdot \chi(|\eta (t)|) ,
\label{7-49}
\end{equation}
where
\begin{eqnarray}
z(t) = (\varepsilon ,\,  Z_1(t, \cdot) , Z_2(t, \cdot), \cdots ) \in \mathbf{R} \times \prod^\infty_{k=1}\mathrm{Exp}_-',
\label{7-49b}
\end{eqnarray}
and $|| \cdot ||_{E_c}$ is a norm on $\hat{\hat{\mathbf{E}}}_c$ defined as follows:
An element $y\in \hat{\hat{\mathbf{E}}}_c$ is denoted by $y = (y_{0}, y_{1}, 0,\cdots )$ with 
$y_{1} = \alpha \mu_0$.
Then, $|| y ||_{E_c}$ is defined to be
\begin{equation}
|| y ||_{E_c} = (|y_{0}|^2 + |\alpha |^2)^{1/2}.
\end{equation}
Thm.5.12 shows that $\Pi_c Z_1$ is given as
\begin{eqnarray}
\Pi_c  Z_1  = \alpha \mu_0, 
\quad \alpha = \frac{K_c}{2D_0}\langle \mu_0 \,|\, Z_1^* \rangle
\label{7-55}
\end{eqnarray}
With this $\alpha$, $|| P_cz ||_{E_c}$ is given by
\begin{equation}
|| P_cz ||_{E_c} = (\varepsilon ^2 + |\alpha |^2 )^{1/2},
\label{7-56}
\end{equation}
for $z = (\varepsilon ,\,  Z_1 , \cdots )\in \mathbf{R} \times i(\mathcal{F})$.
Since $\hat{\hat{\mathbf{E}}}_c$ is a finite dimensional vector space, the topology on $\hat{\hat{\mathbf{E}}}_c$ induced by 
$|| \cdot ||_{E_c}$ is equivalent to that induced by the metric $d_\infty$.
With this $\hat{\chi}(t)$, we can prove the existence of solutions of Eq.(\ref{7-44}).
\\[0.2cm]
\textbf{Proposition 7.3.}
Eq.(\ref{7-44}) with $\hat{\chi}(t)$ given by (\ref{7-49}) generates a $C^1$ semiflow $\tilde{\varphi }_t$ on $\mathbf{R} \times i(\mathcal{F})$.
That is, for a given initial condition $z\in \mathbf{R} \times i(\mathcal{F})$, Eq.(\ref{7-44}) has a unique solution
denoted by $\tilde{\varphi }_t (z)$, which is  $C^1$ in $z$, on $\mathbf{R} \times i(\mathcal{F})$ for any $t\geq 0$.
\\[0.2cm]
\textbf{Proof.}
For a given initial condition $(\hat{Z}_1(0, \, \cdot\,),  \cdots )\in \mathcal{F}$,
there exists a signed measure $\hat{h}(\theta ,\omega )$ satisfying (\ref{7-12}).
Such a $\hat{h}$ is uniquely determined because there is a one to one correspondence between a measure on $S^1$
and its Fourier coefficients (see Shohat and Tamarkin \cite{Sho}).
Thus the existence of a solution of Eq.(\ref{7-11}) follows from the existence of a solution of Eq.(\ref{7-9})
with the initial condition $\hat{\rho}_0 = \hat{h}(\theta , \omega )$.
Recall that Eq.(\ref{7-9}) is rewritten as the integro-ODE (\ref{7-13}) and (\ref{7-14}).
A proof of the existence of solutions of (\ref{7-13}), (\ref{7-14}) for any $t\geq 0$ is done by the standard iteration method
and omitted here (see \cite{Chi4} for the proof for the case $\hat{\chi} \equiv 1$).
We can also prove that $x(t,s; \theta , \omega )$ is analytic in $\theta $ and $\omega $ by the standard method.
Once a solution $\hat{\rho}_t$ of (\ref{7-9}) is obtained, a solution of (\ref{7-11}) is given through (\ref{7-10}).
Then, Prop.7.2 is applied to show that solutions of Eq.(\ref{7-11}) are included in $\mathcal{F}$.
Note that when $\hat{\chi}$ is given as (\ref{7-49}), (\ref{7-11}) becomes an autonomous system.
Therefore, solutions define a semiflow on $\mathcal{F}$.
This implies that the dual Eq.(\ref{7-44}) generates a semiflow $\tilde{\varphi }_t$ on $\mathbf{R} \times i(\mathcal{F})$.
The proof of smoothness of $\tilde{\varphi }_t$ is also proved by the iteration method.
\hfill $\blacksquare$
\\[-0.2cm]

The semiflow is also denoted componentwise as
\begin{equation}
\tilde{\varphi} _t(z) = (z_0,\, \tilde{\varphi} _t^1(z), \tilde{\varphi} _t^2(z),\cdots ),
\quad z= (z_0, z_1, \cdots ).
\label{7-50}
\end{equation}
By virtue of the variation-of-constant formula (see Eq.(\ref{vari})), $\tilde{\varphi} _t^j$ proves to be of the form
\begin{equation}
\tilde{\varphi} _t^j(z_0, z_1, \cdots ) 
= (e^{T_jt})^\times z_j + \tilde{g}_t^j (z_0, z_1, \cdots ), \quad j= 1, 2, \cdots ,
\label{7-51}
\end{equation}
where $\tilde{g}^j_t$ are nonlinear terms.
Now we introduce another localization factor.
Let $\delta _2 >0$ be a sufficiently small positive number.
By using a function $\chi (t)$ above, we multiply the function $\chi (d_\infty (P_s z , 0)/\delta _2)$
to the nonlinearity $\tilde{g}_t^j$, and define a perturbed map 
$\varphi _t = (z_0, \varphi ^1_t, \varphi ^2_t, \cdots )$ to be
\begin{equation}
\varphi_t^j(z_0, z_1, \cdots ) 
= (e^{T_jt})^\times z_j + \tilde{g}_t^j (z_0, z_1, \cdots ) \cdot \chi \Bigl(\frac{d_\infty(P_s z , 0)}{\delta _2} \Bigr),
\end{equation}
for $j = 1, 2, \cdots $. Put
\begin{eqnarray*}
g_t^j (z_0, z_1, \cdots )
 = \tilde{g}_t^j (z_0, z_1, \cdots ) \cdot \chi \Bigl(\frac{d_\infty(P_s z , 0)}{\delta _2}\Bigr).
\end{eqnarray*}
Fix a positive number $\tau > 0$, and put
\begin{eqnarray*}
L = (id_\mathbf{R}, (e^{T_{10}\tau})^\times, (e^{T_2\tau})^\times, \cdots ), 
\quad g = (0,\, g^1_\tau , g^2_\tau ,\cdots ).
\end{eqnarray*}
Then, the time $\tau$ map $\varphi _{\tau}$ of $\varphi _{t}$ is denoted as
\begin{equation}
\varphi _\tau : \mathbf{R} \times i(\mathcal{F}) \to \mathbf{R} \times i(\mathcal{F}), \quad
\varphi _{\tau}(z) = Lz + g(z).
\label{7-57}
\end{equation}
This is the desired localization of the semiflow of the original system (\ref{6-2}).
By Prop.7.3, $\tilde{g}^j_t$ is a $C^1$ mapping on $\mathbf{R} \times i(\mathcal{F})$.
Since $d_\infty (\cdot , 0)$ and $P_s$ are continuous on $\mathbf{R} \times i(\mathcal{F})$, 
$g : \mathbf{R} \times i(\mathcal{F}) \to \mathbf{R} \times i(\mathcal{F})$ is also continuous 
on $\mathbf{R} \times i(\mathcal{F})$. Hence, the map $\varphi _\tau$ 
is continuously extended to the map on the closure $\mathbf{R} \times \overline{i(\mathcal{F})}$.
Unfortunately, the distance $d_\infty (z,0)$ is not $C^1$ in $z$. However, on the region such that
$d_\infty(P_sz, 0) \leq \delta _2$ or $d_\infty(P_sz, 0) \geq 2\delta _2$, $g$ is a $C^1$ mapping because 
$\chi (d_\infty(P_sz, 0)/\delta _2)$ becomes a constant.
It is easy to see that $g(z) \sim O(z^2)$ as $z\to 0$ because the nonlinearity of Eq.(\ref{7-44}) is of $O(z^2)$. 

\begin{figure}
\begin{center}
\includegraphics[]{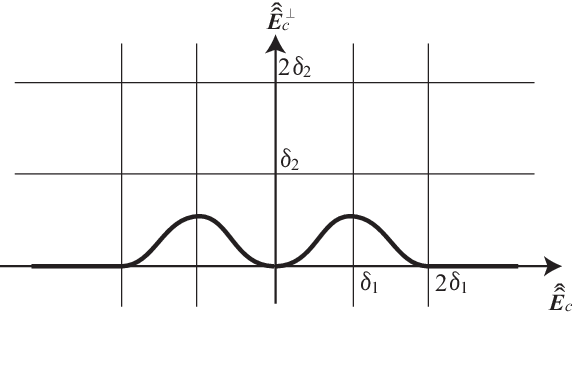}
\caption[]{A center manifold (the black curve) for the map (\ref{7-57}) coincides with that of the semiflow (\ref{7-51}).
The center manifold of the semiflow (\ref{7-51}) coincides with that of the original system (\ref{6-2})
in the region $|| P_c z||_{E_c} < \delta _1$.}
\label{fig11}
\end{center}
\end{figure}

When $|| P_cz ||_{E_c} \leq \delta_1 $, $|\eta (t)| \leq 1$ and $d_\infty(P_sz, 0) \leq \delta_2$, 
then $\chi \left( || P_cz ||_{E_c}/\delta _1 \right) = 1$, $\chi (|\eta (t)|) = 1$ and $\chi (d_\infty(P_s z , 0)/\delta _2) = 1$. 
Thus Eq.(\ref{7-57}) coincides with the time $\tau$ map of the semiflow of the original system (\ref{6-2}).
When $|| P_cz ||_{E_c} \geq 2\delta_1 $ or $d_\infty(P_sz, 0) \geq 2\delta_2$, then $\chi \left(  || P_cz ||_{E_c}/\delta _1 \right)$
$\cdot \,\chi (d_\infty(P_s z , 0) /\delta _2) = 0$.
In this case, $g = 0$ and Eq.(\ref{7-57}) is reduced to the linear map.
Therefore, by taking $\delta _1$ and $\delta _2$ sufficiently small, the Lipschitz constant of $g$
\begin{equation}
\mathrm{Lip}(g) 
:= \sup_{z, z' \in \mathbf{R} \times \overline{i(\mathcal{F})}} \frac{d_\infty(g(z), g(z'))}{d_\infty(z, z')}
\label{7-59}
\end{equation}
can be assumed to be sufficiently small.
\\[0.2cm]
\textbf{Remark.}
We introduced the factors for localization in two steps.
The one $\chi \left( || P_cz ||_{E_c}/\delta _1 \right) $ is multiplied to the nonlinearity of the equation (\ref{7-44}),
and the other $\chi (d_\infty(P_s z , 0)/\delta _2)$ is multiplied to the nonlinearity of the semiflow (\ref{7-51}).
The reason is that if we multiply both of them to the equation (\ref{7-44}), then the proof of the existence of solutions
for (\ref{7-44}) becomes too difficult;
$\chi \left( || P_cz ||_{E_c}/\delta _1 \right)$ is essentially a finite dimensional perturbation, 
although $\chi (d_\infty(P_s z , 0)/\delta _2)$ includes infinite dimensional terms $P_sz$.
On the other hand, if we multiply both of them to the nonlinearity of the semiflow of the original system (\ref{6-2}),
then a center manifold of the resultant perturbed mapping does not coincide with a center manifold of the original system (\ref{6-2})
because the perturbed mapping is not a semiflow for any differential equations in general
(i.e. the property $\varphi _{t+s} = \varphi _t \circ \varphi _s$ is violated because of the perturbation for the semiflow),
see Krisztin \cite{Kri} for details.
However, if we introduce these factors in two steps as above, a local center manifold of the original system is correctly obtained as follows:
In Sec.7.4, we will prove the existence of a center manifold for the map (\ref{7-57}).
We will show that if $\delta _1$ is sufficiently small, the center manifold is included in the ``strip" $d_\infty (P_sz, 0) < \delta _2$.
Since the map (\ref{7-57}) is the same as (\ref{7-51}), which is a semiflow of the system (\ref{7-44}), when $d_\infty (P_sz, 0) < \delta _2$,
the obtained center manifold is a center manifold of the system (\ref{7-44}).
When $|| P_c z||_{E_c} < \delta _1$, (\ref{7-44}) is reduced to the original system (\ref{6-2}).
Therefore, a local center manifold of (\ref{6-2}) is obtained as a restriction of the center manifold of (\ref{7-44})
to the region $|| P_c z||_{E_c} < \delta _1$, see Fig.\ref{fig11}.

\subsection{Proof of the center manifold theorem}

Let us prove that the mapping $\varphi _\tau$ defined in Eq.(\ref{7-57}) has a center manifold,
which gives a local center manifold for the original system (\ref{6-2}).
The strategy of the proof is the same as that in Chen, Hale and Bin \cite{Che}, in which the existence of center manifolds
is proved for mappings on Banach spaces.
At first, we need the next lemma to treat the metric $d_\infty$ as a norm.
\\[0.2cm]
\textbf{Lemma 7.4.}
For $u = (u_0, u_1, \cdots ) \in \mathbf{R} \times \overline{i(\mathcal{F})}$, suppose that there exists a positive constant
$\delta _3$ such that $|| u_j ||^*_{0,1} \leq \delta _3$ for $j=1,2,\cdots $.
If constants $\{ P(\beta)\}^\infty_{\beta =0}$ are sufficiently large,
there exist positive constants $A_c$ and $A_s = A_s(\delta _3)$ such that the inequalities
\begin{equation}
d_\infty (L^m P_c u\,,\, 0) \leq A_c d_\infty(u, 0), \quad m= 0, \pm 1, \pm2, \cdots ,
\label{7-61}
\end{equation}
and
\begin{equation}
d_\infty (L^m P_s u \, ,\, 0) \leq A_se^{-am\tau } d_\infty(u, 0), \quad m= 0, 1, 2, \cdots ,
\label{7-61b}
\end{equation}
hold, where $a > 0$ is the constant appeared in Prop.5.15.

Note that since $P_c$ is a projection to the finite dimensional vector space $\hat{\hat{\mathbf{E}}}_c$,
the linear operator $L$ restricted to $\hat{\hat{\mathbf{E}}}_c$ has the inverse $L^{-1}$ 
on $\hat{\hat{\mathbf{E}}}_c$ and $L^{-m}P_c$ is well-defined.
\\[0.2cm]
\textbf{Proof.}\,
For $u = (u_0, u_1, \cdots ) \in \mathbf{R} \times \overline{i(\mathcal{F})}$,
$L^m P_c u$ is given by
\begin{eqnarray*}
L^m P_c u = 
(u_0\, ,\, (e^{T_{10}m\tau})^\times \Pi_c u_1 \, , \,0 \,, 0\,, \cdots ).
\end{eqnarray*}
Since $u_1 \in \overline{i(V)} \subset \overline{i(W_{3,0})}$, Eq.(\ref{7-41}) is applied to yield
\begin{eqnarray*}
d_\beta ((e^{T_{10}m\tau})^\times \Pi_c u_1\,, \, 0)
 = \sum^\infty_{n=1} \frac{1}{2^n}\frac{|| (e^{T_{10}m\tau})^\times \Pi_c u_1 ||^*_{\beta ,n}}{1+
|| (e^{T_{10}m\tau})^\times \Pi_c u_1 ||^*_{\beta ,n}}
 \leq \sum^\infty_{n=1} \frac{1}{2^n} \frac{L_{\beta} || u_1 ||^*_{\beta ,n}}{1 + L_{\beta} || u_1 ||^*_{\beta ,n}},
\end{eqnarray*}
where $L_\beta = L_{3,0,\beta}$ is independent of $m$ due to $(e^{T_1t})^\times \mu_0 = \mu_0$.
We can assume without loss of generality that $L_{\beta} \geq 1$. Then,
\begin{eqnarray*}
d_\beta ((e^{T_{10}m\tau})^\times \Pi_c u_1\,, \, 0)
\leq L_{\beta} \sum^\infty_{n=1} \frac{1}{2^n} \frac{|| u_1 ||^*_{\beta ,n}}{1 + || u_1 ||^*_{\beta ,n}}
= L_{\beta} d_\beta (u_1, 0).
\end{eqnarray*}
Eq.(\ref{compara}) provides
\begin{eqnarray*}
d_0(u_1, 0) = \sum^\infty_{n=1} \frac{1}{2^n} \frac{|| u_1 ||^*_{0 ,n}}{1 + || u_1 ||^*_{0 ,n}}
\leq \sum^\infty_{n=1} \frac{1}{2^n} \frac{|| u_1 ||^*_{\beta ,n}}{1 + || u_1 ||^*_{\beta ,n}} = d_{\beta} (u_1, 0),
\end{eqnarray*}
and 
\begin{eqnarray*}
d_\beta (u_1, 0) = \sum^\infty_{n=1} \frac{1}{2^n} \frac{|| u_1 ||^*_{\beta ,n}}{1 + || u_1 ||^*_{\beta ,n}}
\leq \sum^\infty_{n=1} \frac{1}{2^n} \frac{Q(\beta)|| u_1 ||^*_{0 ,n}}{1 + || u_1 ||^*_{0 ,n}} = Q(\beta) d_{0} (u_1, 0).
\end{eqnarray*}
By using them, $d ((e^{T_{10}m\tau})^\times \Pi_c u_1\,, \, 0)$ is calculated as
\begin{eqnarray*}
d ((e^{T_{10}m\tau})^\times \Pi_c u_1\,, \, 0)
&=& \sum^\infty_{\beta =0} \frac{1}{P(\beta)} 
     \frac{d_\beta ((e^{T_{10}m\tau})^\times \Pi_c u_1\,, \, 0)}{1 + d_\beta ((e^{T_{10}m\tau})^\times \Pi_c u_1\,, \, 0)} \\
&\leq & \sum^\infty_{\beta =0} \frac{1}{P(\beta)} \frac{L_{\beta} d_\beta (u_1, 0)}{1 + L_{\beta} d_\beta (u_1 , 0)} \\
&\leq & \sum^\infty_{\beta = 0} \frac{L_{\beta} Q(\beta)}{P(\beta)} \frac{d_0(u_1, 0)}{1 + d_0(u_1, 0)},
\end{eqnarray*}
where we choose a sequence $\{ P(\beta)\}^\infty_{\beta = 0}$ so that $\sum^\infty_{\beta = 0}L_{\beta}Q(\beta)/P(\beta)$
converges. Then,
\begin{eqnarray}
d ((e^{T_{10}m\tau})^\times \Pi_c u_1\,, \, 0)
&\leq & \sum^\infty_{\beta = 0}\frac{L_{\beta} Q(\beta)}{P(\beta)}\cdot \Bigl( \sum^\infty_{\beta =0}\frac{1}{P(\beta)} \Bigr)^{-1}
\cdot \sum^\infty_{\beta = 0} \frac{1}{P(\beta)} \frac{d_0(u_1, 0)}{1 + d_0(u_1, 0)} \nonumber \\
&\leq & \sum^\infty_{\beta = 0}\frac{L_{\beta} Q(\beta)}{P(\beta)}\cdot \Bigl( \sum^\infty_{\beta =0}\frac{1}{P(\beta)}\Bigr)^{-1}
\cdot \sum^\infty_{\beta =0} \frac{1}{P(\beta)}  \frac{d_\beta(u_1, 0)}{1 + d_\beta(u_1, 0)} \nonumber \\
&=& A_c d(u_1, 0), \quad A_c := 
\sum^\infty_{\beta = 0}\frac{L_{\beta} Q(\beta)}{P(\beta)}\cdot \Bigl( \sum^\infty_{\beta =0}\frac{1}{P(\beta)}\Bigr)^{-1}.
\label{prop7-9}
\end{eqnarray}
Thus we obtain
\begin{eqnarray*}
d_\infty(L^mP_cu, 0) \leq \sup \{ u_0, A_c d(u_1, 0)\} \leq \sup \{ A_c u_0, A_c d(u_1, 0)\} = A_c d_\infty(u, 0) ,
\end{eqnarray*}
where  we suppose that $A_c \geq 1$. 
Note that we did not use the condition $|| u_j ||^*_{0,1} \leq \delta _3$ for Eq.(\ref{7-61}).

Next, $L^mP_s u$ is given by
\begin{eqnarray*}
L^{m}P_su = (0, \,\, (e^{T_{10}m\tau})^\times \Pi_s u_1\, ,\,\,(e^{T_{2}m\tau})^\times u_2 \,,\,\, \cdots  ).
\end{eqnarray*}
Eq.(\ref{7-42}) is used to yield
\begin{eqnarray*}
d_\beta ((e^{T_{10}m\tau})^\times \Pi_s u_1 , 0)
&=& \sum^\infty_{n=1} \frac{1}{2^n} \frac{|| (e^{T_{10}m\tau})^\times \Pi_s u_1 ||^*_{\beta ,n}}
{1 + || (e^{T_{10}m\tau})^\times \Pi_s u_1 ||^*_{\beta ,n}} \\
&\leq & \sum^\infty_{n=1} \frac{1}{2^n} \frac{M_{\beta} e^{-am\tau} || u_1 ||^*_{\beta ,n}}
{1 + M_{\beta} e^{-am\tau} || u_1 ||^*_{\beta ,n}} 
\leq M_{\beta} e^{-am\tau} \sum^\infty_{n=1} \frac{1}{2^n}\frac{|| u_1 ||^*_{\beta ,n}}{1 + e^{-am\tau} || u_1 ||^*_{\beta ,n}},
\end{eqnarray*}
where $M_\beta : = M_{3,0, \beta}$ is assumed to be larger than $1$.
Note that the condition $|| u_j ||^*_{0,1} \leq \delta _3$ yields $|| u_j ||^*_{\beta, n} \leq Q(\beta ) \delta _3$.
When $|| u_j ||^*_{\beta, n} \leq Q(\beta ) \delta _3$, putting $A_\beta'(\delta _3) = (1 + Q(\beta)\delta _3)$ provides
\begin{eqnarray*}
\frac{1 }{1 + e^{-am\tau} || u_1 ||^*_{\beta ,n}} \leq A'_\beta (\delta _3) \frac{1}{1 + || u_1 ||^*_{\beta ,n}},
\end{eqnarray*}
uniformly in $m= 0,1,\cdots $.
Therefore, we obtain
\begin{eqnarray*}
d_\beta ((e^{T_{10}m\tau})^\times \Pi_s u_1 , 0) \leq A'_\beta (\delta _3) M_{\beta} e^{-am\tau}d_\beta (u_1 , 0).
\end{eqnarray*}
By the same way as in Eq.(\ref{prop7-9}), we can verify that there exists a constant $A_s (\delta _3)> 0$ such that
\begin{eqnarray*}
d ((e^{T_{10}m\tau})^\times \Pi_s u_1 , 0) \leq A_s(\delta _3)e^{-am\tau}d(u_1 , 0).
\end{eqnarray*}
In this calculation, constants $P(\beta)$ are chosen sufficiently large as before.
Similarly, Eq.(\ref{7-43}) shows that
\begin{eqnarray*}
d ((e^{T_{j}m\tau})^\times u_j , 0) \leq A_s (\delta _3)e^{-am\tau}d(u_j , 0),
\end{eqnarray*}
holds for every $j = 2, 3,\cdots $.
Note that the constant $A_s$ can be taken so that it is independent of $j$ because the constant $N_{C, \alpha ,\beta}$ in 
Eq.(\ref{7-43}) is independent of $j$.
Thus $d_\infty(L^m P_s u, 0)$ satisfies Eq.(\ref{7-61b}).
\hfill $\blacksquare$
\\[-0.2cm]

If $|| P_cu ||_{E_c} \geq 2\delta _1$ or $d_\infty(P_su, 0) \geq 2\delta _2$, then $g(u) = 0$.
Thus there exists $D = D(\delta _1, \delta _2)$ such that the $j$-th component of $g$ satisfies $d(g(u)_j, 0) \leq D$
for every $j=1,2,\cdots $.
This shows that there exists $\delta _3 = \delta _3 (\delta _1, \delta _2)$ such that $|| g(u)_j ||^*_{0,1} \leq \delta _3$
for all $j$.
In what follows, we fix $\delta _3$ in Lemma 7.4 so that $g$ satisfies $|| g(u)_j ||^*_{0,1} \leq \delta _3$ for all $j$.
Then, Eq.(\ref{7-61b}) is applicable to $g(u)$.
Note that Eq.(\ref{7-61}) holds for any $\mathbf{R} \times \overline{i(\mathcal{F})}$ without the assumption.
\\[0.2cm]
\textbf{Lemma 7.5.}\, For a sequence 
$\{ u^{(-m)} = (u_{0}^{(-m)}, u_{1}^{(-m)}, \cdots  ) \}^{\infty}_{m=0} \subset \mathbf{R} \times \overline{i(\mathcal{F})}$,
suppose that
\begin{equation}
L^qP_su^{(-q)} \to 0
\label{7-61d}
\end{equation}
as $q\to \infty$, and that there exist constants $C>0$ and a sufficiently small $b >0$ such that
\begin{equation}
d_\infty(u^{(-m)} , 0) \leq C e^{b m \tau}
\label{7-61c}
\end{equation}
holds for every $m = 0,1,\cdots $. Then, $\{ u^{(-m)} \}^{\infty}_{m=0}$ satisfies
\begin{equation}
u^{(-m)} = \varphi _\tau (u^{(-m-1)}) = Lu^{(-m-1)} + g(u^{(-m-1)}),\quad m=0,1,2, \cdots ,
\label{7-62}
\end{equation}
if and only if it satisfies the equation
\begin{equation}
u^{(-m)} = L^{-m}P_cu_0 - \sum^{m}_{k=1}L^{k-m-1}P_cg(u^{(-k)}) + \sum^\infty_{k=m+1}L^{k-m-1}P_sg(u^{(-k)}),
\quad m=0,1,2, \cdots .
\label{7-63}
\end{equation}
Eq.(\ref{7-62}) means that $\{ u^{(-m)}\}^{\infty}_{m=0} = (u^{(0)}, u^{(-1)}, \cdots )$ is a negative semi-orbit of the mapping $\varphi _\tau$.
Eq.(\ref{7-63}) is called the \textit{Lyapunov-Perron equation} \cite{Che}.
\\[0.2cm]
\textbf{Proof.}
Suppose that $\{ u^{(-m)}\}$ satisfies Eq.(\ref{7-62}).
By iterating Eq.(\ref{7-62}), we obtain
\begin{eqnarray*}
u^{(0)} &=& P_cu^{(0)}  + P_su^{(0)}  \\
&=& P_cu^{(0)}  + LP_su^{(-1)}  + P_sg(u^{(-1)} ) \\
&=& P_cu^{(0)}  + L^2P_su^{(-2)}  + LP_s g(u^{(-2)} ) + P_sg(u^{(-1)} ) \\
&\vdots& \\
&=& P_cu^{(0)}  + L^qP_su^{(-q)} + \sum^q_{k=1}L^{k-1}P_sg(u^{(-k)} ), \quad q=0,1,2,\cdots .
\end{eqnarray*}
In a similar manner, we obtain
\begin{equation}
u^{(-m)}  = L^{-m}P_cu_0 - \sum^m_{k=1}L^{k-m-1}P_cg(u^{(-k)} ) + L^{q-m}P_su^{(-q)} + \sum^q_{k=m+1}L^{k-m-1}P_sg(u^{(-k)}),
\label{7-64}
\end{equation}
for $q=m,m+1, \cdots$ and $m=0,1,2,\cdots $.
By the assumption, $L^{q-m}P_su^{(-q)}\to 0$ as $q\to 0$.
Next thing to do is to show that $\sum^q_{k=m+1}L^{k-m-1}P_sg(u^{(-k)})$ converges as $q\to \infty$.
Eq.(\ref{7-61b}) is applicable to $g(u^{(-m)})$ to yield
\begin{eqnarray}
d_\infty(L^{k-m-1}P_sg(u^{(-k)} ) , 0)
&\leq & A_se^{-a (k-m-1) \tau} d_\infty(g(u^{(-k)} ) , 0) \nonumber \\
&\leq & \mathrm{Lip}(g)A_s  e^{-a(k-m-1) \tau} d_\infty(u^{(-k)}  , 0) \nonumber\\
&\leq &  \mathrm{Lip}(g)A_s C e^{-(a- b ) k \tau}e^{a (m+1) \tau},
\label{7-64b}
\end{eqnarray}
which shows that $L^{k-m-1}P_sg(u^{(-k)})$ decays exponentially as $k\to \infty$ when $a>b$.
Thus taking the limit $q\to \infty$ in Eq.(\ref{7-64}) yields Eq.(\ref{7-63}).

Conversely, suppose that $\{ u^{(-m)} \}$ satisfies Eq.(\ref{7-63}).
Because of the assumption Eq.(\ref{7-61c}), the series $\sum^\infty_{k=m+1}L^{k-m-1}P_sg(u^{(-k)})$ exists.
Replacing $m$ by $m+1$ and using $L$ for Eq.(\ref{7-63}), we obtain
\begin{equation}
Lu^{(-m-1)} = L^{-m}P_cu^{(0)}  - \sum^{m+1}_{k=1}L^{k-m-1}P_cg(u^{(-k)}) + \sum^\infty_{k=m+2}L^{k-m-1}P_sg(u^{(-k)} ).
\label{7-65}
\end{equation}
Eq.(\ref{7-65}) is put together with Eq.(\ref{7-63}) to yield Eq.(\ref{7-62}). 
\hfill $\blacksquare$
\\[-0.2cm]

Let $X^{(-m)}\,\, (m=0,1,\cdots )$ be copies of the space $\mathbf{R} \times \overline{i(\mathcal{F})}$.
Define a metric $D_{m}$ on $X^{(-m)}$ to be
\begin{equation}
D_{m}(z, z') = e^{-b m \tau } d_\infty(z, z'), \quad z, z'\in \mathbf{R} \times \overline{i(\mathcal{F})},
\label{7-67}
\end{equation}
with a small positive constant $b$.
Let $X$ be a subspace of the product $\prod^\infty_{m=0}X^{(-m)}$ consisting of elements $u=(u^{(0)} , u^{(-1)} , u^{(-2)} , \cdots )$
such that $\sup_{m} D_{m} (u^{(-m)} ,0)$ is finite. 
With the metric defined by
\begin{equation}
D(u, v) = \sup_{m\geq 0} D_{m}(u^{(-m)} , v^{(-m)} ),
\label{7-68}
\end{equation}
$X$ is a complete metric vector space (see Table 2).
It is easy to verify the inequality
\begin{equation}
D(u, \pm u') = D(u\mp u', 0) \leq D(u, 0) + D(u',0).
\label{7-69}
\end{equation}

Motivated by Eq.(\ref{7-63}), let us define the map $J : X \times \hat{\hat{\mathbf{E}}}_c \to X$ to be
\begin{equation}
\left\{ \begin{array}{l}
J(u, y) = (J^{(0)}(u, y), J^{(-1)}(u, y), J^{(-2)}(u, y) , \cdots ),  \\
\displaystyle J^{(-m)}(u,y) =  L^{-m}y - \sum^m_{k=1}L^{k-m-1}P_cg(u^{(-k)}) + \sum^\infty_{k=m+1}L^{k-m-1}P_sg(u^{(-k)}).
\end{array} \right.
\label{7-70}
\end{equation}
If the map $J(\cdot \,, y)$ has a fixed point $u = q(y) = (q^{(0)}(y), q^{(-1)}(y), \cdots )$, 
$q(y)$ is a solution of the Lyapunov-Perron equation (\ref{7-63}) with $P_cu_0 = y$.
If $q(y)$ satisfies conditions (\ref{7-61d}) and (\ref{7-61c}),
Lemma 7.5 shows that $q(y)$ is a negative semi-orbit (that is, it satisfies Eq.(\ref{7-62})) for each $y\in \hat{\hat{\mathbf{E}}}_c$.
We will see that this $q(y)$ gives a desired center manifold.
At first, let us show that $J$ is well-defined.
\\[0.2cm]
\textbf{Lemma 7.6.}
$J$ is a map from $X \times \hat{\hat{\mathbf{E}}}_c$ into $X$.
\\[0.2cm]
\textbf{Proof.}
Let us show that $D_{m}(J^{(-m)}(u,y), 0)$ is bounded uniformly in $m=0,1,\cdots $. It satisfies
\begin{eqnarray}
& & D_{m}(J^{(-m)}(u,y), 0) = e^{-b m \tau } d_\infty(J^{(-m)}(u,y), 0) \nonumber \\
&\leq & \!\!\! e^{-b m \tau }d_\infty (L^{-m}y,0) + e^{-b m \tau } \!\sum^m_{k=1}d_\infty(L^{k-m-1}P_cg(u^{(-k)}), 0)
 + e^{-b m \tau} \!\! \sum^\infty_{k=m+1} \! d_\infty (L^{k-m-1}P_s g(u^{(-k)}), 0) \nonumber \\[-0.3cm]
\label{7-71}
\end{eqnarray}
Eq.(\ref{7-61}) shows that the first term $e^{-b \tau m}d_\infty (L^{m}y,0)$ is bounded uniformly in $m=0,1,\cdots $.
Similarly, we obtain
\begin{eqnarray*}
d_\infty(L^{k-m-1}P_cg(u^{(-k)}), 0)
\leq A_c d_\infty (g(u^{(-k)}), 0) \leq \mathrm{Lip}(g) A_c d_\infty(u^{(-k)}, 0).
\end{eqnarray*}
Since $u\in X$, there is a constant $C>0$ such that $d_\infty(u^{(-k)} , 0)\leq Ce^{b k \tau}$.
Therefore,
\begin{eqnarray*}
e^{-b m \tau} \!\sum^m_{k=1}d_\infty(L^{k-m-1}P_cg(u^{(-k)}), 0)
\leq \mathrm{Lip}(g) A_c C e^{-b m\tau} \sum^m_{k=1}e^{b k\tau}
\leq \mathrm{Lip}(g) A_c C  \frac{e^{b \tau} - e^{b \tau (1-m)}}{e^{b \tau} -1}
\end{eqnarray*}
is bounded uniformly in $m = 0,1,\cdots $.
The last term in Eq.(\ref{7-71}) obviously tends to zero as $m\to \infty$.
This proves that $J(u, y) \in X$.
\hfill $\blacksquare$
\\[0.2cm]
\textbf{Proposition 7.7.}
If the constants $\delta _1$ and $\delta _2$ are sufficiently small,
$J$ is a contraction map on $X$ for each $y\in \hat{\hat{\mathbf{E}}}_c$.
\\[0.2cm]
\textbf{Proof.}
For $u, u'\in X$, we have
\begin{eqnarray}
& & D_{m}(J^{(-m)}(u,y), J^{(-m)}(u',y)) = e^{-b m \tau } d_\infty(J^{(-m)}(u,y)-J^{(-m)}(u',y), 0) \nonumber \\
&\leq & \!\!\! e^{-b m \tau } \!\sum^m_{k=1}d_\infty(L^{k-m-1}P_c(g(u^{(-k)}) - g(u'^{(-k)})), 0)
 + e^{-b m \tau } \!\! \sum^\infty_{k=m+1} \! d_\infty (L^{k-m-1}P_s (g(u^{(-k)}) - g(u'^{(-k)})), 0). \nonumber \\[-0.3cm]
\label{7-73}
\end{eqnarray}
Eqs.(\ref{7-61}) and (\ref{7-59}) are used to yield
\begin{eqnarray*}
d_\infty(L^{k-m-1}P_c(g(u^{(-k)}) - g(u'^{(-k)})), 0) \leq \mathrm{Lip} (g) A_c d_\infty (u^{(-k)}, u'^{(-k)}).
\end{eqnarray*}
Similarly, we obtain
\begin{eqnarray*}
d_\infty (L^{k-m-1}P_s (g(u^{(-k)}) - g(u'^{(-k)})), 0) \leq \mathrm{Lip}(g) A_s e^{-a(k-m-1)\tau }d_\infty (u^{(-k)}, u'^{(-k)}).
\end{eqnarray*}
Therefore, we obtain
\begin{eqnarray*}
& & D_{m}(J^{(-m)}(u,y), J^{(-m)}(u',y)) \\
&\leq & e^{-b m \tau } \sum^m_{k=1}\mathrm{Lip} (g) A_c d_\infty (u^{(-k)}, u'^{(-k)})
 + e^{-b m \tau } \sum^\infty_{k=m+1} \mathrm{Lip}(g) A_s e^{-a(k-m-1)\tau }d_\infty (u^{(-k)}, u'^{(-k)}) \\
&\leq & \mathrm{Lip} (g) A_c \sum^m_{k=1} e^{b (k-m)\tau } D_{k}(u^{(-k)}, u'^{(-k)})
 + \mathrm{Lip}(g) A_s \sum^\infty_{k=m+1} e^{-a(k-m-1)\tau }e^{b (k-m) \tau } D_{k} (u^{(-k)}, u'^{(-k)}) \\
&\leq & \mathrm{Lip} (g) \left(A_c \sum^m_{k=1} e^{b (k-m) \tau} 
 + A_s \sum^\infty_{k=m+1}  e^{-a(k-m-1)\tau }e^{b (k-m) \tau } \right) \cdot D(u, u').
\end{eqnarray*}
This yields
\begin{eqnarray*}
D(J(u, y), J(u',y)) &=& \sup_{m\geq 0}D_{m}(J^{(-m)}(u, y) , J^{(-m)}(u',y)) \\
&\leq & \mathrm{Lip} (g) \left( A_c \frac{e^{b \tau}}{e^{b \tau } - 1} 
 + A_s \frac{e^{b \tau}}{1 - e^{(b -a)\tau}} \right) \cdot D(u, u').
\end{eqnarray*}
We can take $\delta _1$ and $\delta _2$ sufficiently small so that $\mathrm{Lip} (g)$ becomes sufficiently small
and 
\begin{eqnarray}
\mathrm{Lip} (g) \left( A_c \frac{e^{b \tau}}{e^{b \tau } - 1} 
 + A_s \frac{e^{b \tau}}{1 - e^{(b -a)\tau}} \right) < 1
\label{7-73b}
\end{eqnarray}
holds. This implies that $J(\, \cdot \,,y)$ is a contraction map on $X$ for each $y\in \hat{\hat{\mathbf{E}}}_c$. \hfill $\blacksquare$
\\[0.2cm]
\textbf{Remark.}
The numbers $a$ and $b $ are the same as those in Thm.7.1.
The reason we introduced metrics $D_{m}$ and $D$ is that the center manifold is characterized by the ``slow" dynamics
whose Lyapunov exponent is smaller than $b $, see Eq.(\ref{7-3}).
The above condition for $\mathrm{Lip}(g)$ shows that if we take $b >0$ sufficiently small,
$\mathrm{Lip}(g)$ (and thus $\delta _1$ and $\delta _2$) should be small accordingly.
Since the open set $U$ in Thm.7.1, in which we can prove the existence of the local center manifold, is determined by
$\delta _1$ and $\delta _2$, $U$ also becomes small as a result.
\\[-0.2cm]

By the contraction principle, $J(\, \cdot \,,y)$ has a unique fixed point $u = q(y)$ on $X$:
\begin{equation}
\left\{ \begin{array}{l}
q(y) = (q^{(0)}(y), \, q^{(-1)}(y),\, q^{(-2)}(y), \cdots ), \quad q^{(-m)} : \hat{\hat{\mathbf{E}}}_c \to X^{(-m)},   \\
\displaystyle q^{(-m)}(y) 
= L^{-m}y - \sum^m_{k=1}L^{k-m-1}P_c g(q^{(-k)}(y)) + \sum^\infty_{k=m+1} L^{k-m-1}P_s g(q^{(-k)}(y)).  \\
\end{array} \right.
\label{7-74}
\end{equation}
In particular, $q^{(0)}$ defines a map from $\hat{\hat{\mathbf{E}}}_c$ into $X^{(0)} = \mathbf{R} \times \overline{i(\mathcal{F})}$ given by
\begin{equation}
q^{(0)}(y) = y + \sum^\infty_{k=1} L^{k-1}P_s g(q^{(-k)}(y)).
\label{7-75}
\end{equation}
Since $q(y)\in X$, there exists $C>0$ such that $D_m(q^{(-m)}(y),0) = e^{-bm\tau}d_\infty(q^{(-m)}(y),0) \leq C$,
which verifies the condition (\ref{7-61c}).
Further, Eq.(\ref{7-74}) shows that 
\begin{eqnarray*}
L^mP_sq^{(-m)}(y) = \sum^\infty_{k=m+1}L^{k-1}P_s g(q^{(-k)}(y)).
\end{eqnarray*}
Since this is a convergent series, $L^mP_sq^{(-m)}(y) \to 0$ as $m\to \infty$, which verifies (\ref{7-61d}).
Therefore, Lemma 7.5 is applicable to conclude that $\{ q^{(-m)}(y)\}^\infty_{m=0}$ is a negative semi-orbit for each $y$.
\\[0.2cm]
\textbf{Proposition 7.8.}
For any $m=0,1,\cdots $,
\\
(i) $q^{(-m)}(0) = 0$.
\\
(ii) $q^{(-m)} : \hat{\hat{\mathbf{E}}}_c \to X^{(-m)}$ is Lipschitz continuous.
\\
(iii) if $|| y ||_{E_c} \geq 2\delta _1$, then $q^{(-m)}(y) = L^{-m}y\in \hat{\hat{\mathbf{E}}}_c$.
\\
(iv) $q^{(-m)} : \hat{\hat{\mathbf{E}}}_c \to X^{(-m)}$ is a $C^1$ mapping.
In particular, $q^{(0)} : \hat{\hat{\mathbf{E}}}_c \to \mathbf{R} \times \overline{i(\mathcal{F})}$ is $C^1$.
\\[0.2cm]
\textbf{Proof.}
(i) Since $g(0)=0$, $q^{(-m)}(0) = 0$ satisfies Eq.(\ref{7-74}).
\\
(ii) For $y, y' \in \hat{\hat{\mathbf{E}}}_c $, we estimate $D_{m}(q^{(-m)}(y), q^{(-m)}(y'))$.
By the same calculation as the proof of Prop.7.7, we obtain
\begin{equation}
D(q(y), q(y')) \leq \frac{A_c}{\displaystyle 1- \mathrm{Lip} (g) \left( A_c \frac{e^{b \tau}}{e^{b \tau } - 1} 
 + A_s \frac{e^{b \tau}}{1 - e^{(b -a)\tau}} \right)}d_\infty(y, y').
\label{7-76}
\end{equation}
This means that $q : \hat{\hat{\mathbf{E}}}_c \to X$ is Lipschitz continuous.
In particular, we obtain
\begin{equation}
d_\infty(q^{(-m)}(y), q^{(-m)}(y')) 
\leq \frac{A_ce^{b m\tau}}{\displaystyle 1- \mathrm{Lip} (g) \left( A_c \frac{e^{b \tau}}{e^{b \tau } - 1} 
 + A_s \frac{e^{b \tau}}{1 - e^{(b -a)\tau}} \right)}d_\infty(y, y'),
\label{7-77}
\end{equation}
which proves the proposition.
\\
(iii) Put $y = (y_0, y_1, 0,\cdots )$, and $y_1 = \alpha \mu_0 \in \mathbf{E}_c$.
Then, the assumption implies
\begin{eqnarray*}
|| y ||_{E_c} = |y_0|^2 + |\alpha |^2 \geq 2\delta _1.
\end{eqnarray*}
On the other hand, $L^{-m}y$ is given by 
$L^{-m}y = (y_0, (e^{-T_{10}m\tau })^\times y_1 ,0,\cdots )$,
and
\begin{eqnarray*}
(e^{-T_{10}m\tau })^\times y_1 =  \alpha  (e^{-T_{10}m\tau })^\times \mu_0
 = \alpha \mu_0.
\end{eqnarray*}
Hence, $|| L^{-m}y ||_{E_c} = || y ||_{E_c} \geq 2\delta _1$.
By the construction of the nonlinearity $g$, $g(L^{-m}y) = 0$ if $|| L^{-m}y ||_{E_c} \geq 2\delta _1$.
Therefore, $q^{(-m)}(y) = L^{-m}y$ satisfies Eq.(\ref{7-74}).
\\
(iv) For $y, y^*\in \hat{\hat{\mathbf{E}}}_c$ and $\kappa \in \mathbf{R}$, put
\begin{equation}
\tilde{q}^{(-m)}(y, y^*, \kappa)= \frac{1}{\kappa}(q^{(-m)}(y + \kappa y^*) - q^{(-m)}(y)).
\label{7-78}
\end{equation}
Then, it satisfies the equation
\begin{eqnarray}
\tilde{q}^{(-m)}(y, y^*, \kappa) &=& L^{-m}y^* 
 - \frac{1}{\kappa} \sum^m_{k=1}L^{k-m-1}P_c \bigl(g(q^{(-k)}(y) + \kappa \tilde{q}^{(-k)}(y, y^*, \kappa)) - g(q^{(-k)}(y)) \bigr) \nonumber \\
& &  + \frac{1}{\kappa} \sum^\infty_{k=m+1} L^{k-m-1}P_s \bigl( g(q^{(-k)}(y) + \kappa \tilde{q}^{(-k)}(y, y^*, \kappa)) - g(q^{(-k)}(y)) \bigr),
\label{7-79}
\end{eqnarray}
for $\kappa \neq 0$.
If $\tilde{q}^{(-m)}(y, y^*, 0)$ exists, it should satisfy
\begin{eqnarray}
\tilde{q}^{(-m)}(y, y^*, 0) &=& L^{-m}y^* -
      \sum^m_{k=1}L^{k-m-1}P_c \frac{dg}{dx}(q^{(-k)}(y)) \tilde{q}^{(-k)}(y, y^*, 0) \nonumber \\
& &  +  \sum^\infty_{k=m+1} L^{k-m-1}P_s \frac{dg}{dx} (q^{(-k)}(y))\tilde{q}^{(-k)}(y, y^*, 0).
\label{7-80}
\end{eqnarray}
Motivated by these equations, we define a map $J': X\times \hat{\hat{\mathbf{E}}}_c\times \hat{\hat{\mathbf{E}}}_c
\times \mathbf{R} \to X$ to be $J' = (J'_0 , J'_{-1}, J'_{-2}, \cdots )$ and
\begin{equation}
J'_{-m}(u, y, y^*, \kappa) = \left\{ \begin{array}{l}
\displaystyle L^{-m}y^* - \frac{1}{\kappa} \sum^m_{k=1}L^{k-m-1}P_c \bigl(g(q^{(-k)}(y) + \kappa u^{(-k)}) - g(q^{(-k)}(y)) \bigr) \\
\displaystyle \quad \quad \quad 
 + \frac{1}{\kappa} \sum^\infty_{k=m+1} L^{k-m-1}P_s \bigl( g(q^{(-k)}(y) + \kappa u^{(-k)}) - g(q^{(-k)}(y)) \bigr),\quad \kappa \neq 0, \\
\displaystyle L^{-m}y^* -  \sum^m_{k=1}L^{k-m-1}P_c \frac{dg}{dx}(q^{(-k)}(y)) u^{(-k)} \\
\displaystyle \quad \quad \quad 
 +  \sum^\infty_{k=m+1} L^{k-m-1}P_s \frac{dg}{dx} (q^{(-k)}(y))u^{(-k)}, \quad \kappa = 0.
\end{array} \right.
\label{7-81}
\end{equation}
We can prove that $J'$ is a contraction map on $X$ for each $y, y^*$ and $\kappa$ by the completely same way as the proofs
of Lemma 7.6 and Proposition 7.7.
Hence, there uniquely exists $u^{(-m)} = \tilde{q}^{(-m)}(y, y^*, \kappa)$ satisfying Eq.(\ref{7-79}) and Eq.(\ref{7-80}).
Taking the limit $\kappa \to 0$ in Eq.(\ref{7-79}) yields
\begin{eqnarray}
\lim_{\kappa \to 0 }\tilde{q}^{(-m)}(y, y^*, \kappa ) &=& L^{-m}y^* 
     - \sum^m_{k=1}L^{k-m-1}P_c \frac{dg}{dx}(q^{(-k)}(y)) \lim_{\kappa \to 0 }\tilde{q}^{(-k)}(y, y^*, \kappa ) \nonumber \\
& &  +  \sum^\infty_{k=m+1} L^{k-m-1}P_s \frac{dg}{dx} (q^{(-k)}(y))\lim_{\kappa \to 0 }\tilde{q}^{(-k)}(y, y^*, \kappa).
\end{eqnarray}
This implies that $\lim_{\kappa \to 0 }\tilde{q}^{(-m)}(y, y^*, \kappa )$ is a solution of Eq.(\ref{7-80}).
By the uniqueness of a solution, we obtain
\begin{equation}
\lim_{\kappa \to 0 }\tilde{q}^{(-m)}(y, y^*, \kappa )
 = \lim_{\kappa \to 0 }\frac{1}{\kappa} \left( q^{(-m)}(y + \kappa y^*) - q^{(-m)}(y)\right) = \tilde{q}^{(-m)}(y, y^*, 0 ).
\end{equation}
From Eq.(\ref{7-80}), it turns out that $\tilde{q}^{(-m)}(y, y^*, 0 )$ is linear in $y^*$.
Thus we denote it as
\begin{equation}
\tilde{q}^{(-m)}(y, y^*, 0 )= dq^{(-m)}(y)y^*.
\end{equation}
Then, $dq^{(-m)}(y) : \hat{\hat{\mathbf{E}}}_c \to X^{(-m)}$ defines a linear operator for each $y\in \hat{\hat{\mathbf{E}}}_c$.
The remaining task is to show that $dq^{(-m)} : \hat{\hat{\mathbf{E}}}_c\times \hat{\hat{\mathbf{E}}}_c \to X^{(-m)}$ is continuous.
This is done in the same way as the proof of part (ii) of the proposition.
For $y', y'^*\in \hat{\hat{\mathbf{E}}}_c$, we estimate $d_\infty(\tilde{q}^{(-m)}(y, y^*, 0 ), \tilde{q}^{(-m)}(y', y'^*, 0 ))$.
Then, we can show that $\tilde{q}^{(-m)}$ is Lipschitz continuous in $y$ and $y^*$.
The details are left to the reader.
This means that $dq^{(-m)}(y)$ gives the derivative of $q^{(-m)}$ at $y\in \hat{\hat{\mathbf{E}}}_c$.
\hfill $\blacksquare$
\\[-0.2cm]

Now we define the center manifold $W^c$ of the map $\varphi _\tau$ by
\begin{equation}
W^c = \{q^{(0)}(y) = y + \hat{q}(y) \, | \, y\in \hat{\hat{\mathbf{E}}}_c\},
\end{equation}
where
\begin{equation}
\hat{q}(y) = \sum^\infty_{k=1} L^{k-1}P_s g(q^{(-k)}(y)) \in \hat{\hat{\mathbf{E}}}_c^\bot.
\end{equation}
\\
\textbf{Proposition 7.9.}
(i) $W^c$ is a dim-$\hat{\hat{\mathbf{E}}}_c$ dimensional $C^1$ manifold, which is tangent to the space $\hat{\hat{\mathbf{E}}}_c$.
In particular, $q^{(0)}(y)$ is expanded as $q^{(0)}(y) = y + O(y^2)$ as $y\to 0$.
\\
(ii) $W^c$ is $\varphi _\tau$ invariant; that is, $\varphi _\tau (W^c) \subset W^c$.
\\
(iii) For any $\xi_0\in W^c$, there exists a negative semi-orbit $\{ u^{(-m)}\}^{\infty}_{m=0} \subset W^c$
satisfying $u_0 = \xi_0$ and
\begin{eqnarray*}
d_\infty(u^{(-m)}, 0) \leq Ce^{b m\tau},
\end{eqnarray*}
where $b >0$ as above and $C$ is a positive constant.
\\
(iv) if $\delta _1 > 0$ is sufficiently small, the center manifold $W^c$
is included in the strip region $\{ z\in \mathbf{R} \times \overline{i(\mathcal{F})} \, | \, 
d_\infty(P_sz , 0) \leq \delta _2\}$ (see Fig.\ref{fig11}).
\\[0.2cm]
\textbf{Proof.}
(i) Since $q^{(-k)}(0) = 0$ and $q^{(-k)}(y)$ is $C^1$, $q^{(-k)}(y)$ is expanded as $q^{(-k)}(y) \sim O(y)$.
This shows that $\hat{q}(y) \sim O(y^2)$ because $g(z) \sim O(z^2)$ as $z\to 0$.
\\
(ii) Recall that $\{ q^{(-m)}(y)\}^{\infty}_{m=0}$ is a negative semi-orbit satisfying Eqs.(\ref{7-61d}) and (\ref{7-61c}).
Define $q^{(1)}(y):= \varphi _\tau (q^{(0)}(y))$.
Obviously $\{ q^{(-m+1)}(y)\}^{\infty}_{m=0}$ is also a negative semi-orbit satisfying (\ref{7-61d}) and (\ref{7-61c}) with some $C>0$.
Then, Lemma 7.5 implies that $\{ q^{(-m+1)}(y)\}^{\infty}_{m=0}$ is a solution of the Lyapunov-Perron equation (\ref{7-63}).
By the uniqueness of a solution, there exists $y'\in  \hat{\hat{\mathbf{E}}}_c$ such that $q^{(-m+1)}(y) = q^{(-m)}(y')$ for $m=0,1,\cdots $.
In particular, we obtain $\varphi _\tau (q^{(0)}(y)) = \varphi _\tau (q^{(-1)}(y')) = q^{(0)}(y')$, which proves 
$\varphi _\tau (q^{(0)}(y))\in W^c$.
\\
(iii) This is obvious from the definition: if $\xi_0 = q^{(0)}(y)$, $\{ q^{(-m)}(y)\}^{\infty}_{m=0}$ is a negative semi-orbit
included in $W^c$.
\\
(iv) Prop.7.8 (iii) implies that $P_sq^{(0)}(y) = 0$ if $|| y ||_{E_c} \geq 2\delta _1$.
Thus $\sup_{y\in \hat{\hat{\mathbf{E}}}_c } d_\infty(P_sq^{(0)}(y), 0)$ becomes sufficiently small
if $\delta _1$ is sufficiently small.
\hfill $\blacksquare$
\\[-0.2cm]

If restricted to a small neighborhood of the origin, $W^c$ gives the desired local center manifold for Eq.(\ref{6-2}).
\\[0.2cm]
\textbf{Proof of Theorem 7.1.}
If $\delta _1 > 0$ is sufficiently small, $W^c$ is included in the region $\{ z\in \mathbf{R} \times \overline{i(\mathcal{F})} \, | \, 
d_\infty(P_sz , 0) \leq \delta _2\}$, on which $\chi (d_\infty(P_s z, 0)/\delta _2) = 1$.
Thus $\varphi _\tau$-invariant manifold $W^c$ is also invariant under the map $\tilde{\varphi }_\tau$ given by Eq.(\ref{7-51}),
which is a time $\tau$ map of the semiflow of the system (\ref{7-44}).
Take $u_0\in W^c$.
By Prop.7.9 (iii), there is a negative semi-orbit $\{ u^{(-m)}\}^{\infty}_{m=0} \subset W^c$ of $\tilde{\varphi }_\tau$
satisfying Eqs.(\ref{7-61d}) and (\ref{7-61c}).
Since $\tilde{\varphi }_\tau$ is a semiflow, we have
\begin{equation}
\tilde{\varphi }_\tau \circ \tilde{\varphi }_t (u^{(-m)}) = \tilde{\varphi }_t (\tilde{\varphi }_\tau (u^{(-m)}))
 = \tilde{\varphi }_t (u^{(-m+1)}),
\end{equation}
for each $t>0$. This means that $\{ \tilde{\varphi }_t (u^{(-m)})\}^{\infty}_{m=0}$ is a negative semi-orbit of $\tilde{\varphi }_\tau$.
Since $\tilde{\varphi }_t$ is $C^1$ with respect to the metric $d_\infty$, there is a positive number $\tilde{C}$ such that
\begin{equation}
d_\infty(\tilde{\varphi }_t (u^{(-m)}), 0) \leq \tilde{C} d_\infty(u^{(-m)} , 0) \leq \tilde{C}C e^{bm\tau}.
\end{equation}
Further, $L^mP_s \tilde{\varphi }_t(u^{(-m)})$ is estimated as
\begin{eqnarray*}
L^mP_s \tilde{\varphi }_t(u^{(-m)})=  L^{m+1} P_su^{(-m)} + L^mP_s\tilde{g}_t(u^{(-m)}).
\end{eqnarray*}
Since $u^{(-m)}$ satisfies (\ref{7-61d}), $L^{m+1} P_su^{(-m)}$ tends to zero as $m\to \infty$.
By the same calculation as Eq.(\ref{7-64b}), we see that the second term $L^mP_s\tilde{g}_t(u^{(-m)})$ also
tends to zero as $m\to \infty$.
This shows that $\{ \tilde{\varphi }_t (u^{(-m)})\}^{\infty}_{m=0}$ satisfies the conditions (\ref{7-61d}) and (\ref{7-61c}).
Therefore, it is a solution of the Lyapunov-Perron equation (\ref{7-63}).
By the uniqueness of a solution, there is $y\in \hat{\hat{\mathbf{E}}}_c$ such that $\tilde{\varphi }_t (u^{(0)}) = q^{(0)}(y)\in W^c$,
which proves that $W^c$ is $\tilde{\varphi }_t $-invariant.

In Eq.(\ref{7-44}), since $\varepsilon $ is a constant which is independent of $t$, $W^c(\varepsilon ) :=
W^c \cap (\{\varepsilon \} \times \overline{i(\mathcal{F})})$ is also $\tilde{\varphi }_t$-invariant for each $\varepsilon $.

On the region $\hat{U} = \{ z \, | \, || P_cz ||_{E_c} \leq \delta _1,\, |\eta (t)| \leq 1 \}$,
$\chi (|| P_cz ||_{E_c}/\delta _1)\cdot \chi (|\eta (t)|) = 1$ and Eq.(\ref{7-44}) is reduced to the original system (\ref{6-2}).
Thus $W^c(\varepsilon ) \cap \hat{U}$ is invariant under the semiflow generated by (\ref{6-2}),
which gives a local center manifold stated in Thm.7.1 with $U = \hat{U} \cap (\{\varepsilon \} \times \overline{i(\mathcal{F})})$
and $W^c_{loc} = W^c(\varepsilon ) \cap U$.
Parts (I) and (II) in Thm.7.1 immediately follows from Prop.7.9.
It remains to show the part (III) of Thm.7.1.
This is proved in the same way as Chen, Hale and Tan \cite{Che}.
In \cite{Che}, the existence of invariant foliations is proved for dynamical systems on Banach spaces.
Though our phase space $\overline{i(\mathcal{F})}$ is not a Banach space, the distance from the origin
$d_\infty(z, 0)$ plays the same role as a norm.
Thus with the aid of the estimates (\ref{7-61}) and (\ref{7-61b}), we can prove the existence of invariant foliations 
by the same way as \cite{Che}.
The details are left to the reader. \hfill $\blacksquare$


\subsection{Reduction to the center manifold}

Let us derive the dynamics on the center manifold and prove the Kuramoto's conjecture.
Recall that for the continuous limit (\ref{conti}) of the Kuramoto model,
Putting $Z_j(t, \omega ) = \int^{2\pi}_{0}\! e^{j\sqrt{-1}\theta }\rho_t(\theta , \omega )d\theta $ yields
the system of equations (\ref{4-0}) and (\ref{4-0b}).
Since solutions are included in $V_{1,0} \subset \mathrm{Exp}_+$ (Thm.5.10 (iii)), the canonical inclusion is applied to
rewrite Eq.(\ref{4-0}) and (\ref{4-0b}) as equations of the form (\ref{6-2}) defined on $\mathrm{Exp}_-'$.
The order parameter $\eta$ is defined as $\eta (t) = (Z_1, P_0) = \langle Z_1 \,|\, P_0 \rangle$.
For this system, we have proved that when $0<K < K_c$, the trivial solution (de-synchronous state) is asymptotically stable
because of the existence of resonance poles on the left half plane.
In particular, $\eta (t) \to 0$ as $t\to \infty$.
When $K>K_c$, we have proved that the trivial solution is unstable because of the existence of eigenvalues on the 
right half plane.
Thus a bifurcation from the trivial solution may occur at $K = K_c$.
In Sec.7.1 to Sec.7.4, we have proved that there exists a smooth local center manifold near the origin in 
$\prod^\infty_{k=1}\mathrm{Exp}_-'$ if $K$ is sufficiently close to $K_c$.
Our purpose is to obtain a differential equation describing the dynamics on the center manifold
to reveal a bifurcation structure of the Kuramoto model.

Since $g(\omega )$ is the Gaussian, there exists only one resonance pole $\lambda _0 = 0$ on the imaginary axis when $K=K_c$.
Thus the center subspace $\mathbf{E}_c$ is of one dimensional.
Let $\mu_0$  be the generalized eigenfunction associated with $\lambda _0 = 0$.
By the definition, $\mu_0$ is given by
\begin{equation}
\langle \mu_0 \,|\, \phi^* \rangle
 = \lim_{x\to +0}\int_{\mathbf{R}} \! \frac{1}{x - \sqrt{-1}\omega }\phi (\omega )g(\omega ) d\omega.
\label{8-1}
\end{equation}
This is also written as
\begin{equation}
\mu_0 = \lim_{x\to +0}\, i(\frac{1}{x - \sqrt{-1}\omega }),
\label{8-2}
\end{equation}
where the limit is taken with respect to the weak dual topology on $\mathrm{Exp}_-'$.
The main theorem in this section, which confirms the Kuramoto's conjecture, is stated as follows:
\\[0.2cm]
\textbf{Theorem 7.10.}\, For the continuous model (\ref{conti}) of the Kuramoto model,
there exist positive constants $\varepsilon _0$ and $\delta $ such that
if $K_c < K < K_c + \varepsilon _0$ and if the initial condition $h(\theta )$ satisfies
\begin{equation}
\left| \int^{2\pi}_{0} \! e^{\sqrt{-1}j \theta }h(\theta ) d\theta  \right| < \delta 
\label{8-3}
\end{equation}
for $j = 1,2,\cdots $, then the order parameter $\eta (t)$ tends to the constant expressed as
\begin{equation}
r(t) = |\eta (t) |= \sqrt{\frac{-16}{\pi K_c^4 g''(0)}}\sqrt{K - K_c} + O(K-K_c), 
\label{8-4}
\end{equation}
as $t\to \infty$. In particular, the bifurcation diagram of the order parameter is given as Fig.\ref{fig2} (a).
\\[0.2cm]
\noindent \textbf{Proof.} 
Suppose that an initial condition $h(\theta )$ satisfies Eq.(\ref{8-3}).
Then, we have
\begin{eqnarray*}
|| \, Z_j(0, \cdot)  \, ||^*_{\beta,n}
 &=& \sup_{|| \phi ||_{\beta ,n} = 1} \Bigl| \int_{\mathbf{R}}\! Z_j(0, \omega ) \phi (\omega )g(\omega )d\omega  \Bigr| \\
&=& \sup_{|| \phi ||_{\beta ,n} = 1} \Bigl| \int_{\mathbf{R}}\! \int^{2\pi}_{0}\!
         e^{\sqrt{-1}j\theta } h(\theta ) \phi (\omega )g(\omega ) d\theta d\omega \Bigr| \\
&=& \Bigl| \int^{2\pi}_{0}\! e^{\sqrt{-1}j\theta } h(\theta ) d\theta \Bigr| \cdot
      \sup_{|| \phi ||_{\beta ,n} = 1} \Bigl| \int_{\mathbf{R}}\! \phi (\omega )g(\omega ) d\omega \Bigr| \\
&\leq & \delta \cdot || P_0 ||^*_{\beta ,n},
\end{eqnarray*}
for every $j, \beta$ and $n$.
Thus we can take $\delta $ sufficiently small so that the initial condition
$(Z_1(0, \cdot),\, Z_2(0, \cdot) ,\cdots )$ for Eq.(\ref{6-2}) is included in the neighborhood
$U$ (with respect to the metric $d_\infty$) of the origin given in Thm.7.1.
Then, the center manifold theorem is applicable.
Let us derive the dynamics on the center manifold.

Since we are interested in a bifurcation at $K = K_c$,
put $\varepsilon  = K - K_c$ and divide the operator $T_1$ as
\begin{equation}
T_1\phi (\omega ) = T_{10} \phi (\omega ) + \frac{\varepsilon }{2} \langle   \phi  \,|\, P_0 \rangle P_0(\omega ),
\label{8-5}
\end{equation}
where
\begin{equation}
T_{10} \phi (\omega ) 
= \sqrt{-1}\omega \phi (\omega ) + \frac{K_c}{2} \langle \phi  \,|\, P_0 \rangle P_0(\omega ).
\label{8-6}
\end{equation}
Then, the operator $T_{10}$ has a resonance pole at the origin and all other resonance poles are on the left half plane.
Eq.(\ref{6-2}) is rewritten as
\begin{equation}
\frac{d}{dt} Z_1
  = T_{10}^\times  Z_1  
       + \frac{\varepsilon }{2} \langle  Z_1 \,|\, P_0 \rangle P_0 
           - \frac{K}{2} \langle \overline{Z_1 \,|\,P_0}  \rangle Z_2.
\label{8-7}
\end{equation}
To obtain the dynamics on the center manifold, by using the spectral decomposition, we put
\begin{equation}
Z_1 = \frac{K_c}{2}\alpha (t) \mu_0 +  Y_1,
\label{8-8}
\end{equation}
where $\mu_0 $ is defined by Eq.(\ref{8-2}), $Y_1 $ is included in the complement $\mathbf{E}_c^\bot$
of $\mathbf{E}_c$, and where
\begin{equation}
\alpha (t) = \frac{1}{D_0} \langle \mu_0 \, | \, Z_1^* \rangle.
\label{8-9}
\end{equation}
We will derive the dynamics of $\alpha $.
Since $\langle \mu_0 \,|\, P_0 \rangle = 2/K_c$ by the definition of resonance poles,  we obtain
\begin{equation}
\langle Z_1 \,|\, P_0 \rangle = \alpha (t) + \langle Y_1  \,|\, P_0 \rangle,
\label{8-10}
\end{equation}
and $P_0 $ is decomposed as
\begin{equation}
P_0 = \frac{1}{D_0} \mu_0 + Y_0,
\label{8-11}
\end{equation}
where $Y_0 \in \mathbf{E}_c^\bot$.
By Thm.7.1 (I), on the local center manifold, we can suppose that
\begin{equation}
\langle Y_1 \,|\,\phi^* \rangle,\, \langle Z_j \,|\, \phi^* \rangle
\sim O(\alpha ^2, \alpha \varepsilon , \varepsilon ^2),
\label{8-12}
\end{equation}
for $j=2,3,\cdots $ and for every $\phi \in \mathrm{Exp}_+$.
Let us calculate the expression of the center manifold.
Substituting Eqs.(\ref{8-10}),(\ref{8-8}) into Eq.(\ref{6-2}) for $j=2$ yields
\begin{equation}
\frac{d}{dt}Z_2 = T_2^\times  Z_2
 + K \left( (\alpha + \langle Y_1 \,|\, P_0 \rangle) 
         \left( \frac{K_c}{2} \alpha \mu_0 + Y_1 \right)
             - (\overline{\alpha }+ \langle \overline{Y_1 \,|\, P_0} \rangle)  Z_3\right).
\label{8-13}
\end{equation}
We suppose that $d\alpha /dt \sim O(\alpha ^2, \alpha \varepsilon , \varepsilon ^2)$, which will be justified later.
Then, the above equation yields
\begin{equation}
T_2^\times  Z_2 = -\frac{KK_c}{2} \alpha ^2 \mu_0 
  + O(\alpha ^3, \alpha ^2\varepsilon , \alpha \varepsilon ^2, \varepsilon ^3).
\label{8-14}
\end{equation}
\\[0.2cm]
\textbf{Lemma 7.11.}\, Define the operator $(T_2^\times)^{-1} : i(\mathrm{Exp}_+) \to \mathrm{Exp}_-'$ to be
\begin{equation}
\langle (T_2^\times)^{-1}\psi \,|\, \phi^* \rangle
  = -\frac{1}{2}\lim_{x\to +0} \int_{\mathbf{R}} \! \frac{1}{x-\sqrt{-1}\omega} \phi(\omega ) \psi(\omega ) g(\omega )d\omega . 
\label{8-15}
\end{equation}
Then,
\begin{equation}
(T_2^\times) (T_2^\times)^{-1} \psi = (T_2^\times)^{-1}(T_2^\times) \psi = \psi
\label{8-16}
\end{equation}
for any $\psi \in i(\mathrm{Exp}_+)$, and it is continuous on $i(V)$.
\\[0.2cm]
\textbf{Proof.}\, The straightforward calculation shows that 
\begin{eqnarray*}
\langle (T_2^\times) (T_2^\times)^{-1} \psi \,|\,\phi^* \rangle &=& 
     \langle(T_2^\times)^{-1} \psi \,|\, T_2^* \phi^* \rangle \\
&=& -\frac{1}{2}\lim_{x\to +0} 
        \int_{\mathbf{R}} \! \frac{2\sqrt{-1}\omega }{x-\sqrt{-1}\omega} \phi(\omega ) \psi(\omega ) g(\omega )d\omega \\
&=& \int_{\mathbf{R}} \! \phi (\omega ) \psi (\omega ) g(\omega ) d\omega  
 - \lim_{x\to +0} \int_{\mathbf{R}} \! \frac{x}{x-\sqrt{-1}\omega} \phi(\omega ) \psi(\omega ) g(\omega )d\omega.
\end{eqnarray*}
Since the limit
\begin{eqnarray*}
\lim_{x\to +0} \int_{\mathbf{R}} \! \frac{1}{x-\sqrt{-1}\omega} \phi(\omega ) \psi(\omega ) g(\omega )d\omega
 = \langle  \mu_0  \,|\, \phi^* \cdot \psi^* \rangle
\end{eqnarray*}
exists, the second term in the right hand side above is zero.
Thus we obtain
\begin{eqnarray*}
\langle (T_2^\times) (T_2^\times)^{-1} \psi \,|\,\phi^* \rangle =  \langle \psi \,|\, \phi^*  \rangle.
\end{eqnarray*}
In the same way,
\begin{eqnarray*}
\langle (T_2^\times)^{-1}(T_2^\times)  \psi \,|\, \phi^* \rangle 
  &=& \langle (T_2^\times)^{-1} \cdot 2\sqrt{-1}\omega \psi \,|\, \phi^* \rangle \\
&=&  -\frac{1}{2}\lim_{x\to +0} 
         \int_{\mathbf{R}} \! \frac{2\sqrt{-1}\omega }{x-\sqrt{-1}\omega} \phi(\omega ) \psi(\omega ) g(\omega )d\omega \\
&=& \langle \psi \,|\, \phi^*  \rangle.
\end{eqnarray*}
Note that the right hand side of (\ref{8-15}) is also written as $-\langle \mu_0 \,|\, \phi^* \cdot \psi^* \rangle/2$.
Thus the continuity of $(T_2^\times)^{-1}$ on $i(V)$ follows from Lemma 5.11 (iii).
\hfill $\blacksquare$
\\[-0.2cm]

Since $(T_2^\times)^{-1}$ is continuous on $i(V)$, its domain is continuously extended to the closure $\overline{i(V)}$. 
Since $\mu_0 \in \overline{i(V)}$ (see Prop.5.13) and it is given as
Eq.(\ref{8-2}), $(T_2^\times)^{-1}\mu_0$ is calculated as
\begin{equation}
\langle (T_2^\times)^{-1} \mu_0 \,|\, \phi^*  \rangle 
    = \lim_{x\to +0} \Bigl\langle (T_2^\times)^{-1} \frac{1}{x- \!\sqrt{-1}\omega } \,\Bigl|\, \phi^* \Bigr\rangle
 = -\frac{1}{2} \lim_{x\to +0} \int_{\mathbf{R}} \! \frac{1}{(x- \!\sqrt{-1}\omega )^2} \phi (\omega ) g(\omega )d\omega . 
\label{8-17}
\end{equation}
Then, Eq.(\ref{8-14}) provides
\begin{equation}
\langle Z_2 \,|\, \phi^* \rangle
   = \frac{KK_c}{4}\alpha ^2 \lim_{x\to +0} \int_{\mathbf{R}} \! \frac{1}{(x- \sqrt{-1}\omega )^2}
\phi (\omega ) g(\omega ) d\omega + O(\alpha ^3, \alpha ^2\varepsilon , \alpha \varepsilon ^2, \varepsilon ^3), 
\label{8-18}
\end{equation}
which gives the expression of the center manifold to the $Z_2$ direction.
The projection of it to the center subspace is given as
\begin{eqnarray*}
\Pi_0 Z_2
   &=& \frac{K_c}{2D_0} \langle  \mu_0 \,|\, Z_2^* \rangle  \cdot \mu_0 
 = \frac{K_c}{2 D_0} 
       \lim_{x\to +0} \Bigl\langle \frac{1}{x - \sqrt{-1}\omega } \,\Bigl|\, Z_2^* \Bigl\rangle \cdot \mu_0,\\
&=& \frac{K_c}{2 D_0} 
       \lim_{x\to +0} \Bigl\langle Z_2 \,\Bigl|\, \left( \frac{1}{x-\sqrt{-1}\omega }\right)^* \Bigl\rangle \cdot \mu_0
\end{eqnarray*}
where
\begin{eqnarray*}
\lim_{x\to +0}\Bigl\langle Z_2 \,\Bigl|\, \left( \frac{1}{x-\sqrt{-1}\omega }\right)^* \Bigl\rangle
 &=& \frac{KK_c}{4} \alpha ^2 \lim_{x\to +0} \int_{\mathbf{R}} \! \frac{1}{(x- \sqrt{-1}\omega )^3} g(\omega )d\omega 
          + O(\alpha ^3, \alpha ^2\varepsilon , \alpha \varepsilon ^2, \varepsilon ^3) \\
&=& -\frac{KK_c}{8} \alpha ^2 \lim_{x\to +0} \int_{\mathbf{R}} \! \frac{1}{x- \sqrt{-1}\omega } g''(\omega )d\omega 
          + O(\alpha ^3, \alpha ^2\varepsilon , \alpha \varepsilon ^2, \varepsilon ^3) \\
&=& -\frac{KK_c}{8} \alpha ^2 \cdot \pi g''(0) + O(\alpha ^3, \alpha ^2\varepsilon , \alpha \varepsilon ^2, \varepsilon ^3).
\end{eqnarray*}
Thus we obtain
\begin{equation}
\Pi_0 Z_2 = -\frac{KK_c^2}{16D_0} \alpha ^2 \cdot \pi g''(0) \cdot \mu_0
 + O(\alpha ^3, \alpha ^2\varepsilon , \alpha \varepsilon ^2, \varepsilon ^3).
\label{8-19}
\end{equation}
Finally, the projection of Eq.(\ref{8-7}) to the center subspace is given by
\begin{eqnarray*}
\frac{d}{dt} \Pi_0 Z_1
  = T_{10}^\times  \Pi_0 Z_1  
       + \frac{\varepsilon }{2} \langle  Z_1  \,|\, P_0 \rangle  \Pi_0P_0  
           - \frac{K}{2} \langle \overline{Z_1\,|\, P_0} \rangle  \Pi_0 Z_2.
\end{eqnarray*}
By using Eqs.(\ref{8-8}),(\ref{8-10}),(\ref{8-11}) and (\ref{8-19}), we obtain
\begin{eqnarray*}
\frac{d}{dt} \frac{K_c}{2} \alpha \mu_0
   &=& \frac{K_c}{2} \alpha T^\times_{10} \mu_0 + \frac{\varepsilon }{2} 
                     \left( \alpha  + \langle Y_1 \,|\, P_0 \rangle \right) \frac{1}{D_0} \mu_0 \\
 & & - \frac{K}{2}( \overline{\alpha} + \langle \overline{Y_1 \,|\, P_0} \rangle) \cdot
            \left( -\frac{\pi g''(0)KK_c^2}{16D_0} \alpha ^2 \mu_0 
                   + O(\alpha ^3, \alpha ^2\varepsilon , \alpha \varepsilon ^2, \varepsilon ^3) \right), \\
&=& \frac{\varepsilon }{2D_0} \alpha \mu_0
     + \frac{\pi g''(0) K_c^4}{32 D_0} \alpha |\alpha |^2 \mu_0
          + O(\varepsilon \alpha ^2, \varepsilon ^2\alpha , \varepsilon ^3, \alpha ^4),
\end{eqnarray*}
which yields the dynamics on the center manifold as
\begin{equation}
\frac{d}{dt}\alpha  = \frac{\alpha }{D_0 K_c} \left( 
   \varepsilon  + \frac{\pi g''(0)K_c^4}{16} |\alpha |^2 \right) 
         + O(\varepsilon \alpha ^2, \varepsilon ^2\alpha , \varepsilon ^3, \alpha ^4).
\label{8-20}
\end{equation}
Since $g''(0) <0$, this equation has a fixed point expressed as Eq.(\ref{8-4})
when $\varepsilon  = K- K_c >0$. Note that the order parameter $\eta (t) = (Z_1, P_0)$ is rewritten as
\begin{equation}
\eta (t) = (Z_1, P_0)= \langle Z_1  \,|\, P_0 \rangle
 = \frac{K_c}{2} \alpha \langle \mu_0   \,|\, P_0 \rangle + \langle Y_1  \,|\, P_0 \rangle
 = \alpha + O(\alpha ^2, \alpha \varepsilon , \varepsilon ^2).
\label{8-21}
\end{equation}
Thus the dynamics of the order parameter is also given by Eq.(\ref{8-20}).
To prove that the fixed point (\ref{8-4}) is asymptotically stable, it is sufficient to show the following.
\\[0.2cm]
\textbf{Lemma 7.12.} \, $D_0 >0$.
\\[0.2cm]
\textbf{Proof} \, Put 
\begin{eqnarray*}
f(\lambda ) = 1 - \frac{K_c}{2} \int_{\mathbf{R}} \! \frac{1}{\lambda - \sqrt{-1}\omega } g(\omega )d\omega 
     -\pi K_c g(-\sqrt{-1}\lambda ). 
\end{eqnarray*}
By the definition of $D_0$,
\begin{eqnarray*}
D_0 = f'(0) &=& \lim_{\lambda \to 0} \frac{K_c}{2} \int_{\mathbf{R}} \! \frac{1}{(\lambda - \sqrt{-1}\omega )^2} g(\omega )d\omega 
          + \sqrt{-1}\pi K_c g'(0) \\
&=& \frac{\sqrt{-1}K_c}{2} \lim_{\lambda \to 0} \int_{\mathbf{R}} \! \frac{1}{\lambda - \sqrt{-1}\omega } g'(\omega )d\omega 
          + \sqrt{-1}\pi K_c g'(0). 
\end{eqnarray*}
Since $g(\omega )$ is even,
\begin{eqnarray*}
D_0 = -\frac{K_c}{2} \lim_{x\to 0} \int_{\mathbf{R}} \! \frac{\omega }{x^2 + \omega ^2} g'(\omega ) d\omega 
 = -K_c \lim_{x\to 0} \int^\infty_{0} \! \frac{\omega }{x^2 + \omega ^2} g'(\omega ) d\omega.
\end{eqnarray*}
Since $g(\omega )$ is unimodal, $g'(\omega ) \leq 0$ when $\omega >0$, which proves that $D_0 >0$.
\hfill $\blacksquare$
\\[-0.2cm]

Since $D_0>0, K_c>0, g''(0)<0$, the fixed point $\alpha =0$ (de-synchronous state) is unstable and
the fixed point Eq.(\ref{8-4}) (synchronous state) is asymptotically stable when $\varepsilon  = K- K_c >0$.
This completes the proof.
 \hfill $\blacksquare$


\vspace*{0.5cm}
\textbf{Acknowledgements.}

This work was supported by Grant-in-Aid for Young Scientists (B), No.22740069 from MEXT Japan.



\begin{thebibliography}{99}
\setlength{\baselineskip}{0pt}

\bibitem{Ace}
J. A. Acebr\'{o}n, L. L. Bonilla, C. J. P. Vicente, F. Ritort, R. Spigler,
The Kuramoto model: A simple paradigm for synchronization phenomena,
Rev. Mod. Phys., Vol. 77 (2005), pp. 137-185

\bibitem{Ace2}
J. A. Acebr\'{o}n, L. L. Bonilla,
Asymptotic description of transients and synchronized states of globally coupled oscillators,
Phys. D 114 (1998), no. 3-4, 296--314

\bibitem{Ahl}
L. V. Ahlfors,
Complex analysis. An introduction to the theory of analytic functions of one complex variable,
McGraw-Hill Book Co., New York, 1978

\bibitem{Bal}
N. J. Balmforth, R. Sassi,
A shocking display of synchrony, 
Phys. D 143 (2000), no. 1-4, 21--55

\bibitem{Bat}
P. W. Bates, C. K. R. T Jones,
Invariant manifolds for semilinear partial differential equations,
Dynamics reported, Vol. 2, 1--38, 1989

\bibitem{Ber}
W. Bergweiler, A. Eremenko,
On the singularities of the inverse to a meromorphic function of finite order,
Rev. Mat. Iberoamericana 11 (1995), no. 2, 355--373

\bibitem{Bon1}
L. L. Bonilla, John C. Neu, R. Spigler,
Nonlinear stability of incoherence and collective synchronization in a population of coupled oscillators,
J. Statist. Phys. 67 (1992), no. 1-2, 313--330

\bibitem{Bon2}
L. L. Bonilla, C. J. P\'{e}rez Vicente, R. Spigler,
Time-periodic phases in populations of nonlinearly coupled oscillators with bimodal frequency distributions,
Phys. D 113 (1998), no. 1, 79--97

\bibitem{Che}
X-Y. Chen, J. K. Hale, T. Bin,
Invariant foliations for $C^1$ semigroups in Banach spaces,
J. Differential Equations 139 (1997), no. 2, 283--318

\bibitem{ChiNi}
H.Chiba, I.Nishikawa, 
Center manifold reduction for a large population of globally coupled phase oscillators,
Chaos, 21, 043103 (2011)

\bibitem{Chi4}
H. Chiba,
Continuous limit and the moments system for the globally coupled phase oscillators,
Discret. Contin. Dyn. S.-A, Vol.33, pp.1891-1903 (2013)

\bibitem{ChiPa}
H. Chiba, D. Paz\'{o},
Stability of an $[N/2]$-dimensional invariant torus in the Kuramoto model at small coupling,
Physica D, Vol.238, 1068-1081 (2009)

\bibitem{Cra1}
J. D. Crawford,
Amplitude expansions for instabilities in populations of globally-coupled oscillators,
J. Statist. Phys. 74 (1994), no. 5-6, 1047--1084

\bibitem{Cra2}
J. D. Crawford,
Scaling and Singularities in the Entrainment of Globally Coupled Oscillators,
Phys. Rev. Lett. 74, 4341 (1995)

\bibitem{Cra3}
J. D. Crawford, K. T. R. Davies,
Synchronization of globally coupled phase oscillators: singularities and scaling for general couplings,
Phys. D 125 (1999), no. 1-2, 1--46

\bibitem{Dai}
H. Daido,
Onset of cooperative entrainment in limit-cycle oscillators with uniform all-to-all interactions:
bifurcation of the order function,
Phys. D 91 (1996), no. 1-2, 24--66

\bibitem{Ere}
A. E. Eremenko,
The set of asymptotic values of a finite order meromorphic function,
(Russian) Mat. Zametki 24 (1978), no. 6, 779--783

\bibitem{Gel0}
I. M. Gelfand, G. E. Shilov,
Generalized functions. Vol. 2. Spaces of fundamental and generalized functions,
Academic Press, New York-London, 1968

\bibitem{Gel1}
I. M. Gelfand, N. Ya. Vilenkin,
Generalized functions. Vol. 4. Applications of harmonic analysis,
Academic Press, New York-London, 1964

\bibitem{Gross}
W. Gross,
Eine ganze Funktion, fur die jede komplexe Zahl Konvergenzwert ist,
Math. Ann. 79 (1918), no. 1-2, 201--208

\bibitem{Gro}
A. Grothendieck,
Topological vector spaces,
Gordon and Breach Science Publishers, New York-London-Paris, 1973

\bibitem{Hil}
E. Hille, R. S. Phillips,
Functional analysis and semigroups,
American Mathematical Society, 1957 

\bibitem{Kato}
T. Kato,
Perturbation theory for linear operators,
Springer-Verlag, Berlin, 1995

\bibitem{Kom}
H. Komatsu,
Projective and injective limits of weakly compact sequences of locally convex spaces,
J. Math. Soc. Japan, 19, 1967 366--383

\bibitem{Kri}
T. Krisztin,
Invariance and noninvariance of center manifolds of time-$t$ maps with respect to the semiflow,
SIAM J. Math. Anal. 36 (2004/05), no. 3, 717--739

\bibitem{Kura1}
Y. Kuramoto,
Self-entrainment of a population of coupled non-linear oscillators,
International Symposium on Mathematical Problems in Theoretical Physics,
pp. 420--422. Lecture Notes in Phys., 39. Springer, Berlin, 1975

\bibitem{Kura2}
Y. Kuramoto,
Chemical oscillations, waves, and turbulence,
Springer Series in Synergetics, 19. Springer-Verlag, Berlin, 1984

\bibitem{Mai1}
Y.~Maistrenko, O.~Popovych, O.~Burylko, P.~A. Tass, Mechanism of
  desynchronization in the finite-dimensional Kuramoto model, Phys. Rev.
  Lett. 93 (2004) 084102

\bibitem{Mai2}
Y.~L. Maistrenko, O.~V. Popovych, P.~A. Tass, Chaotic attractor in the
  Kuramoto model, Int. J. of Bif. and Chaos 15 (2005) 3457--3466

\bibitem{Mau}
K. Maurin,
General eigenfunction expansions and unitary representations of topological groups,
Polish Scientific Publishers, Warsaw, 1968

\bibitem{Mar}
E. A. Martens, E. Barreto, S. H. Strogatz, E. Ott ,P. So, T. M. Antonsen,
Exact results for the Kuramoto model with a bimodal frequency distribution,
Phys. Rev. E 79, 026204 (2009)

\bibitem{Marv}
S. A. Marvel, R. E. Mirollo, S. H. Strogatz,
Identical phase oscillators with global sinusoidal coupling evolve by Mobius group action,
Chaos 19, 043104 (2009)

\bibitem{Mir0}
R. E. Mirollo, S. H. Strogatz,
Amplitude death in an array of limit-cycle oscillators,
J. Statist. Phys. 60 (1990), no. 1-2, 245--262

\bibitem{Mir2}
R. Mirollo, S. H. Strogatz,
The spectrum of the partially locked state for the Kuramoto model,
J. Nonlinear Sci. 17 (2007), no. 4, 309--347

\bibitem{New}
D. J. Newman,
A simple proof of Wiener's $1/f$ theorem. 
Proc. Amer. Math. Soc. 48 (1975), 264--265

\bibitem{Ott1}
E. Ott, T. M. Antonsen,
Low dimensional behavior of large systems of globally coupled oscillators,
Chaos 18 (2008), no. 3, 037113

\bibitem{Ott2}
E. Ott, T. M. Antonsen,
Long time evolution of phase oscillator systems,
Chaos 19 (2009), no. 2, 023117

\bibitem{Pik}
A. Pikovsky, M. Rosenblum, J. Kurths, Synchronization: A Universal Concept
  in Nonlinear Sciences, Cambridge University Press, Cambridge, 2001

\bibitem{Reed}
M. Reed, B. Simon,
Methods of modern mathematical physics IV. Analysis of operators,
Academic Press, New York-London, 1978

\bibitem{San}
J. A. Sanders, F. Verhulst,
Averaging methods in nonlinear dynamical systems,
Springer-Verlag, New York, 1985

\bibitem{Sho}
J. A. Shohat, J. D. Tamarkin,
The Problem of Moments,
American Mathematical Society, New York, 1943

\bibitem{Ste}
E. M. Stein, G. Weiss,
Introduction to Fourier analysis on Euclidean spaces,
Princeton University Press, Princeton, 1971

\bibitem{Str1}
S. H. Strogatz,
From Kuramoto to Crawford: exploring the onset of synchronization in populations of coupled oscillators,
Phys. D 143 (2000), no. 1-4, 1--20

\bibitem{Str3}
S. H. Strogatz, R. E. Mirollo, 
Stability of incoherence in a population of coupled oscillators,
J. Statist. Phys. 63 (1991), no. 3-4, 613--635

\bibitem{Str2}
S. H. Strogatz, R. E. Mirollo, P. C. Matthews,
Coupled nonlinear oscillators below the synchronization threshold:
relaxation by generalized Landau damping,
Phys. Rev. Lett. 68 (1992), no. 18, 2730--2733

\bibitem{Tit}
E. C. Titchmarsh, 
Introduction to the theory of Fourier integrals,
Chelsea Publishing Co., New York, 1986

\bibitem{Tre}
F. Tr\'{e}ves,
Topological vector spaces, distributions and kernels,
Academic Press, New York-London, 1967

\bibitem{Van}
A. Vanderbauwhede, G. Iooss,
Center manifold theory in infinite dimensions,
Dynam. Report. Expositions Dynam. Systems, 1, Springer, Berlin, 1992

\bibitem{Vil}
N. J. Vilenkin,
Special functions and the theory of group representations,
American Mathematical Society, 1968

\bibitem{Yos}
K. Yosida,
Functional analysis,
Springer-Verlag, Berlin, 1995

\bibitem{Zha}
G. H. Zhang,
Theory of entire and meromorphic functions,
Deficient and asymptotic values and singular directions,
American Mathematical Society, 1993

\end{thebibliography}
\end{document}